\documentclass[11pt]{amsart}
\usepackage{graphicx}
\usepackage{amssymb}
\usepackage{epstopdf}
\usepackage{amsfonts}
\usepackage{amsmath}
\usepackage{esint,color}
\usepackage{amsthm}
\usepackage{enumerate}
\usepackage{mathtools}
\usepackage{esvect}

\newcommand{\R}{{\mathbb R}}

\newcommand{\Z}{{\mathbb Z}}

\newcommand{\cS}{{\mathcal S}}

\newcommand{\cG}{{\mathcal G}}

\newcommand{\cP}{{\mathcal P}}

\newcommand{\cF}{{\mathcal F}}

\newcommand{\e}{\varepsilon}
\newcommand{\vp}{\varphi}
\newcommand{\osc}{\operatornamewithlimits{osc}}

\newcommand{\spt}{\operatorname{spt}}

\newcommand{\dist}{\operatorname{dist}}

\newcommand{\diam}{\operatorname{diam}}
\newcommand{\rad}{\operatorname{rad}}

\theoremstyle{plain}
\newtheorem{theorem}{Theorem}[section]

\newtheorem{lemma}[theorem]{Lemma}
\newtheorem{proposition}[theorem]{Proposition}

\theoremstyle{definition}
\newtheorem{definition}[theorem]{Definition}

\theoremstyle{remark}

\newtheorem{case}{Case}

\newtheorem{remark}[theorem]{Remark}

\numberwithin{equation}{subsection}

\title[]{Uniform Integrability in Periodic Homogenization of Fully Nonlinear Equations}

\author{Sunghan Kim}
\address{Department of Mathematics, KTH Royal Institute of Technology, 100 44 Stockholm, Sweden}
\email{sunghan@kth.se}
\thanks{The author was supported by Knut and Alice Wallenberg Foundation.}

\begin{document}
    
\begin{abstract}
This paper is devoted to the study of uniform $W^{1,\frac{np}{n-p}}$- and $W^{2,p}$-estimates for viscosity solutions to fully nonlinear, uniformly elliptic, periodic homogenization problems, up to boundaries, subject to Dirichlet boundary conditions. We characterize the size of ``effective'' Hessian and gradient of viscosity solutions to homogenization problems, and prove its uniform integrability without any regularity assumption on the governing functionals. Our estimates are new even for the standard problems. Our analysis applies to a large class of non-convex functionals. 
\end{abstract}

\maketitle

\tableofcontents


\section{Introduction}\label{section:intro}

This paper is devoted to the study of uniform integrability of the Hessian and gradient of viscosity solutions $u^\e\in C(\overline\Omega)$ to fully nonlinear periodic homogenization problems, of the type 
\begin{equation}\label{eq:main}
\begin{dcases}
F \left( D^2 u^\e,\frac{\cdot}{\e}\right) =f  & \text{in }\Omega,\\
u^\e = g & \text{on }\partial \Omega.
\end{dcases}
\end{equation} 

In \cite{Caff}, the size of the ``Hessian'' of a continuous function at a point is characterized by the smallest opening of touching convex and concave paraboloids at that point. Here we extend this concept by allowing the touching to take place in a neighborhood of size $\e$ around the reference point. In addition, we characterize the size of the ``gradient'' in a similar way, by replacing paraboloids with cones. We denote by $H_\Omega^\e$ and $G_\Omega^\e$ the resulting functional for the Hessian and respectively gradient for continuous functions on $\Omega$. See Definition \ref{definition:large-grad-Hess} for more detail.

The first main result for this paper is the uniform integrability of $H_\Omega^\e (u^\e)$; see Definition \ref{definition:W2p} for domains of $W^{2,p}$-type.

\begin{theorem}[$W^{2,p}$-estimates]\label{theorem:W2p}
Let $F\in C(\cS^n\times\R^n)$ be a functional satisfying \eqref{eq:F-ellip} -- \eqref{eq:Fb-W2VMO}, $\Omega\subset\R^n$ be a bounded domain, $f\in C(\Omega)\cap L^p(\Omega)$ for some finite $p> p_0$, $g\in C(\partial\Omega)\cap W^{2,p}(\Omega)$ and $u^\e\in C(\overline\Omega)$ be a viscosity solution to \eqref{eq:main}. Suppose either of the following: (i) $\Omega$ is a $W^{2,n}$-type domain, and $p_0 < p < n$, (ii) $\Omega$ is a $W^{2,n+\sigma}$-type domain for some $\sigma>0$, and $p = n$, (iii) $ \Omega$ is a $W^{2,p}$-type domain and $p > n$, all with size $(\delta,R)$. Then $H_\Omega^\e (u^\e)\in L^p(\Omega_\e)$, with $\Omega_\e = \{x\in\Omega: \dist(x,\partial\Omega)>\e\}$,  and 
$$
\left( \int_{\Omega_\e} ( H_\Omega^\e (u^\e)  )^p dx  \right)^{\frac{1}{p}}\leq C\bigg( \| u^\e \|_{L^\infty(\Omega)} + \| f + |D^2 g| \|_{L^p(\Omega)} \bigg),
$$
where $C>0$ depends only on $n$, $\lambda$, $\Lambda$, $\psi$, $\kappa$, $\sigma$, $\delta$, $R$ and $p$. Assume further that \eqref{eq:F-C} and \eqref{eq:F-W2q} for some $q > p$. Then the same assertion holds for $H_\Omega(u^\e)$, hence $|D^2 u^\e|$, and $\Omega$, in place of $H_\Omega^\e( u^\e)$ and respectively $\Omega_\e$. In this case the constant $C$ may depend also on $q$.
\end{theorem}

Let us provide some motivation for the assumption \eqref{eq:Fb-W2VMO}. Roughly speaking, the assumption says that the effective problem $\bar F(D^2 v ) = 0$ has interior $VMO$-estimates for the Hessian of its (viscosity) solutions. Note that $u^\e$ converges to its effective profile $\bar u$, only, uniformly. This is too weak to ensure any closeness between their Hessian. Under $VMO$-condition on $D^2 \bar u$, however, $D^2 \bar u$ satisfies small $BMO$-condition at an intermediate scale, and we observe that $\cP_\pm (D^2 ( u^\e - \bar u - \e^2 w(\frac{\cdot}{\e}))) = o(1)$, at that scale, with $w$ being an interior corrector. This is a key observation for an approximation lemma (Lemma \ref{lemma:apprx-msr-bdry}) for interior $W^{2,p}$-estimates. 

It should also be addressed that mainly due to \cite[Theorem 3.4]{H}, there is a large class of non-convex functionals satisfying \eqref{eq:Fb-W2VMO}; more specifically, if there exists a functional $F_* :\cS^n\times\R^n\to\R$, which is convex in the first argument and satisfies \eqref{eq:F-ellip} -- \eqref{eq:F-0}, such that $|(F - F_*)(P,y) - (F - F_*)(Q,y)| \leq \theta | P- Q|$ for all $P,Q\in\cS^n$ and all $y\in\R^n$, for some small constant $\theta$,  then the effective functional $\bar F$ satisfies \eqref{eq:Fb-W2VMO}. In particular, assuming \eqref{eq:Fb-W2VMO} is strictly relaxed than assuming interior $C^{2,\alpha}$-estimates, unless the governing functional is continuously differentiable \cite{H2}. 

Moreover, the uniform $L^p$-estimate for $H_\Omega^\e (u^\e)$ above does not require any regularity for $F$ other than mere continuity. The continuity assumption (on $F$ as well as $f$) is imposed only to have homogenization. It is worth mentioning that homogenization of viscosity solutions has not yet been justified for equations with measurable ingredients, and the justification does not seem to be a trivial matter either. The uniform estimate for the ``full'' Hessian, $H_\Omega(u^\e)$, is obtained under the additional hypotheses \eqref{eq:F-C} and \eqref{eq:F-W2q}, which are typical assumptions in the framework of standard problems. 

The above estimates are sharp not only in terms of the regularity of the data, but also of the regularity of the boundary layer. The major challenge here arises from the fact that boundary flattening maps destroy the pattern of the rapid oscillation. For this reason, our analysis here is quite different from, and in fact more complicated than, the argument for standard problems, c.f. \cite{W}. 

As a matter of fact, the boundary estimates for the case $p_0 < p< n$ are even new in the context of standard problems. The analysis is based on the following sharp $W^{1,\frac{np}{n-p}}$-estimates up to the boundary, with $\frac{np}{n-p}$ being the critical Sobolev exponent.

\begin{theorem}[$W^{1,\frac{np}{n-p}}$-estimates]\label{theorem:W1p}
Let $F\in C(\cS^n\times\R^n)$ be a functional satisfying \eqref{eq:F-ellip} -- \eqref{eq:F-0}, $\Omega\subset\R^n$ be a $(\delta,R)$-Reifenberg flat, bounded domain, $f\in C(\Omega)\cap L^p(\Omega)$ for some $p_0 < p< n$, $g\in C(\partial\Omega)\cap W^{1,\frac{np}{n-p}}(\Omega)$ and $u^\e\in C(\overline\Omega)$ be a viscosity solution to \eqref{eq:main}. Then $G_\Omega^\e (u^\e) \in L^{\frac{np}{n-p}}(\Omega_\e)$, and 
$$
\left(\int_{\Omega_\e} ( G_\Omega^\e(u^\e))^{\frac{np}{n-p}} dx\right)^{\frac{1}{p} - \frac{1}{n}} \leq C \bigg( \| u^\e \|_{L^\infty(\Omega)} + \| f + |Dg |\|_{L^{\frac{np}{n-p}}(\Omega)} \bigg),
$$
where $C$ depends only on $n$, $\lambda$, $\Lambda$, $\delta$, $R$ and $p$. Assume further that \eqref{eq:F-C} holds. Then the same assertion holds for $G_\Omega(u^\e)$, hence $|Du^\e|$, and $\Omega$, in place of $G_\Omega^\e(u^\e)$ and respectively $\Omega_\e$. 
\end{theorem}

To the best of author's knowledge, the closest result under the framework of standard problems is \cite[Theorem 1.4]{DKM}, where interior gradient estimates are established in the Lorentz space that $f\in L^{p,\gamma}$ implies $|Du| \in L_{loc}^{\frac{np}{n-p},\gamma}$ for any $p\in(p_0,n)$ and any $\gamma>0$. As for estimates for subcritical Sobolev exponent, the interior and boundary estimates (for $C^2$-domains and $C^{1,\alpha}$-data on boundaries) are obtained in \cite{S} and respectively \cite{W}. 

Our analysis is based on the decay estimate of the set of large ``gradient'', in the spirit of Caffarelli's approach in \cite{Caff}. We present a parallel study for the set of large gradient by replacing the touching paraboloids of the Hessian with cones. The proof relies on the general maximum principle as well as an elementary observation that the slope of supporting hyperplanes for convex envelopes to viscosity solutions in the Pucci class can be universally bounded from below. 

It is worthwhile to mention that the above estimates for the uniform integrability of the gradient are  sharp in terms of the data, and that the domains are only required to be Reifenberg flat. On the other note, as a byproduct via the Sobolev embedding theorem, we obtain a uniform interior $C^{0,2-\frac{n}{p}}$-estimate, which is rather well-understood in the setting of linear homogenization problems \cite{AL2} and was recently proved in standard fully nonlinear problems \cite{Sir}. Nevertheless, this is again a new result in the framework of fully nonlinear homogenization problems.

Let us briefly summarize the recent development of uniform estimates in the homogenization theory. Needless to say, the study has gained its interests due to a series of papers by Avellaneda and Lin. In particular, $W^{1,p}$-estimates are established in \cite{AL3} for linear divergence-type equations, based on the study of Green functions. Later  in \cite{CP}, Caffarelli and Peral proved $W^{1,p}$-estimates for nonlinear divergence-type equations, via Calder\'on-Zygmund cube decomposition argument. In \cite{MS}, Melcher and Schweizer proved the estimates via a more direct approach, based on the observation that $\e$-difference quotients solve the same class of equations. We would also like to mention \cite{KM}, where uniform integrability estimates are established for nonlinear systems in divergence form. More recently, Byun and Kang proved, in \cite{BJ}, uniform $W^{1,p}$-estimates for linear divergence-type systems, under small $BMO$-condition on the periodically oscillating operators, up to Reifenberg flat domains. Some sharp ``large-scale'' estimates for linear divergence-type equations, without any regularity assumption on the governing operators can be found in \cite{She}. All the above results are concerned with periodic homogenization of solutions to either interior problems or Dirichlet problems. As for Neumann problems, some important sharp estimates can be found in \cite{KLS}. There is also a large amount of literature concerning uniform pointwise estimates for random homogenization, for which we would like to refer readers to a recent book \cite{AKM} and the references therein. 

Most of the literature is concerned with weak solutions to divergence-type problems. Uniform estimates for viscosity solutions to non-divergence type equations was done only recently in a collaboration \cite{KL} by Lee and the author, where pointwise $C^{1,\alpha}$- and $C^{1,1}$-estimates are proved for a class of non-convex functionals. The uniform integrability estimates established in this paper are new, even for linear equations in non-divergence form. 

The paper is organized as follows. In the next section, we collect the notation, main assumptions and some preliminaries. Section \ref{section:tech} is devoted to several technical tools used in the subsequent analysis, yet of their own independent interests. In Section \ref{section:decay}, we study universal decay estimates for the set of large Hessian and gradient that will play an important role in the subsequent analysis. In Section \ref{section:W1p}, we establish the uniform $W^{1,\frac{np}{n-p}}$-estimates for both interior (Theorem \ref{theorem:int-W1p}) and boundaries (Theorem \ref{theorem:bdry-W1p}). Finally, in Section \ref{section:W2p}, we prove the uniform $W^{2,p}$-estimates, whose proof is again divided in to the case of interior (Theorem \ref{theorem:int-W2p}) and of boundaries (Theorem \ref{theorem:bdry-W2p}).


\section{Preliminaries}\label{section:prelim}

We shall denote by $B_r(x)$ the $n$-dimensional ball centered at $x$ with radius $r$, and by $Q_r(x)$ the $n$-dimensional cube centered at $x$ with side-length $r$. By $S_\delta(\nu)$ we denote the slab centered at the origin with width $\delta$ in direction $\nu$, i.e., $S_\delta(\nu) = \{x\in \R^n: |x\cdot\nu| < \delta\}$. By $H_t(\nu)$ we denote the half-space in direction $\nu$ with the lowest level being $t$, i.e., $H_\delta(\nu) = \{ x\in\R^n: x\cdot \nu > t\}$. Also we shall write $H_0(\nu)$ simply by $H(\nu)$. Moreover, by $\cS^n$ we denote the space of all symmetric $(n\times n)$-matrices. 

Throughout this paper, $\lambda$ and $\Lambda$ will be fixed as some positive constants, with $\lambda\leq\Lambda$, and will also denote the lower and respectively upper ellipticity bound for the governing functional. By $p_0$ we shall denote Escauriaza's constant such that the generalized maximum principle holds for all $p> p_0$; note $\frac{n}{2} < p_0< n$ and it depends only on $n$ and the ellipticity bounds, $\lambda$ and $\Lambda$. For more details, we refer readers to \cite{Esc} and \cite{CCKS}. In addition, we denote by $\cP_-$ and $\cP_+$ the Pucci minimal and respectively maximal functional on $\cS^n$, associated with ellipticity bounds $\lambda$ and $\Lambda$, such that
$$
\cP_-(P) = \lambda \sum_{e_i> 0} e_i + \Lambda \sum_{e_i< 0} e_i, \quad \cP_+(P) = - \cP_-(-P),
$$
where $e_i$ is the $i$-th eigenvalue of $P$. 

Let $\psi:(0,\infty)\to(0,\infty)$ be a nondecreasing, strictly concave function such that $\psi(0+) = 0$, and $\Omega$ be a domain in $\R^n$. We say $g\in BMO_\psi(\Omega)$, if $g\in L^1(\Omega)$ and 
$$
[g]_{BMO_\psi(\Omega)} = \sup_{B\subset\Omega} \frac{1}{|B|\psi(\rad B)} \int_B |g - (g)_B| dx < \infty,
$$
where the supremum is taken over all balls $B\subset\Omega$. Given any $g\in L^1(\R^n)$ with $g\geq 0$, $M(g)$ denotes the maximal function of $g$, i.e., 
$$
M(g) (x) = \sup_{B} \frac{1}{|B|} \int_B g(y)dy, 
$$
where the supremum is taken over all balls $B$ containing $x$. Given any $g\in L_{loc}^1(\R^n)$, $I_\alpha(g)$ denotes the Riesz potential of $g$, i.e., 
$$
I_\alpha(g) (x) = c_\alpha \int_{\R^n} \frac{g(y)}{|x-y|^{n-\alpha}}dy,
$$
with $c_\alpha$ being a suitable normalization constant. 

For definiteness, we shall assume that $F\in C(\cS^n\times\R^n)$ is a functional satisfying, for any $P,Q\in\cS^n$ and any $y\in\R^n$, 
\begin{equation}\label{eq:F-ellip}
\cP_- (P - Q) \leq F(P,y) - F(Q,y) \leq \cP_+(P-Q),
\end{equation}
\begin{equation}\label{eq:F-peri}
F(P, y +k) = F(P,y),
\end{equation}
\begin{equation}\label{eq:F-0}
F(0,y) = 0.
\end{equation}
Under the first two assumptions, there exists a unique functional $\bar F:\cS^n\to\R$ (the so-called effective functional), according to \cite[Theorem 3.1]{E}, such that any limit $\bar u\in C(\Omega)$ of the sequence of viscosity solutions $u^\e\in C(\Omega)$ to $F(D^2 u^\e,\frac{\cdot}{\e}) = f$ in $\Omega$, with $f\in C(\Omega)$ and $\e>0$, under locally uniform convergence as $\e\to 0$ is a viscosity solution to $\bar F(D^2\bar u) = f$ in $\Omega$. We shall suppose that 
\begin{equation}\label{eq:Fb-W2VMO}
\text{$\bar F(D^2 v) = 0$ has interior $W^{2,BMO_\psi}$-estimates with constant $\kappa$.}
\end{equation}
where $\psi:(0,\infty)\to (0,\infty)$ is a concave, non-decreasing function satisfying $\psi(0+) = 0$ and $\kappa>0$ is a fixed constant. More specifically, by \eqref{eq:Fb-W2VMO}, we indicate the following: given any ball $B_R\subset\R^n$ and any function $v_0 \in C(\partial B_R)$, there exists a viscosity solution $v\in C(\overline{B_R})\cap W^{2,1}(B_R)$ to $\bar F(D^2 v ) = 0$ in $B_R$, $v = v_0$ on $\partial B_R$ such that for any $r\in(0,R)$, 
$$\frac{1}{(R-r)^n}\int_{B_r} |D^2 v |\,dx + [D^2 v]_{BMO_\psi(B_r)} \leq \frac{\kappa}{(R-r)^2} \| v_0 \|_{L^\infty(\partial B_R)}.$$
In order to obtain uniform estimates for the ``full'' gradient, we shall assume that $F(P,\cdot) \in BMO_\psi(Q_2)$ (with $Q_1$ being the unit periodic cell) for all $P\in\cS^n$, and 
\begin{equation}\label{eq:F-C}
\sup_{B\subset Q_2} \frac{1}{|B|\psi(\rad B)} \int_B |F(P, y) - (F)_B(P)|^n\,dy \leq \kappa, 
\end{equation}
where the supremum is taken over all balls $B\subset Q_2$ and $(F)_B(P)$ denotes the integral average of $F(P,\cdot)$ over $B$, i.e., $(F)_B(P) = \frac{1}{|B|} \int_B F(P,y)\,dy$. As for the estimates for the full Hessian, we will additionally suppose that 
\begin{equation}\label{eq:F-W2q}
\text{$F(D^2 v , y_0) = 0 $ has interior $W^{2,q}$-estimates with constant $\kappa$,}
\end{equation} 
uniformly for all $y_0\in Q_1$. 

Let us collect by now standard results regarding periodic homogenization for fully nonlinear problems in the following lemma.

\begin{lemma}[{\cite[Lemma 3.1-3.2]{E}}]\label{lemma:hom}
Let $F\in C(\cS^n\times\R^n)$ be a functional satisfying \eqref{eq:F-ellip} and \eqref{eq:F-peri}. Then there exists a unique functional $\bar F :\cS^n\to\R$, satisfying \eqref{eq:F-ellip}, such that for each $P\in\cS^n$, $\bar F(P)$ is the unique constant for which there exists a viscosity solution to 
$$
\begin{dcases}
F(D^2 w + P,\cdot\,) = \bar F(P) &\text{in }\R^n,\\
w(y + k ) = w(y) & \text{for all } y\in\R^n, k\in\Z^n. 
\end{dcases}
$$
In particular, if $w\in C(\R^n)$ is a viscosity solution to this problem, $w\in C^\alpha(\R^n)$ and 
$$
\| w - w(0)\|_{C^\alpha(\R^n)}  + |\bar F(P)| \leq c(|P| + \|F (P,\cdot)\|_{L^\infty(\R^n)}),
$$ 
where $c>0$ and $\alpha\in(0,1)$ depend only on $n$, $\lambda$ and $\Lambda$. 
\end{lemma}

Next, we introduce the set of large ``Hessian'', as well as the set of large ``gradient'', with room for errors of order $\e^2$ and respectively $\e$. These sets will play the main role throughout this paper, as our primary goal is to establish ``large-scale'' estimates. The set of large Hessian without any room for error has played a central role in the $W^{2.p}$-theory for standard fully nonlinear problems. Nevertheless, the set of large gradient seems to appear in this paper for the first time, in the literature. Needless to say, the sets with room for errors are entirely new, as far as the author is concerned. It should be stressed that one can generalize this concept to homogenization problems under various oscillating structures, such as quasi-periodic, almost-periodic or random environment.

\begin{definition}\label{definition:large-grad-Hess}
Let $\Omega\subset\R^n$ be a bounded domain, $\e\geq 0$, $t>0$ and $u\in C(\Omega)$ be given. Let $A_t^\e(u,\Omega)$ (and $L_t^\e(u,\Omega)$) be defined as a subset of $\Omega$ such that $x_0 \in \Omega\setminus A_t^\e(u,\Omega)$ (resp., $\Omega\setminus L_t^\e(u,\Omega)$) if and only if  there exists a linear polynomial $\ell$ (resp., a constant $a$) for which  $| u(x)-\ell(x)| \leq \frac{t}{2} ( |x-x_0|^2 + \e^2)$ (resp., $| u(x)-  a | \leq t (|x - x_0| + \e)$) for all $x\in\Omega$. Denote by $A_t(u,\Omega)$ (and $L_t(u,\Omega)$) the set $A_t^0(u,\Omega)$ (resp., $L_t^0(u,\Omega)$). Let $H_\Omega^\e(u) \, (\text{resp., }G_\Omega^\e(u)) \,: \Omega\to\R\cup\{\pm\infty\}$  be defined as 
$$
\begin{aligned}
&H_\Omega^\e(u) (x) = \inf\{t > 0 : x\in A_t^\e(u,\Omega)\},\\
&(\text{resp., } G_\Omega^\e(u)(x) = \inf\{t>0 : x\in L_t^\e(u,\Omega)\} ),
\end{aligned}
$$
and denote by $H_\Omega(u)$ (and $G_\Omega(u)$) the function $H_\Omega^0(u)$ (resp., $G_\Omega^0(u)$).
\end{definition}

The following is the definition for the Reifenberg flat sets, which will appear in uniform boundary $W^{1,\frac{np}{n-p}}$-estimates. 

\begin{definition}\label{definition:reifenberg}
Let $\Omega$ be a domain, and $U$ be a neighborhood of a point at $\partial\Omega$. Set $\partial\Omega\cap U$ is said to be $(\delta,R)$-Reifenberg flat from exterior (or interior), if for any $x_0\in\partial\Omega\cap U$ and any $r \in(0,R]$, $\Omega\cap B_r(x_0)\subset \{ x: (x- x_0)\cdot \nu_{x_0} > -  \delta r\}$ (resp., $B_r(x_0)\setminus \Omega \subset \{ x: (x- x_0)\cdot \nu_{x_0} <  -  \delta r\}$) for some unit vector $\nu_{x_0,r}$; here $\nu_{x_0,r}$ may vary upon both $x_0$ and $r$. The set $\partial\Omega\cap U$ is said to be $(\delta,R)$-Reifenberg flat, if it is $(\delta,R)$-Reifenberg flat from both exterior and interior. The domain $\Omega$ is said to be $(\delta,R)$-Reifenberg flat, if $\partial\Omega\cap B_R(x_0)$ is $(\delta,R)$-Reifenberg flat for each $x_0\in\partial\Omega$. 
\end{definition} 

Next, we define domains of $W^{2,p}$-type.

\begin{definition}\label{definition:W2p}
Let $p > 1$ be a constant, $\Omega$ be a domain, and $U$ be a neighborhood of a point at $\partial\Omega$. Set $\partial\Omega\cap U$ is said to be of $W^{2,p}$-type with size $\kappa$, if there exists a neighborhood $V\subset\R^n$ and a diffeomorphism $\Phi\in C^1(U;V)\cap W^{2,p}(U;V)$ such that $\Phi(\Omega\cap U) = H(e_n)\cap V$, $\Phi(\partial\Omega\cap U) = H(e_n)\cap V$, $\osc_U D\Phi \leq \delta$ and $\| D^2 \Phi\|_{L^p(U)} \leq \delta$. The domain $\Omega$ is said to be of $W^{2,p}$-type with size $(\kappa,R)$, if $\partial\Omega\cap B_R(x_0)$ is of $W^{2,p}$-type with size $\kappa$ for each $x_0\in\partial\Omega$.
\end{definition} 

We shall also need some covering lemmas. As for the analysis for interior estimates, we shall use the classical Calder\'on-Zygmund cube decomposition lemma, c.f. \cite[Section 9.2]{GT} and \cite[Lemma 4.1]{CC}:

\begin{lemma}[Calder\'on-Zygmund cube decomposition]\label{lemma:cal-zyg}
Let $A\subset Q_1$ be a measurable set such that $|A|\leq \eta$ for some $\eta\in (0,1)$. Then there exists a finite collection $\cF$ of cubes from the dyadic subdivision of $Q_1$, such that $|A\cap Q| > \eta |Q|$ for all $Q\in\cF$, and $|A\cap \tilde Q|\leq \eta |\tilde Q|$ for the predecessor of $Q$.

Let $B$ be a measurable set such that $A\subset B\subset Q_1$. If $\tilde Q\subset B$ for the predecessor $\tilde Q$ of any dyadic cube $Q$ satisfying $|Q\cap A| >\eta |Q|$, then $|A| \leq \eta |B|$. 
\end{lemma}

As for the boundary estimates, we shall resort to a Vitali-type covering lemma, for Reifenberg flat domains, namely \cite[Theorem 2.8]{BW}. We present its statement for the reader's convenience.

\begin{lemma}[Vitali-type covering {\cite[Theorem 2.8]{BW}}]\label{lemma:vitali}
Let $\Omega\subset\R^n$ be a domain such that $\partial\Omega\cap B_1$ is $(\delta,1)$-Reifenberg flat, with some $\delta\in(0,\frac{1}{8})$, and contains the origin. Let $D\subset E\subset \Omega\cap B_1$ be two measurable sets. Suppose that $|D|\leq \eta |B_1|$, for some $\eta\in(0,1)$, and that for any ball $B\subset B_1$ whose center lies in $\overline\Omega\cap B_1$ and radius is at most $1$, $|D\cap B| > \eta |B|$ implies that $\Omega\cap B\subset E$. Then $|D| \leq (10/(1-\delta))^n \eta |E|$. 
\end{lemma}


\section{Some Technical Tools}\label{section:tech}

Let us begin with an assertion that the difference between a viscosity solution and a viscosity sub- or super-solution belongs to the Pucci class in the viscosity sense. It is particularly important that one of two must be a solution. This assertion might be known for some experts. Still, we intend to present a proof because the assertion is not as simple as it sounds, apart from the fact that the author was not able to find a proof in the literature. It should be stressed that the assertion is yet to be known if we replace viscosity solution with viscosity super- or sub-solution (depending on what it is compared with).

\begin{lemma}\label{lemma:visc}
Let $\Omega\subset\R^n$ be a bounded domain, and $F\in C(\cS^n\times\Omega)$ be a functional satisfying $\cP_-(P - Q)\leq F(P,x) - F(Q,x) \leq \cP_+(P-Q)$ for all $P,Q\in\cS^n$ and $x\in\Omega$. Let $u,v\in C(\Omega)$ be such that $F(D^2 u,\cdot\,) = f$ in $\Omega$ and $F(D^2 v, \cdot )\geq g$ in $\Omega$  in the viscosity sense, for some $f,g\in C(\Omega)$. Then $\cP_+(D^2 (u - v)) \leq f -g$ in $\Omega$ in the viscosity sense. 
\end{lemma}

\begin{proof}
Fix $\delta>0$, and denote $\Omega_\delta =\{x\in\Omega:\dist(x,\partial\Omega)>\delta\}$. Clearly, $u,v\in C(\overline{\Omega_\delta})$, and $\Omega_\delta$ satisfies the uniform exterior sphere condition with radius at most $\delta^{-1}$. Given any pair $(\tau,\sigma)$ of real parameters such that $0<\sigma<\tau<\delta$, let $v_\tau:\Omega\to\R$ be the sup-inf convolution of $v$ over $\Omega$; i.e.,
$$
v_\tau (x) = \inf_{z\in\Omega} \sup_{y\in\Omega} \left( v(y) - \frac{|z-y|^2}{2\tau} + \frac{|x-z|^2}{2\sigma} \right).
$$
Such a regularization is by now considered standard. Among other important properties, we shall use the following, which can be found in \cite[Theorem]{LL} that $v_\tau \to v$ uniformly in $\overline{\Omega_\delta}$, $v_\tau \in W^{2,\infty}(\Omega)$, $|D v_\tau| \leq c_\delta \tau^{-1/2}$ in $\Omega$, with $c_\delta>0$ being a constant depending only on the sup-norm and the modulus of continuity of $v$ over $\overline{\Omega_\delta}$, and
$$
F_\tau (D^2 v_\tau, \cdot\,) \geq g_\tau\quad\text{a.e.\ in }\Omega_\tau,
$$
where $F_\tau:\cS^n\times\Omega_\tau\to \R$ and $g_\tau : \Omega_\tau\to\R$ are defined by $F_\tau(P,x) = F(P, x - (\tau-\sigma) Dv_\tau (x))$ and $g_\tau (x) = g( x- (\tau-\sigma) Dv_\tau(x))$. In particular, since $\tau|Dv_\tau| \to 0$ uniformly on $\Omega$, $F_\tau \to F$ locally uniformly in $\cS^n\times\overline{\Omega_\delta}$ and $g_\tau \to g$ uniformly on $\overline{\Omega_\delta}$. 

Consider an auxiliary Dirichlet boundary value problem,
$$
\begin{dcases}
F_\tau (D^2 u_\tau, \cdot\,) = f & \text{in }\Omega_\delta\\
u_\tau = u & \text{on }\partial\Omega_\delta.
\end{dcases}
$$
Since $F_\tau\in C(\cS^n\times\Omega_\delta)$ is uniformly elliptic, $f\in C(\Omega_\delta)$, $u\in C(\partial\Omega_\delta)$ and $\Omega_\delta$ satisfies the uniform exterior sphere condition, there exists a unique viscosity solution $u_\tau\in C(\overline{\Omega_\delta})$ to this problem, according to e.g., \cite[Theorem 4.1]{CIL}. Moreover, as the radius for the uniform exterior sphere condition for $\Omega_\delta$ being independent of $\tau$,  it follows from the global regularity of viscosity solutions \cite[Proposition 4.14]{CC} that $\{ u_\tau : 0<\tau<\delta\}$ is a uniformly bounded and equicontinuous family on $\overline{\Omega_\delta}$. However, since $F_\tau \to F$ locally uniformly on $\cS^n\times\Omega_\delta$, by the stability \cite[Proposition 4.11]{CC} and the comparison principle \cite[Theorem 3.3]{CIL} for viscosity solutions, one can easily deduce that $u_\tau \to u$ uniformly on $\overline{\Omega_\delta}$ as $\tau \to 0$. 

On the other hand, as $v_\tau \in W^{2,\infty}(\Omega)$ and $F_\tau (D^2 v_\tau, \cdot\,) \geq g_\tau$ a.e.\ in $\Omega_\delta$, we can compute that
$$
\cP_- (D^2 (u_\tau - v_\tau)) \leq F_\tau (D^2 u_\tau, \cdot\,) - F_\tau (D^2 v_\tau, \cdot\,) \leq f - g_\tau\quad\text{in }\Omega_\delta,
$$
in the $L^\infty$-viscosity sense, but then in the usual ($C$-)viscosity sense as $f - g_\tau\in C(\Omega_\delta)$; 
here, we refer to \cite{CCKS} the notion of $L^\infty$-viscosity solutions. Now letting $\tau \to 0$, and recalling that $u_\tau \to u$, $v_\tau \to v$ and $g_\tau \to g$ uniformly on $\overline{\Omega_\delta}$, we may conclude from the stability theory again that 
$$
\cP_- (D^2 (u - v)) \leq f - g\quad\text{in }\Omega_\delta. 
$$
As $\delta>0$ was an arbitrary constant, the assertion of the lemma follows by sending $\delta\to 0$. 
\end{proof}

Let us close this section with a few results that are essentially due to \cite{KL}, but extended so as to be adoptable in our subsequent analysis. We shall start with an interior $L^\infty$-approximation of viscosity solutions to periodic homogenization problems by those to the corresponding effective problems. The assertion is a slight generalization of \cite[Lemma 3.1, 4.2]{KL}, which was established for bounded data. Here we shall extend the result to $L^p$-integrable data. 

\begin{lemma}[Due to {\cite[Lemma 3.1, 4.2]{KL}}]\label{lemma:apprx}
Let $F\in C(\cS^n\times\R^n)$ be a functional satisfying \eqref{eq:F-ellip} -- \eqref{eq:F-0}, $f\in C(B_R)\cap L^p(B_R)$ for some $p > p_0$ and some $R>0$, $u^\e \in C(B_R)$ be a viscosity solution to 
$$
F \left( D^2 u^\e ,\frac{\cdot}{\e} \right) = f \quad\text{in }B_R,
$$
for some $\e>0$. Let $r\in(0,R)$ be given. Then for each $\eta>0$, one can find some $\e_\eta> 0$, depending only on $n$, $\lambda$, $\Lambda$ and $\eta$, such that if $0< \e < r\e_\eta$ and $0<r<R$, then there exists a viscosity solution $\bar u\in C(\overline{B_r})$ to 
$$
\bar F(D^2 \bar u) = 0\quad\text{in }B_r,
$$
for which 
$$
\| \bar u\|_{L^\infty(B_r)}+ \frac{ \| u^\e - \bar u \|_{L^\infty(B_r)}}{\eta} \leq \frac{C }{(R-r)^\alpha} \left( \| u^\e \|_{L^\infty(B_R)}+ R^{2-\frac{n}{p}} \|f \|_{L^p(B_R)}\right),
$$
where $\alpha\in(0,1)$ depends only on $n$, $\lambda$ and $\Lambda$, and $C>1$ may depend further on $p$. 
\end{lemma}

\begin{proof}
The assertion for $f = 0$ is a direct consequence of \cite[Lemma 3.1, 4.2]{KL}. As for the general case, we consider an auxiliary boundary value problem,
$$
\begin{dcases}
F\left(D^2 \hat u^\e, \frac{\cdot}{\e}\right) = 0 &\text{in }B_R,\\
\hat u^\e = u^\e &\text{on }\partial B_R,
\end{dcases}
$$
which admits a unique viscosity solution. By the general maximum principle, one may easily construct a barrier function to verify that 
$$
\| u^\e - \hat u^\e \|_{L^\infty(B_R)} \leq c R^{2-\frac{n}{p}} \|f \|_{L^p(B_R)},
$$
for some $c>0$ depending only on $n$, $\lambda$, $\Lambda$ and $p$. This combined with the assertion with $f= 0$ yields the conclusion. 
\end{proof} 

With the above lemma at hand, we can extend the uniform pointwise $C^{1,\alpha}$-estimates for fully nonlinear homogenization problems, established in \cite{KL}, to a more general setting. 

\begin{lemma}[Due to {\cite[Theorem 4.1]{KL}}]\label{lemma:int-C1a}
Let $F\in C(\cS^n\times\R^n)$ be a functional satisfying \eqref{eq:F-ellip} -- \eqref{eq:F-0}, $\Omega\subset\R^n$ be a bounded domain, $f\in C(\Omega)\cap L^p(\Omega)$ for some $p> p_0$, and $u^\e \in C(\Omega)$ be a viscosity solution to 
$$
F\left( D^2 u^\e,\frac{\cdot}{\e}\right) = f\quad\text{in }\Omega,
$$
for some $\e>0$. Suppose that $\bar F(D^2 v) = 0$ admits an interior $C^{1,\bar\alpha}$-estimates with constant $\kappa$, for some $\bar\alpha \in(0,1)$ and $\kappa>0$. Let $\alpha \in (0,\bar\alpha)$ be given. Given any $x_0\in \Omega_\e$ for which $I_{(1-\alpha)p}(|f|^p\chi_\Omega)(x_0) < \infty$, there exists a linear polynomial $\ell_{x_0}^\e$ such that for any $x\in\Omega$, 
$$
\begin{aligned}
& |D\ell_{x_0}^\e| + \sup_{x\in\Omega} \frac{| (u^\e - \ell_{x_0}^\e)(x)| }{|x - x_0|^{1+\alpha} + \e^{1+\alpha}} \\
&\leq C \left( \frac{\| u^\e \|_{L^\infty(\Omega)}}{\dist(x_0,\partial\Omega)^{1+\alpha}} +  ((I_{(1-\alpha)p}(|f|^p\chi_\Omega)(x_0))^{\frac{1}{p}}\right),
\end{aligned}
$$ 
where $C>0$ depends only on $n$, $\lambda$, $\Lambda$, $\kappa$, $\bar\alpha$, $\alpha$, $p$ and $\diam(\Omega)$. Moreover, if $F(D^2 v,y_0) = 0$ has interior $C^{1,\bar\alpha}$-estimate with constant $\kappa$, for any $y_0\in\R^n$, and \eqref{eq:F-C} holds for some modulus of continuity $\psi:(0,\infty)\to(0,\infty)$, then the same inequality holds for all $x_0\in\Omega$, instead of $\Omega_\e$, without the additional term $\e^{1+\alpha}$, and with $C$ depending further on $\psi$. 
\end{lemma}

\begin{proof}
This assertion is proved for the case $f\in L^\infty$ in \cite[Theorem 4.1]{KL}. However, the same proof works equally well for any point $x_0\in\Omega$ featuring  $I_{(1-\alpha)p}(|f|^p\chi_\Omega)(x_0) < \infty$, since this implies that 
$$
\sup_{r>0} r^{(1-\alpha)p-n} \int_{B_r(x_0)\cap\Omega} |f|^p\,dx <\infty.
$$
With the latter observation, the iteration argument in \cite[Lemma 4.3]{KL} works, without any notable modification, once we invoke Lemma \ref{lemma:apprx} as the approximation lemma, in place of \cite[Lemma 4.2]{KL}, in the proof there. Let us remark that the iteration technique for standard problems is by now understood as standard, c.f. \cite[Remark 2.5]{S}. Therefore, we shall not repeat the detail here. 
\end{proof}

The following lemma is a uniform boundary $C^{1,\alpha}$-estimates, which extends \cite[Theorem 5.1]{KL} to $L^p$-integrable datum, and $C^{1,\alpha}$-domains. 

\begin{lemma}[Due to {\cite[Theorem 5.1]{KL}}]\label{lemma:bdry-C1a}
Let $F\in C(\cS^n\times\R^n)$ be a functional satisfying \eqref{eq:F-ellip} -- \eqref{eq:F-0}, $\Omega\subset\R^n$ be a bounded domain, $U\subset\R^n$ be a neighborhood of a point on $\partial\Omega$ such that $\partial\Omega\cap U$ is a $C^{1,\alpha}$-graph, whose norm is bounded by $\kappa$,  for some $\alpha\in(0,1)$ and some $\kappa>0$, $f\in C(\Omega)\cap L^p(\Omega\cap U)$ for some $p> n$, $g\in C^{1,\alpha}(\partial\Omega\cap U)$, and $u^\e \in C(\overline\Omega\cap U)$ be a viscosity solution to 
$$
\begin{dcases}
F\left( D^2 u^\e,\frac{\cdot}{\e}\right) = f & \text{in }\Omega\cap U,\\
u^\e = g &\text{on }\partial \Omega\cap U,
\end{dcases}
$$
for some $\e>0$. Set $\alpha_p = \min\{\alpha,1-\frac{n}{p}\}$. Given any $x_0\in \partial\Omega\cap U_\e$, there exists a linear polynomial $\ell_{x_0}^\e$ such that for any $x\in\Omega\cap U$, 
$$
\begin{aligned}
& |D\ell_{x_0}^\e| + \sup_{x\in\Omega} \frac{| (u^\e - \ell_{x_0}^\e)(x)| }{|x - x_0|^{1+\alpha_p} + \e^{1+\alpha_p}} \\
& \leq C \left( \frac{\| u^\e \|_{L^\infty(\Omega\cap U)}}{\dist(x_0,\partial U)^{1+\alpha_p}} + \|f \|_{L^p(\Omega\cap U)}+ \| g\|_{C^{1,\alpha}(\partial\Omega\cap U)}\right),
\end{aligned}
$$ 
where $C>0$ depends only on $n$, $\lambda$, $\Lambda$, $\kappa$, $\alpha$, $p$ and $\diam(U)$. Moreover, under the same additional hypothesis in Lemma \ref{lemma:int-C1a}, the same inequality holds for all $x_0\in\Omega\cap U$, instead of $\Omega\cap U_\e$, without the additional term $\e^{1+\alpha}$, and with $C$ depending further on $\psi$. 
\end{lemma} 

\begin{remark}\label{remark:bdry-C1a}
The reason that we state the above lemma for the case $p>n$ is only because the other case, i.e., $p_0< p\leq n$, holds in a much general setting, namely for standard problems in the Pucci class, c.f. \cite[Theorem 1.6]{LZ}; of course, in the latter case, one needs to replace $\alpha_p$ with $\alpha$ and $\|f \|_{L^p(\Omega\cap U)}$ with $((I_{(1-\alpha)p}(|f|^p\chi_{\Omega\cap U})(x_0))^{1/p}$, as in the case of Lemma \ref{lemma:int-C1a}. 
\end{remark}

\begin{proof}[Proof of Lemma \ref{lemma:bdry-C1a}]
With a similar modification shown in the proof of Lemma \ref{lemma:apprx}, we can extend $f\in L^\infty(\Omega\cap U)$, required in \cite[Lemma 5.2 and 5.3]{KL}, to $L^p(\Omega\cap U)$. To extend the lemmas to $\partial\Omega\cap U \in C^{1,\alpha}$ (from $C^{1,1}$), we may combine the compactness argument in \cite[Lemma 3.1]{LZ3} with \cite[Lemma 3.1, 5.2]{KL}. We skip the detail. 
\end{proof}


\section{Universal Decay Estimates}\label{section:decay}

This section is devoted to a global, universal decay estimate of the measure of the set for large ``gradient'' and ``Hessian'' of viscosity solutions to fully nonlinear equations, up to Reifenberg flat boundaries. As surprising as it may sound, our estimates would not see the boundary value of solutions, as long as the solutions are bounded up to the boundaries. Roughly speaking, this is because of the fact that the boundary layer, as a Reifenberg flat set, can be trapped in between a thin slab, which already has small measure and thus can be neglected. Of course, at a cost, the decay rate we establish here could be extremely slow, yet universal. Let us remark that such a global estimate is hinted in \cite{W}, which proves global universal decay estimate for Hessians, up to flat boundaries. 


\subsection{Set of Large Gradient}\label{section:W1d}

Let us begin with estimates for the set of large gradient. Throughout this section, given any $u\in C(\Omega)$ and $t>0$, $\underline L_t(u,\Omega)$ is the subset of $\Omega$ such that $x_0\in \Omega\setminus \underline L_t(u,\Omega)$ if and only if there exists a constant $a$ for which $u(x) \geq a - t|x-x_0|$ for all $x\in\Omega$. Clearly, $L_t(u,\Omega) = \underline L_t(u,\Omega)\cap \underline (-u,\Omega)$, with $L_t(u,\Omega)$ as in Definition \ref{definition:large-grad-Hess}.

\begin{proposition}\label{proposition:W1d}
Let $\Omega\subset\R^n$ be a bounded, $(\delta,R)$-Reifenberg flat domain with $R\in(0,1]$, $f\in L^p(\Omega)$ be given, with some $p>p_0$, and $u\in C(\Omega)\cap L^\infty(\Omega)$ be an $L^p$-viscosity solution to $\cP_- (D^2 u) \leq f$ in $\Omega$. Then for any $t>0$,
$$
|\underline L_t (u,\Omega)| \leq \frac{C|\Omega|}{R^\mu t^\mu} (\|u \|_{L^\infty(\Omega)} +\| f\|_{L^p(\Omega)})^\mu,
$$
where $C>1$, $\delta >0$ and $\mu>0$ depend only on $n$, $\lambda$, $\Lambda$ and $p$. 
\end{proposition}

We shall split our analysis into two parts, each concerning interior and respectively boundary layer. Let us begin with the interior case first. As the analysis below will be of local character, we shall confine ourselves to the case $\Omega = B_{4\sqrt n}$ and $\Omega' = Q_1$. The following lemma is the gradient-counterpart of \cite[Lemma 7.5]{CC}.  

\begin{lemma}\label{lemma:apprx-W1d}
Let $\Omega\subset\R^n$ be a bounded domain with $B_{4\sqrt n}\subset\Omega$, $u\in C(\Omega)$ be an $L^p$-viscosity solution to $\cP_-(D^2 u) \leq f$ in $B_{4\sqrt n}$, for some $f \in L^p(\Omega)$ satisfying $\| f \|_{L^p(B_{4\sqrt n})} \leq \delta_0$, with $p>p_0$, such that $\inf_{B_{4\sqrt n}} u \geq -1$, $\inf_{Q_3} u \leq 0$ and $u(x) \geq - |x|$ for all $x\in\Omega\setminus B_{2\sqrt n}$. Then 
$$
|\underline{L}_m (u,\Omega)\cap Q_1| \leq \sigma, 
$$
where $m>1$, $\delta_0>0$ and $\sigma \in(0,1)$ are constants depending at most on $n$, $\lambda$, $\Lambda$ and $p$. 
\end{lemma}

\begin{proof}
As in the proof of \cite[Lemma 7.5]{CC}, we consider an auxiliary function $w = u + 1 +\vp$ on $B_{4\sqrt n}$, where $\vp \in C^2(B_{4\sqrt n})$ is the barrier function found in \cite[Lemma 4.1]{CC}. Due to the general maximum principle (which is available as $f\in L^p$ with $p>p_0$), one can argue as in the proof of \cite[Lemma 4.5]{CC} (where the smallness of $\delta_0 \geq \|f \|_{L^p(B_{4\sqrt n})}$ is determined) and observe that
$$
|\{ w = \Gamma_w\}\cap Q_1| \geq 1-\sigma, 
$$
for some constant $\sigma>0$, depending only on $n$, $\lambda$, $\Lambda$ and $p$, where $\Gamma_w$ is the convex envelope of $-w^-$ in $B_{4\sqrt n}$. Our claim is, as again in the proof of \cite[Lemma 7.5]{CC}, that 
$$
\underline{L}_m(u,\Omega) \subset Q_1\setminus \{w =\Gamma_w\},
$$
for some large constant $m>1$, depending only on $n$, $\lambda$ and $\Lambda$. 

The main observation here is that the gradient of the supporting hyperplanes for the convex envelope $\Gamma_w$ at the contact set is universally bounded. This actually follows from a simple fact that by construction, $\Gamma_w = 0$ on $\partial B_{4\sqrt n}$, $w\geq 0$ in $B_{4\sqrt n }\setminus B_{2\sqrt n}$ and $-m \leq \inf_{Q_3} w \leq -1$, for some $m>0$ depending only on $n$, $\lambda$ and $\Lambda$. These inequalities follows from the specific choice of the barrier function $\vp$ in \cite[Lemma 4.1]{CC}, that $\vp \geq 0$ on $\R^n\setminus B_{2\sqrt n}$ and $\inf_{B_{2\sqrt n}} \vp \geq - m$. 

In what follows, we shall let $m$ denote a generic positive constant depending only on $n$, $\lambda$ and $\Lambda$ and allow it to vary at each occurrence. 

Keeping in mind of these properties of $\vp$, let $x_0\in \{w=\Gamma_w\}\cap Q_1$ be any. As $\Gamma_w$ being the convex envelope of $w$ in $B_{4\sqrt n}$, we can find a linear polynomial $\ell$ (as one of the supporting hyperplanes of $\Gamma_w$ at $x_0$) such that $w\geq \ell$ in $B_{4\sqrt n}$ and $w(x_0) = \ell(x_0)$. Since $\Gamma_w = 0$ on $\partial B_{4\sqrt n}$, $\ell \leq 0$ on $\partial B_{4\sqrt n}$. However, as $\ell(x_0) = w(x_0) \geq -m$ with $x_0\in Q_1$, we deduce that $|D\ell| \leq m$, where $c_n>0$ depends only on $n$. Thus, 
\begin{equation}\label{eq:w-Lip}
w(x) \geq w(x_0) - m |x-x_0|, 
\end{equation} 
for all $x\in\Omega$. 

Next, we observe that we have freedom to choose $\vp$ in such a way that $\|D\vp\|_{L^\infty(B_{4\sqrt n})} \leq M$. To see this, note that $\vp$ is constructed in such a way that $\vp (x) = m_1 - m_2 |x|^{-\alpha}$ for $x\in B_{4\sqrt n}\setminus B_{1/4}$, and it is extended smoothly inside $B_{1/4}$ such that $\cP_- (D^2\vp ) \leq c_0\xi$ in $B_{4\sqrt n}$, where $\xi \in C(B_{4\sqrt n})$ is a continuous function with $0\leq \xi\leq 1$ and $\spt \xi\subset \overline{Q_1}$; here all constants $m_1$, $m_2$, $\alpha$ and $c_0$ depend only on $n$, $\lambda$ and $\Lambda$. Since $|D\vp|\leq m$ in $B_{4\sqrt n}\setminus B_{1/4}$, and the extension leaves the gradient free, we can find an extension such that 
\begin{equation}\label{eq:vp-Lip}
\sup_{B_{4\sqrt n}}|D\vp| \leq m,
\end{equation}
by taking $m$ larger if necessary. 

Finally, by the definition of $w$, we deduce from \eqref{eq:w-Lip}, \eqref{eq:vp-Lip} and the assumption that $u(x) \geq -|x|^2$ for all $x\in\Omega\setminus B_{2\sqrt n}$ that 
$$
u(x) \geq  u(x_0) - m|x-x_0|,
$$
for all $x\in\Omega$. This proves that $x_0\in Q_1\setminus \underline L_m(u,\Omega)$, as desired.
\end{proof} 

Now we may argue as in \cite[Lemma 7.7]{CC} to deduce a universal decay of the measure of the set with large gradient ``from below". 

\begin{lemma}\label{lemma:iter-W1d}
Let  $u\in C(B_{4\sqrt n})$ be an $L^p$-viscosity solution to $\cP_- (D^2 u ) \leq f$ in $B_{4\sqrt n}$, for some $f\in L^p(B_{4\sqrt n})$. Suppose that $\inf_{B_{4\sqrt n}} u \geq -1$, $\inf_{Q_3} u \leq 0$, and $\|f \|_{L^p(B_{4\sqrt n})} \leq \delta_0$. Let $L_k$ and $B_k$ denote $\underline L_{m^k} (u,B_{4\sqrt n})\cap Q_1$ and respectively $\{ M(|f|^p\chi_{B_{4\sqrt n}}) > \delta_0^p m^{kp}\}$. Then for each integer $k\geq 1$, 
$$
|L_{k+1}| \leq \sigma |L_k\cup B_k|, 
$$
where $m>1$, $\delta_0>0$ and $\sigma\in(0,1)$ depend only on $n$, $\lambda$, $\Lambda$ and $p$. 
\end{lemma}

\begin{proof}
The proof is almost the same with that of \cite[Lemma 7.7]{CC}. The involvement of the Riesz potential, which replaces the maximal function in the statement of the latter lemma, is due to the linear rescaling of the solutions. 

Fix an integer $k\geq 1$. Due to Lemma \ref{lemma:apprx-W1d}, $|L_1| \leq \sigma$. As $L_{k+1}\subset L_k \subset\cdots\subset L_1$, we have $|L_{k+1}| \leq \sigma$. Hence, due to the Calder\'on-Zygmund cube decomposition lemma, it suffices to prove that given any dyadic cube $Q\subset Q_1$, if $|L_k \cap Q| > \sigma |Q|$, then $\tilde Q\subset L_k\cup B_k$ where $\tilde Q$ is the predecessor of $Q$. 

Suppose, by way of contradiction, that $\tilde Q\setminus (L_k\cup B_k)\neq\emptyset$. Let $x_Q$ and $s_Q$ be the center and respectively the side-length of $Q$, i.e., $Q = Q_{s_Q} (x_Q)$. Let us consider the rescaled version of $u$ and respectively $f$,
$$
u_Q (x) = \frac{u(x_Q + s_Q x) - u(x_Q)}{cs_Q m^k},\quad f_Q (x) = \frac{s_Q}{cm^k} f(x_Q + s_Q x),
$$
with $c>1$ being a constant to be determined solely by $n$ and $p$, and let $\Omega_Q = s_Q^{-1}( -x_Q + B_{4\sqrt n})$. Choose any $\tilde x_Q \in \tilde Q\setminus (L_k\cup B_k)$. Then since $\tilde x_Q \in Q_{3s_Q}(x_Q)$, one can easily verify from $\tilde x_Q \not\in L_k$ that 
\begin{equation}\label{eq:iter-W1d-1}
\inf_{B_{4\sqrt n}} u_Q \geq -1, \quad \inf_{Q_3} u_Q \leq u_Q(0) = 0,
\end{equation}
where the first inequality is ensured by choosing $c>1$ large, depending only on $n$. Moreover, since $B_{4s_Q\sqrt n}(x_Q)\subset B_{6s_Q\sqrt n}(\tilde x_Q)$, $s_Q < 1$ and $\tilde x_Q\not\in B_k$, we also obtain 
\begin{equation}\label{eq:iter-W1d-2}
\|f_Q\|_{L^p(B_{4\sqrt n})} \leq \frac{c_0s_Q}{cm^k}(M(|f|^p\chi_{B_{4\sqrt n}})(\tilde x_Q))^{\frac{1}{p}}  \leq  \delta_0, 
\end{equation}
provided that we choose $c > c_0$. Furthermore,  
\begin{equation}\label{eq:iter-W1d-3}
\cP_- (D^2 u_Q) \leq f_Q\quad\text{in }B_{4\sqrt n},
\end{equation} 
in the $L^p$-viscosity sense.

Thanks to \eqref{eq:iter-W1d-1} -- \eqref{eq:iter-W1d-3}, $u_Q$ and $f_Q$ fall under the setting of Lemma \ref{lemma:apprx-W1d}, from which we deduce that 
$$
|\underline L_{c^{-1}m}(u_Q,\Omega_Q)\cap Q_1| \leq \sigma
$$
by choosing $m>1$ larger from the beginning so that $c^{-1} m$ becomes the constant appearing in the latter lemma. Rescaling back, we arrive at $|L_{k+1} \cap Q|\leq \sigma |Q|$, a contradiction to the choice of $Q$. Thus, the proof is finished. 
\end{proof}

As a corollary, we obtain a universal decay estimate in the interior. We shall only present the statement and omit the proof, as it being essentially the same with that of \cite[Lemma 7.8]{CC}. 

\begin{lemma}\label{lemma:int-W1d}
Under the setting of Lemma \ref{lemma:iter-W1d}, for any $t>0$, 
$$
|\underline L_t (u,\Omega)\cap Q_1| \leq ct^{-\mu},
$$
where $c>1$ and $\mu>0$ depend at most on $n$, $\lambda$, $\Lambda$ and $p$. 
\end{lemma}

From now on, we shall study the estimates near boundaries. As mentioned earlier, the idea to combine the interior estimate with the small measure of the thin slab that contains the boundary layer is originally from \cite[Lemma 2.9]{W}; here we simply extend the argument to the framework of Reifenberg flat domains. 

\begin{lemma}\label{lemma:apprx-W1d-bdry}
Let $\Omega\subset\R^n$ be a bounded domain with $0\in\partial\Omega$ such that $\partial\Omega\cap B_2$ is $(\delta,1)$-Reifenberg flat, and $u\in C(\Omega)$ be an $L^p$-viscosity solution to $\cP_-(D^2 u)\leq f $ in $\Omega$ for some $f\in L^p(\Omega)$, for some $p>p_0$. Suppose that $\|f \|_{L^p(\Omega\cap B_2)} \leq 1$, $\inf_{\Omega\cap B_2} u\geq - 1$, $\inf_{\Omega\cap B_1} u \leq 0$ and $u(x) \geq -|x|$ for all $x\in\Omega\setminus B_1$. Then
$$
| \underline L_t(u,\Omega)\cap B_1| \leq c(\delta^{-\mu} t^{-\mu} + \delta), 
$$
for any $t>0$, where $c>0$ and $\mu>0$ depend at most on $n$, $\lambda$, $\Lambda$ and $p$.  
\end{lemma} 

\begin{proof}
Since $\partial\Omega\cap B_2$ is $(\delta,1)$-Reifenberg flat and $0\in\partial\Omega$, there exists a unit vector $\nu\in\R^n$ such that $\partial\Omega\cap B_1\subset\{ x \in B_1 :| x\cdot \nu | < \delta\}$. 
Due to (a properly rescaled form of) Lemma \ref{lemma:int-W1d}, we have, for any $t>1$, 
$$
| \underline L_t(u,\Omega) \cap \{x\in B_1: x \cdot \nu > 2\delta\}| \leq c \delta^{-\mu} t^{-\mu}.
$$
The conclusion follows easily from the observation that $|\{x\in B_1: |x\cdot \nu| < 2\delta\}| \leq c_n\delta$, and $\underline L_t(u,\Omega)\subset \Omega$ (hence $L_t(u,\Omega)\cap B_1 = L_t(u,\Omega)\cap \{x\in B_1: x\cdot \nu > -\delta\}$).
\end{proof}

Next, we obtain universal decay estimates near boundary layers. 

\begin{lemma}\label{lemma:iter-W1d-bdry}
Let $\Omega\subset \R^n$ be a bounded domain with $0\in\partial\Omega$ such that $\partial\Omega\cap B_2$ is $(\delta_0,2)$-Reifenberg flat, and $u\in C(\Omega\cap B_2)$ be an $L^p$-viscosity solution to $\cP_-(D^2 u) \leq f$ in $\Omega\cap B_2$ for some $f\in L^p(\Omega\cap B_2)$, for some $p>p_0$. Denote by $L_k$ and $B_k$ the sets $\underline L_{m^k}(u,\Omega\cap B_2)\cap B_1$ and respectively $\{ M (|f|^p\chi_{\Omega\cap B_2}) > m^{kp}\}$, where $I_p$ is the Riesz potential of order $p$. Suppose that $\| u\|_{L^\infty(\Omega\cap B_2)} \leq 1$ and $\|f \|_{L^p(\Omega\cap B_2)}\leq 1$. Then for each integer $k\geq 1$, 
$$
| L_{k+1} | \leq \frac{1}{2} | L_k \cup B_k |,
$$
where $m>1$ and $\delta_0 \in(0,\frac{1}{2})$ are constants depending at most on $n$, $\lambda$, $\Lambda$ and $p$. 
\end{lemma}

\begin{proof}
Let $m>1$ be a constant to be determined later, and set $\eta = c_1 c_0(\delta_0^{-\mu} m^{-\mu} + \delta_0)$, where $c_0>1$ and $\mu>0$ are as in Lemma \ref{lemma:apprx-W1d-bdry}, and $c_1 > 1$ is a constant to be determined later, by $n$, $\lambda$, $\Lambda$ and $p$ only. Fix any integer $k\geq 1$. Then it follows from the latter lemma, as well as the relation $L_{k+1}\subset L_k\subset\cdots \subset L_1$, that $|L_{k+1}| \leq \eta$. 

Fix any ball $B\subset B_1$ with center in $\overline\Omega\cap B_1$ and $\rad B\leq 1$. Suppose that $|L_{k+1}\cap B| > \eta |B|$. We claim that 
\begin{equation}\label{eq:iter-W1d-bdry-claim}
\Omega\cap B\subset L_k\cup B_k. 
\end{equation}
Assume for the moment that the claim is true. Then it follows from Lemma \ref{lemma:vitali} (along with $\delta_0 < \frac{1}{2}$), that $|L_{k+1}| \leq c_n \eta |L_k\cup B_k|$. Then we first choose $\delta_0$ sufficiently small such that $4c_nc_0 \delta_0 \leq 1$. Selecting $m$ accordingly large such that $4c_nc_0\delta^{-\mu} m^{-\mu} \leq 1$, we obtain that $c_n\eta = c_nc_0(\delta_0^{-\mu } m^{-\mu} + \delta_0) \leq \frac{1}{2}$, which finishes the proof. 

Henceforth, we shall prove the claim \eqref{eq:iter-W1d-bdry-claim}. Suppose by way of contradiction that $\Omega\cap B\setminus (L_k\cup B_k) \neq\emptyset$. Here it suffices to consider the case $2B\setminus (\Omega\cap B_2)\neq\emptyset$, since the other case can be handled as in the interior analysis (see the proof of Lemma \ref{lemma:iter-W1d}). 

Set $r_B = \rad B$ and choose $x_B\in \partial \Omega\cap B_1$ in such a way that $2B\subset B_{4r_B}(x_B)$. Set $\Omega_B = -\frac{1}{2r_B} (-x_B +\Omega\cap B_2)$ and rescale $u$ and $f$ as follows,
$$
u_B(x) = \frac{u(x_B + 2r_B x) - u(x_B)}{cm^k r_B},\quad f_B (x) = \frac{r_B}{cm^k} f(x_B +2r_Bx). 
$$
Arguing analogously as in the proof of Lemma \ref{lemma:iter-W1d}, we may deduce from $\Omega\cap B\setminus (L_k\cup B_k) \neq \emptyset$ that  $\inf_{\Omega\cap B_2} u_B \geq -1$, $ \inf_{\Omega\cap B_1} u_B \leq 0$, $u_B (x) \geq -|x|\quad\text{for all }x\in\Omega_B\setminus B_1$, and $\|f _B\|_{L^p(\Omega_B\cap B_2)} \leq 1$, provided that $c>1$ is a large constant, depending at most on $n$ and $p$. Moreover, 
$$
\cP_- (D^2 u_B ) \leq f_B\quad\text{in }\Omega_B\cap B_2,
$$ 
in the $L^p$-viscosity solution. Thanks to the scaling invariance of the Reifenberg flatness, $\partial\Omega_B\cap B_2$ is $(\delta_0,2)$-Reifenberg flat and contains the origin. Thus, we can employ Lemma \ref{lemma:apprx-W1d-bdry} to deduce that 
$$
|\underline L_t (u_B ,\Omega_B) \cap B_1| \leq c_0(\delta_0^{-\mu} t^{-\mu} + \delta_0), 
$$
for any $t>0$. Rephrase the above inequality in terms of $u$, and deduce that $|\underline L_{cm^kt} (u,\Omega) \cap B_{2r_B}(x_B)| \leq c_0 (\delta_0^{-\mu} t^{-\mu} + \delta_0) |2r_B|^n$. As $B \subset B_{2r_B}(x_B)$ and $r_B = \rad B$, we derive that 
$$
|\underline L_{cm^kt} (u,\Omega)\cap B| \leq \frac{2^nc_0}{\omega_n}  (\delta_0^{-\mu} t^{-\mu} + \delta_0) |B|,
$$
where $\omega_n$ is the volume of the $n$-dimensional unit ball. Evaluating this inequality at $t= c^{-1} m$, we reach contradiction against $|L_{k+1}\cap B_r(x_0)| > \eta |B_r(x_0)|$ with $\eta = 2^n\omega_n^{-1} c^\mu c_0(\delta_0 m^{-\mu} + \delta_0)$ (i.e., $c_1 = 2^n\omega_n^{-1} c^\mu$ from the notation in the beginning of the proof). 
\end{proof}

Now we have a boundary-analogue of Lemma \ref{lemma:int-W1d}. Let us skip the proof for the same reason as mentioned above the statement of the latter lemma. 

\begin{lemma}\label{lemma:bdry-W1d}
Under the same hypothesis of Lemma \ref{lemma:iter-W1d-bdry}, for any $t>0$, 
$$
|\underline L_t (u,\Omega)\cap B_1| \leq ct^{-\mu}, 
$$
where $c>1$ and $\mu>0$ depend at most on $n$, $\lambda$, $\Lambda$ and $p$. 
\end{lemma}

Finally, we are ready to prove the global universal decay asserted in the beginning of this section.

\begin{proof}[Proof of Proposition \ref{proposition:W1d}]
With Lemma \ref{lemma:int-W1d} and Lemma \ref{lemma:bdry-W1d}, the assertion of this proposition follows easily via a standard covering argument. The exponent $\mu$ can be taken as the minimum between those in both lemmas. We omit the detail.
\end{proof} 


\subsection{Set of Large Hessian}\label{section:W2d}

Here we shall study universal decay estimates for the set of large Hessian. Note that an interior estimate is by now considered classical, and can be found in \cite[Lemma 7.8]{CC}, while an estimate near flat boundary is established rather recently in \cite{W}. Here we extend the result to Reifenberg flat boundaries. 

Here given any $u\in C(\Omega)$ and $t>0$, $\underline A_t(u,\Omega)$ is defined, as in \cite[Section 7]{CC}, as a subset of $\Omega$ such that $x_0\in \Omega\setminus \underline A_t(u,\Omega)$ if and only if there exists a linear polynomial $\ell$ for which $u(x) \geq \ell(x) - \frac{t}{2} |x-x_0|^2$ for all $x\in\Omega$. 

\begin{proposition}\label{proposition:W2d}
Let $\Omega\subset\R^n$ be a bounded, $(\delta,R)$-Reifenberg flat domain with $R\in(0,1]$, $f\in L^p(\Omega)$ be given, with some $p>p_0$, and $u\in C(\Omega)\cap L^\infty(\Omega)$ be an $L^p$-viscosity solution to $\cP_- (D^2 u) \leq f$ in $\Omega$. Then for any $t>0$,
$$
|\underline A_t (u,\Omega)| \leq \frac{C|\Omega|}{R^\mu t^\mu} (\|u \|_{L^\infty(\Omega)} +\| f\|_{L^p(\Omega)})^\mu,
$$
where $C>1$, $\delta >0$ and $\mu>0$ depend at most on $n$, $\lambda$, $\Lambda$ and $p$. 
\end{proposition}

\begin{proof}
Since the proof repeats most of the argument presented in the section above, we shall pinpoint the difference and skip the detail. First, observe that we can replace $\underline L_t$ in Lemma \ref{lemma:apprx-W1d-bdry} with $\underline A_t$, by simply applying \cite[Lemma 7.5]{CC} (which holds equally well for $L^p$-viscosity solutions with $L^p$-integrable right hand side, due to \cite{CCKS}) in place of Lemma \ref{lemma:apprx-W1d} in the proof. Now the assertion of Lemma \ref{lemma:iter-W1d-bdry} holds true with $L_k$ now denoting the set $\underline A_{m^k}(u,\Omega)\cap B_1$, since the proof only uses Lemma \ref{lemma:apprx-W1d}. Finally, iterating the modified version of Lemma \ref{lemma:iter-W1d-bdry} would yield Lemma \ref{lemma:bdry-W1d}, again with $\underline L_t(u,\Omega)$ replaced by $\underline A_t(u,\Omega)$. Thus, a standard covering argument along with the $L^p$-variant of \cite[Lemma 7.5]{CC} and the modified version of Lemma \ref{lemma:bdry-W1d} would yield the conclusion of this proposition. 
\end{proof}


\section{Uniform $W^{1,\frac{np}{n-p}}$-Estimates}\label{section:W1p}


\subsection{Estimates in the Interior}\label{section:int-W1p}

This section is devoted to uniform interior $W^{1,\frac{np}{n-p}}$-estimates in fully nonlinear homogenization problems, for any $p \in(p_0,n)$; note that $\frac{np}{n-p}$ is the (critical) Sobolev exponent of $p$. The estimate is optimal, and is even new in the context of standard fully nonlinear problems.

\begin{theorem}\label{theorem:int-W1p}
Let $F\in C(\cS^n\times\R^n)$ be a functional satisfying \eqref{eq:F-ellip} -- \eqref{eq:F-0}, $\Omega\subset\R^n$ be a bounded domain, $f \in C(\Omega)\cap L^p(\Omega)$, for some $p \in (p_0,n)$, and $u^\e\in C(\Omega)$. for some $\e>0$, be a viscosity solution to
\begin{equation}\label{eq:main-int-W1p}
F\left( D^2 u^\e , \frac{\cdot}{\e} \right) = f \quad\text{in }\Omega.
\end{equation}
Then $G_\Omega^\e (u^\e) \in L_{loc}^{\frac{np}{n-p}}(\Omega)$, and for any subdomain $\Omega'\Subset\Omega$, 
$$
\| G_\Omega^\e (u^\e) \|_{L^{\frac{np}{n-p}}(\Omega')} \leq C \left( \frac{\| u^\e \|_{L^\infty(\Omega)}}{\dist(\Omega',\partial\Omega)^{2-\frac{n}{p}}} +  \|f \|_{L^p(\Omega)}\right),
$$
where $C>0$ depends only on $n$, $\lambda$, $\Lambda$, $q$ and $p$. Assume further that \eqref{eq:F-C} holds for some constant $\kappa>0$ and some modulus of continuity $\psi: (0,\infty)\to(0,\infty)$. Then the same assertion holds with $G_\Omega(u^\e)$, hence with $|Du^\e|$, in place of $G_\Omega^\e( u^\e)$, in which case the constant $C$ may depend additionally on $\kappa$ and $\psi$. 
\end{theorem}

In what follows, we shall present our argument with $\Omega = B_{4\sqrt n}$ and $\Omega' = Q_1$, as our analysis will be of local character. Also, unless stated otherwise, we shall always assume that $F$ is a continuous functional on $\cS^n\times\R^n$ satisfying \eqref{eq:F-ellip} -- \eqref{eq:F-0}, $u^\e$ is a viscosity solution to \eqref{eq:main-int-W1p} with $\Omega$ replaced by $B_{4\sqrt n}$, for some $\e>0$, and $f\in C(B_{4\sqrt n})\cap L^p(B_{4\sqrt n})$ for some $p>p_0$. 

Let us begin with an approximation lemma for the measure of the set with large ``gradient". 

\begin{lemma}\label{lemma:apprx-W1p}
Let $\Omega\subset\R^n$ be a bounded domain such that $B_{4\sqrt n}\subset \Omega$, and suppose that $0<\e<1$, $\| u^\e \|_{L^\infty(B_{4\sqrt n})} \leq 1$ and $|u^\e(x)|\leq |x|$ for all $x\in\Omega\setminus B_{2\sqrt n}$. Then for any $s>N$, 
$$
|L_s^\e (u^\e,\Omega)\cap Q_1| \leq cs^{-\mu} \| f\|_{L^p(B_{4\sqrt n})}^\mu,
$$
where $N>1$ and $\mu>0$ depend at most on $n$, $\lambda$, $\Lambda$ and $p$. 
\end{lemma} 

\begin{proof}
Consider an auxiliary boundary value problem,
$$
\begin{dcases}
F\left(D^2 h^\e , \frac{\cdot}{\e} \right) = 0 &\text{in }B_{4\sqrt n},\\
h^\e = u^\e &\text{on }\partial B_{4\sqrt n}.
\end{dcases}
$$
As $u^\e \in C(\partial B_{4\sqrt n})$ and $F\in C(\cS^n\times\R^n)$ satisfying \eqref{eq:F-ellip}, there exists a unique viscosity solution $h^\e\in C(\overline B_{4\sqrt n})$ to the above problem. By the maximum principle, $\| h^\e\|_{L^\infty(B_{4\sqrt n})} \leq \| u^\e \|_{L^\infty(\partial B_{4\sqrt n})} \leq 1$. Now it follows from Lemma \ref{lemma:int-C1a}, along with the Kyrlov theory \cite[Corollary 5.7]{CC}, that for any $x_0\in B_{3\sqrt n}$, there exists a linear polynomial $\ell_{x_0}^\e$ such that $|D\ell_{x_0}^\e |\leq c$ and $|(h^\e - \ell_{x_0}^\e )(x)|\leq c (|x - x_0|^{1+\alpha} + \e^{1+\alpha})$ for all $x\in B_{4\sqrt n}$, for some $c>1$ and $\bar\alpha\in(0,1)$ depending only on $n$, $\lambda$ and $\Lambda$. Thus, 
\begin{equation}\label{eq:Lehe-msr}
L_c^\e (h^\e, B_{4\sqrt n})\cap Q_1 = \emptyset,
\end{equation}
by taking $c>1$ slightly larger if necessary. 

On the other hand, due to Lemma \ref{lemma:visc}, one can compute that $w^\e =  \delta^{-1} (u^\e - h^\e)$ satisfies 
$$
\begin{dcases}
\cP_-(D^2 w^\e ) \leq \frac{f}{\delta} \leq \cP_+ (D^2 w^\e) & \text{in }B_{4\sqrt n}, \\
w^\e = 0 &\text{on }\partial B_{4\sqrt n},
\end{dcases}
$$
in the viscosity sense. Since we assume $\|f \|_{L^p(B_{4\sqrt n})} \leq \delta$, it follows from the general maximum principle that $\| w^\e \|_{L^\infty(B_{4\sqrt n})} \leq c$. Therefore, one can deduce from Proposition \ref{proposition:W1d} to both $w^\e$ and $-w^\e$ that $|L_t (w^\e,B_{4\sqrt n})| \leq ct^{-\mu}$ for all $t>0$. Rephrasing this inequality in terms of $u^\e$, we deduce that 
\begin{equation}\label{eq:Lsue-he-msr}
|L_s (u^\e - h^\e, B_{4\sqrt n})| \leq c\delta^\mu s^{-\mu},
\end{equation}
for any $s>0$. Thus, the conclusion follows from \eqref{eq:Lehe-msr} and \eqref{eq:Lsue-he-msr}, as well as the assumption that $|u^\e(x)|\leq |x|$ for all $x\in\Omega\setminus B_{2\sqrt n}$. 
\end{proof} 

The following lemma is an analogue of \cite[Lemma 7.12]{CC}.

\begin{lemma}\label{lemma:iter-W1p}
Suppose that $\|f \|_{L^p(B_{4\sqrt n})} \leq \delta$ for some $\delta \in (0,1)$, and that $\| u^\e \|_{L^\infty(B_{4\sqrt n})} \leq 1$. Let $L_k^\e$ and $B_k$ denote $L_{m^k}^\e (u^\e,\Omega)\cap Q_1$ and respectively $\{ I_p(|f|^p \chi_{B_{4\sqrt n}}) > \delta^p m^{kp}\}$. Then for any integer $k\geq 1$, 
$$
|L_{k+1}^\e| \leq \delta^\mu |L_k^\e \cup B_k|,
$$
where $M>1$ and $\mu>0$ are constants depending at most on $n$, $\lambda$, $\Lambda$ and $p$. 
\end{lemma}

\begin{proof}
Fix an integer $k\geq 1$. Let $M>N$ be a large constant such that $c_0M^{-\mu} < 1$, with $c_0,N>1$ and $\mu>0$ as in Lemma \ref{lemma:apprx}; note that $M$ depends only on $n$, $\lambda$, $\Lambda$ and $p$. By Lemma \ref{lemma:apprx} (along with $c_0M^{-\mu} < 1$), $|L_1^\e|\leq \delta^\mu$. Since $L_{k+1}^\e\subset L_k^\e\subset\cdots\subset L_1^\e$, we have $|L_{k+1}^\e|\leq \delta^\mu < 1$. 

The rest of our proof will resemble that of \cite[Lemma 7.12]{CC}. Let $Q\subset Q_1$ be a dyadic cube such that $|L_{k+1}^\e \cap Q| > \delta^\mu |Q|$. We claim that $\tilde Q\subset L_k^\e\cup B_k$, where $\tilde Q$ is the predecessor of $Q$. Once this claim is justified, the conclusion is ensured by the Calder\'on-Zygmund cube decomposition lemma. 

Suppose, by way of contradiction, that $\tilde Q\setminus (L_k^\e\cup B_k)\neq \emptyset$. Denote by $x_Q$ and $s_Q$ the center and respectively the side-length of $Q$; i.e., $Q = Q_{s_Q}(x_Q)$. Choose any point $\tilde x_Q\in \tilde Q\setminus (L_k^\e\cup B_k)$. Now since $|\tilde x_Q - x_Q|\leq \frac{3}{2}s_Q\sqrt n$, we have $B_{4s_Q\sqrt n}(x_Q)\subset B_{6s_Q\sqrt n}(\tilde x_Q)$. 

Let us first remark that $s_Q > \e$, since $L_{k+1}^\e\cap Q\neq\emptyset$. The reason is as follows. Suppose that $s_Q \leq \e$. Then since $\tilde x_Q\in \tilde Q\setminus L_k^\e$, there exists some constant $a\in\R$ for which $|u^\e (x) - a|\leq m^k (|x-\tilde x_Q| + \e)$ for all $x\in B_{4\sqrt n}$. Therefore, as $\diam(\tilde Q) = 2s_Q\sqrt n < 2\e\sqrt n$, it follows from the latter inequality that for any $x_0\in\tilde Q$, $|u^\e (x) - a| \leq m^k (|x-x_0| + (2\sqrt n + 1)\e) \leq m^{k+1} (|x-x_0| + \e)$ for all $x\in B_{4\sqrt n}$, provided that $m > 2\sqrt n + 1$. This implies that $\tilde Q\cap L_{k+1}^\e = \emptyset$, a contradiction.

Since $\tilde x_Q\not\in L_k$, we can choose a constant $a$ for which $|u^\e(x) - a |\leq m^k(|x-\tilde x_Q| + \e)$ for all $x\in B_{4\sqrt n}$. This combined with $|\tilde x_Q - x_Q|\leq \frac{3}{2}s_Q\sqrt n$ and $\e<s_Q$ yields that 
\begin{equation}\label{eq:iter-W1p-1}
|u^\e (x) - a| \leq m^k \left(\frac{5}{2}s_Q\sqrt n + |x-x_Q| \right),
\end{equation}
for all $x\in B_{4\sqrt n}$. In addition, due to the assumption $\tilde x_Q\not\in B_k$, as well as $B_{4s_Q\sqrt n}(x_Q)\subset B_{6s_Q\sqrt n}(\tilde x_Q)$, one can compute that
\begin{equation}\label{eq:iter-W1p-2}
\begin{aligned}
\int_{B_{4s_Q \sqrt n}(x_Q)\cap B_{4\sqrt n}} |f(x)|^p\,dx & \leq \int_{B_{6s_Q \sqrt n}(\tilde x_Q)\cap B_{4\sqrt n}} |f(x)|^p\,dx  \\
&\leq (6s_Q\sqrt n)^{n-p} \int_{B_{4\sqrt n}} \frac{|f(x)|^p}{|x-\tilde x_Q|^{n-p}}\,dx\\
&\leq c^p s_Q^{n-p} \delta^p,
\end{aligned}
\end{equation}
where $c>1$ depends only on $n$ and $p$. 

In what follows, we shall use $c > 1$ to denote a positive constant depending at most on $n$ and $p$, and allow it to vary at each occurrence.

Let us consider the following rescaled versions of $u^\e$ and $f$, 
$$
\begin{aligned}
&u_Q^{\e_Q} (x) = \frac{u^\e (x_Q + s_Q x) - a}{cm^ks_Q}, \quad \e_Q = \frac{\e}{s_Q}, \\
&f_Q (x) = \frac{s_Q}{cm^k} f(x_Q + s_Q x).
\end{aligned}
$$
Setting $\Omega_Q = s_Q^{-1}( - x_Q + B_{4\sqrt n})$, we have $B_{4\sqrt n}\subset \Omega_Q$, and thus in view of \eqref{eq:main-int-W1p}, $u_Q^{\e_Q}$ is a viscosity solution to 
$$
F_Q \left( D^2 u_Q^{\e_Q} , \frac{\cdot}{\e_Q}\right) = f_Q\quad\text{in }B_{4\sqrt n},
$$
where $F_Q (P,y) = \frac{s_Q}{cm^k} F( \frac{cm^k}{s_Q} P,y + \frac{x_Q}{\e})$. Clearly, $F_Q\in C(\cS^n\times\R^n)$ and it satisfies \eqref{eq:F-ellip} -- \eqref{eq:F-0}. On the other hand, it follows immediately from \eqref{eq:iter-W1p-1} and \eqref{eq:iter-W1p-2} that with $c>1$ sufficiently large, $|u^{\e_Q}(x)|\leq 1$ for all $x\in B_{4\sqrt n}$, $|u^{\e_Q}(x)| \leq |x|$ for all $\Omega_Q\setminus B_{2\sqrt n}$, and $\|f_Q \|_{L^p(B_{4\sqrt n})} \leq \delta$. 

In all, $\e_Q$, $F_Q$, $u_Q^{\e_Q}$ and $f_Q$ fall into the setting of Lemma \ref{lemma:apprx-W1p}, from which we obtain 
$$
|L_s^{\e_Q} (u_Q^{\e_Q}, \Omega_Q)\cap Q_1| \leq c_0\delta^\mu s^{-\mu}, 
$$
for any $s>N$ with $N>1$ as in the latter lemma. Thus, taking $M> N$ larger if necessary such that $c_0c^{-\mu} M^{-\mu} < 1$, and then rescaling back to $u^\e$, we arrive at $|L_{k+1}^\e \cap Q| \leq \delta^\mu |Q|$, a contradiction. 
\end{proof} 

We are ready to prove the uniform ``large-scale'' interior $W^{1,\frac{np}{n-p}}$-estimate, which amounts to the first part of Theorem \ref{theorem:int-W1p}. 

\begin{proof}[Proof of Theorem \ref{theorem:int-W1p}; the first part]
After some suitable rescaling argument, it suffices to consider the case where $\Omega = B_{4\sqrt n}$, $\Omega' = Q_1$, $\|u^\e\|_{L^\infty(B_{4\sqrt n})} \leq 1$ and $\|f \|_{L^p(B_{4\sqrt n})} \leq \delta$, where $\delta\in(0,1)$ is to be determined by $n$, $\lambda$, $\Lambda$ and $p$ only. 

Set $p' = \frac{p + p_0}{2} \in (p_0, p)$, and apply Lemma \ref{lemma:iter-W1p}, with $p$ replaced by $p'$. Observe that $\|f \|_{L^{p'}(B_{4\sqrt n})} \leq c \|f\|_{L^p(B_{4\sqrt n})}\leq c\delta$ for some $c>0$ depending only on $n$ and $p$. Hence, with $\eta = (c\delta)^\mu < 1$ (by choosing $\delta$ smaller if necessary), $\alpha_k = |L_k^\e|$ and $\beta_k = |B_k| = |\{ I_{p'} (|f|^{p'}\chi_{B_{4\sqrt n}}) > (c\delta m^k)^{p'}\}\cap B_{4\sqrt n}|$, we obtain that 
\begin{equation}\label{eq:int-W1p-re}
\alpha_k \leq \eta (\alpha_{k-1}^\e + \beta_{k-1}) \leq \cdots \leq \eta^k + \sum_{i=1}^k \eta^i \beta_{k-i}.
\end{equation}
Now as $f\in L^p(B_{4\sqrt n})$, we have $|f|^{p'}\chi_{B_{4\sqrt n}} \in L^{\frac{p}{p'}}(\R^n)$ with $\frac{p}{p'} > 1$. Thus, according to the embedding theorem for the Riesz potential, $I_{p'} (|f|^{p'}\chi_{B_{4\sqrt n}}) \in L^{\frac{np}{p'(n-p)}}(\R^n)$, from which we can compute that 
\begin{equation}\label{eq:Riesz-Lp*}
\begin{aligned}
\sum_{k=1}^\infty M^{\frac{np}{n-p}k} \beta_k &\leq c\int_0^\infty t^{\frac{np}{n-p} - 1} |\{ I_{p'}(|f|^{p'}\chi_{B_{4\sqrt n}}) >  t^{p'}\}\cap B_{4\sqrt n}| \\
& = c\int_{B_{4\sqrt n}} (I_{p'} (|f|^{p'}\chi_{B_{4\sqrt n}}))^{\frac{np}{p'(n-p)}}\,dx \\
&\leq c.
\end{aligned}
\end{equation}

To this end, we choose $\delta\in(0,1)$ as a sufficiently small constant such that $M^{\frac{np}{n-p}} \eta = M^{\frac{np}{n-p}} (c\delta)^\mu \leq \frac{1}{2}$; clearly, it is the set of parameters $n$, $\lambda$, $\Lambda$ and $p$ that determines how small $\delta$ should be. Then it follows from \eqref{eq:int-W1p-re} and \eqref{eq:Riesz-Lp*} that 
$$
\sum_{k=1}^\infty M^{\frac{np}{n-p}k} \alpha_k \leq  \sum_{k=1}^\infty M^{\frac{np}{n-p}k}\eta^k + \sum_{k=1}^\infty \sum_{i=1}^k (M^{\frac{np}{n-p}i}\eta^i )(M^{\frac{np}{n-p}(k-i)} \beta_{k-i}) \leq c.
$$
Since $\{ |D_\e u^\e | >t \}\subset L_t^\e(u^\e, B_{4\sqrt n})$, we have proved that $\| D_\e u^\e \|_{L^{\frac{np}{n-p}}(Q_1)} \leq c$, as desired. 
\end{proof}

The proof for the second part of the theorem is more or less the same. We shall only point out where the argument changes and leave out the detail to the reader.

\begin{proof}[Proof of Theorem \ref{theorem:int-W1p}; the second part]
With the additional assumption in the statement of the theorem, we can replace $L_s^\e$ in Lemma \ref{lemma:apprx-W1p} with $L_s$. For we can now invoke the uniform interior $C^{1,\alpha}$-estimate below $\e$-scale, \cite[Theorem 4.1 (ii)]{KL}, to deduce that the approximating solution $h^\e$ belongs to $C^{1,\alpha}(B_{3\sqrt n})$ with $\| D h^\e \|_{C^\alpha(B_{3\sqrt n})} \leq c$, whence we can replace $L_N^\e(h^\e,\cdot)$ with $L_N (h^\e,\cdot)$ in \eqref{eq:Lehe-msr}. 

Now with $L_s^\e$ replaced with $L_s$ in the assertion of Lemma \ref{lemma:apprx-W1p}, we can derive the same conclusion in Lemma \ref{lemma:iter-W1p} with $L_k^\e = L_{m^k} (u^\e,B_{4\sqrt n})\cap Q_1$ instead. The proof actually becomes shorter in this case, since we no longer need to care about the case where the side-length of a given cube is smaller than $\e$. 

Finally, we can repeat the proof above for the first part of this theorem, verbatim, to obtain that $\sum_{k=1}^\infty M^{\frac{np}{n-p}k} \alpha_k \leq c$ with $\alpha_k = |L_{m^k} (u^\e,B_{4\sqrt n})\cap Q_1|$ now. This proves $\| Du^\e \|_{L^{\frac{np}{n-p}}(Q_1)} \leq c$. We omit the detail. 
\end{proof} 


\subsection{Estimates near Boundary Layers}\label{section:bdry-W1p}

Let us now turn to uniform $W^{1,\frac{np}{n-p}}$-estimates near boundary layers. 

\begin{theorem}\label{theorem:bdry-W1p}
Let $F \in C(\cS^n\times\R^n)$ be a functional satisfying \eqref{eq:F-ellip} -- \eqref{eq:F-0}, $\Omega\subset\R^n$ be a domain, $U\subset\R^n$ be an open neighborhood of a point of $\partial\Omega$ such that $\partial\Omega\cap U$ is $(\delta,R)$-Reifenberg flat, for some $\delta>0$ and some $R \in(0,1]$, $f\in L^p(\Omega\cap U)$ and $g\in W^{1,\frac{np}{n-p}}(U)$, for some $p \in (p_0,n)$. Let $u^\e\in C(\overline\Omega\cap U)$ be a viscosity solution to
\begin{equation}\label{eq:main-bdry-W1p}
\begin{dcases}
F\left(D^2 u^\e ,\frac{\cdot}{\e}\right) = f &\text{in }\Omega\cap U, \\
u^\e = g &\text{on }\partial\Omega\cap U.
\end{dcases}
\end{equation}
Then $G_{\Omega\cap U}^\e (u^\e) \in L_{loc}^{\frac{np}{n-p}}(\overline\Omega_\e\cap U)$, and for any $U'\Subset U$, 
$$
\begin{aligned}
\| G_{\Omega\cap U}^\e (u^\e) \|_{L^{\frac{np}{n-p}}(\Omega_\e\cap U')} \leq C \left( \frac{\| u^\e  \|_{L^\infty(\Omega\cap U)}}{\dist(U', \partial U)^{2-\frac{n}{p}}} + \| f \|_{L^p(\Omega\cap U)} + \|Dg \|_{L^{\frac{np}{n-p}}(U)} \right),
\end{aligned}
$$
where $\delta >0$ depends only on $n$, $\lambda$, $\Lambda$ and $p$, and $C>1$ may depend further on $R$ and $\diam(U)$. Moreover, under the same assumption in Theorem \ref{theorem:int-W1p}, the same assertion holds with $G_{\Omega\cap U}(u^\e)$, hence with $|Du^\e|$, and $\Omega\cap U'$ in place of $G_{\Omega\cap U}^\e (u^\e)$ and respectively $\Omega_\e\cap U'$. 
\end{theorem}

The above estimate is optimal as the power of integrability reaches the critical Sobolev exponent. This estimate is even new in the setting of the standard problems. To the best of the author's knowledge, the boundary estimate is proved up to the subcritical Sobolev exponent, i.e., $W^{1,q}$ with $q<\frac{np}{n-p}$, in \cite{W}. Following the spirit of \cite{S}, the proof in \cite{W} relies heavily on pointwise $C^{1,\alpha}$-approximation, and hence the estimate could not reach the critical exponent. 

We shall set the starting point of our analysis, however, at a sub-optimal estimate below. We shall provide some motivations and remarks after the statement. 

\begin{proposition}\label{proposition:bdry-W1p}
Let $F \in C(\cS^n\times\R^n)$ be a functional satisfying \eqref{eq:F-ellip} -- \eqref{eq:F-0}, $\Omega\subset\R^n$ be a domain, $U\subset\R^n$ be an open neighborhood of a point of $\partial\Omega$ such that $\partial\Omega\cap U$ is $(\delta,R)$-Reifenberg flat from exterior,\footnote{That is, $\partial\Omega\cap B_r(x_0)\subset \{x\in B_r(x_0): (x-x_0)\cdot \nu > -\delta r\}$ for any $r\in(0,R)$ and any $x_0\in\partial\Omega\cap U$.} for some $\delta>0$ and some $R \in (0,1]$ $f\in L^p(\Omega\cap U)$ and $g\in C^{0,2-\frac{n}{p}}(\Omega\cap U)$, for some $p \in (p_0,n)$. Let $u^\e\in C(\overline\Omega\cap U)$ be a viscosity solution to
$$
\begin{dcases}
F\left(D^2 u^\e ,\frac{\cdot}{\e}\right) = f &\text{in }\Omega\cap U, \\
u^\e = g &\text{on }\partial\Omega\cap U.
\end{dcases}
$$
Then $G_{\Omega\cap U}^\e (u^\e) u^\e \in L_{loc}^q(\overline\Omega_\e\cap U)$ for any $q\in[1,\frac{np}{n-p})$, and for any $U'\Subset U$, 
$$
\begin{aligned}
\frac{\| G_{\Omega\cap U}^\e (u^\e) \|_{L^q(\Omega_\e\cap U')}}{\dist(U',\partial U)^{1-\frac{n}{p} + \frac{n}{q}}} \leq C \left( \frac{\| u^\e  \|_{L^\infty(\Omega\cap U)}}{\dist(U', \partial U)^{2-\frac{n}{p}}} + \| f \|_{L^p(\Omega\cap U)} + [g ]_{C^{0,2-\frac{n}{p}}(\partial\Omega\cap U)} \right),
\end{aligned}
$$
where $\delta>0$ depends only on $n$, $\lambda$, $\Lambda$ and $p$, and $C>1$ may depend further on $q$, $R$ and $\diam(U)$. Moreover, under the same assumption in Theorem \ref{theorem:int-W1p}, the same assertion holds with $G_{\Omega\cap U}(u^\e)$, hence with $|Du^\e|$, and $\Omega\cap U'$ in place of $G_{\Omega\cap U}^\e (u^\e)$ and respectively $\Omega_\e\cap U'$. 
\end{proposition} 

We shall provide this estimate, mainly because of its independent interests. Of course, we will use this proposition in our subsequent analysis to ensure a better regularity for the approximating solutions. Still, this step could be conveniently replaced by the uniform boundary $C^{1,\alpha}$-estimates \cite{KL}, as the approximating solutions solve a ``clean'' version of the homogenization problems (i.e., $f = 0$ and $g = a$, with $ a$ a constant). 

In view of the compact embedding of $W^{1,\frac{np}{n-p}}$ to $C^{0,2-\frac{n}{p}}$, the above estimate seems to be the best one can expect with a $C^{0,2-\frac{n}{p}}$-regular boundary data. The estimate is quite interesting in the sense that a H\"older regular function may not be (weakly) differentiable a.e., that is, one cannot expect $C^{0,2-\frac{n}{p}}$-regular boundary data to be extended to a $W^{1,\frac{np}{n-p}}$-regular function in a neighborhood of the boundary layer. In particular, one cannot deduce the above sub-optimal estimate from the former optimal estimate (Theorem \ref{theorem:bdry-W1p}). Hence, the above proposition shows certain regularizing effect arising from the presence of boundary layers. 

Let us also remark that the sub-optimal estimate above improves the one in \cite{W}, in terms of the regularity of the boundary data, as the latter estimate requires $C^{1,\alpha}$-regularity.

We shall first present a complete proof for Proposition \ref{proposition:bdry-W1p} and then move onto that of Theorem \ref{theorem:bdry-W1p}. The proof for the former is based on a boundary $C^{0,2-\frac{n}{p}}$-estimate for any viscosity solution belonging to the Pucci class, up to an $L^p$-regular, with $p > p_0$, right hand side. In particular, the functional may not oscillate under certain pattern in small scales, whence it has nothing to do with homogenization. 

\begin{proof}[Proof of Proposition \ref{proposition:bdry-W1p}]
Let us prove the first part of the assertion, and then mention the changes in the argument for the second part, as the latter is almost the same with the former. 

After some standard rescaling procedure, one may prove the assertion for the case where $0\in\partial\Omega$, $U = B_1$, $U' = B_{1/2}$, $\| u^\e \|_{L^\infty(\Omega\cap B_1)} \leq 1$, $\|f \|_{L^p(\Omega\cap B_1)} \leq 1$ and $[g ]_{C^{0,2-\frac{n}{p}}(\partial\Omega\cap B_1)} \leq 1$. Applying the boundary H\"older regularity \cite[Theorem 1.1]{LZ2} to $u^\e$ around each point $x_0\in\partial\Omega\cap B_{1/2}$, we deduce that
\begin{equation}\label{eq:ue-g-Ca}
| u^\e (x) - g(x_0) | \leq c|x-x_0|^{2-\frac{n}{p}},
\end{equation}
for all $x\in\overline\Omega\cap B_1$; let us remark that although the statement of \cite[Theorem 1.1]{LZ2} involves $\sup_{r>0} r^{\alpha-2} \|f \|_{L^n(\Omega\cap B_r(x_0))}$ on the right hand side, one can easily replace this norm with $\|f \|_{L^p(\Omega\cap B_1)}$ by taking $\alpha = 2 - \frac{n}{p}$, as the modification in the proof there is straightforward. 

Let $B\subset\Omega\cap B_{1/2}$ be a ball for which $\partial (2B)\cap \partial\Omega\cap B_{1/2} \neq \emptyset$. Let $x_B$ and $\rho_B$ denote the center and respectively the radius of $B$. Also let $x_{B,0}$ be a point of intersection between $\partial (2B)$ and $\partial\Omega\cap B_{1/2}$. Consider the rescaling 
$$
u_B^{\e_B} (x) = \frac{u^\e (x_B + \rho_B x) - g(x_{B,0})}{\rho_B^{2-\frac{n}{p}}},
$$ 
of $u^\e$, where $\e = \rho_B^{-1} \e$. As $B_{2\rho_B}(x_B)\subset B_1$, one may observe that 
$$
F_B \left(D^2 u_B^{\e_B},\frac{\cdot}{\e_B} \right) = f_B \quad\text{in }B_1,
$$
in the viscosity sense, where we wrote $F_B (P,y) = \rho_B^{n/p}F( \rho_B^{-n/p} P,y)$ and $f_B (x) = \rho_B^{n/p}f (x_B + \rho_B x)$. Clearly, $F_B$ is a continuous functional on $\cS^n\times\R^n$ satisfying \eqref{eq:F-ellip} -- \eqref{eq:F-0}, and $f_B \in L^p(B_2)\cap C(B_2)$ with 
$$
\| f_B \|_{L^p(B_2)} \leq \|f \|_{L^p(\Omega\cap B_{\rho_B}(x_B))} \leq \delta.
$$
Thanks to \eqref{eq:ue-g-Ca}, we also have $\| u_B^{\e_B} \|_{L^\infty(B_2)} \leq c$. Therefore, we can apply Theorem \ref{theorem:int-W1p} (i) to observe that $\| G_{\Omega_B}^{\e_B} u_B^{\e_B} \|_{L^q(B_1)} \leq c$, for each $q\in[1,\frac{np}{n-p})$, where $\Omega_B = \rho_B^{-1}(-x_B + \Omega\cap B_1)$. The latter estimate can be translated in term of $u^\e$ as
\begin{equation}\label{eq:Deue-Lq-re}
\| G_{\Omega\cap B_1}^\e( u^\e) \|_{L^q(B_{\rho_B}(x_B))} \leq c\rho_B^{\frac{n}{q} + 1 -\frac{n}{p}}. 
\end{equation} 

To this end, consider a covering $\cF$ of $\Omega_\e\cap B_{1/2}$ by balls $B$ for which $\partial (2B)\cap \partial\Omega\cap B_{1/2}\neq\emptyset$; i.e., $\cF \subset \{B_{\dist(x,\partial\Omega\cap B_{1/2})}(x): x\in\Omega_\e\cap B_{1/2}\}$. Thanks to the Besicovitch covering theorem, one can find a finite subcovering $\cG\subset\cF$ such that $\#\{ \hat B\in\cG: B\cap \hat B\neq \emptyset\} \leq c_n$ for all $B\in\cG$, where $c_n>0$ depends only on $n$. One can split $\cG$ into a union of $\cG_k$, with $k=1,2,\cdots,k_\e$, such that $B\in \cG_k$ if $2^{-k-1} < \rho_B \leq 2^{-k}$, where $k_\e$ is a large positive integer of order $-\log_2\e$. Then 
\begin{equation}\label{eq:Deue-Lq-bdry-re}
\begin{aligned}
\int_{\Omega_\e\cap B_{1/2}} |G_{\Omega\cap B_1}^\e( u^\e)|^q \,dx &\leq \sum_{k=1}^{k_\e} \sum_{B\in \cG_k} \int_B |G_{\Omega\cap B_1}^\e( u^\e) |^q\,dx \\
&\leq c \sum_{k=1}^{k_\e} \sum_{B\in\cG_k} \rho_B^{n + q -\frac{nq}{p}} \\
&\leq c \sum_{k=1}^{k_\e} 2^{-n - q + \frac{nq}{p}} \\
&\leq c,
\end{aligned}
\end{equation} 
where the last inequality is ensured from the fact that $q < \frac{np}{n-p}$. This finishes the proof, for the first part of the assertion of the theorem. 

As for the second part of the assertion, one may have already noticed that under the additional assumption, one can invoke Theorem \ref{theorem:int-W1p} (ii) in place of (i) above, such that one can replace $D_\e u^\e$ with $D u^\e$ first in \eqref{eq:Deue-Lq-re}, and then in \eqref{eq:Deue-Lq-bdry-re}, leaving everything else untouched. Thus, the conclusion follows. We leave out the detail to the reader.
\end{proof} 

With Proposition \ref{proposition:bdry-W1p} at hand, we may now proceed with the proof for Theorem \ref{theorem:bdry-W1p}. The idea is basically the same with the interior case (Lemma \ref{lemma:apprx-W1p}), but the presence of Reifenberg flat boundaries yields some additional technical difficulties. 

As in the analysis in the interior case, we shall confine ourselves to the case $U = B_2$ and $U' = B_1$, as the analysis here being of local character around a boundary point. Unless specified otherwise, we shall always assume, from now on, that $F\in C(\cS^n\times\R^n)$ verifies  \eqref{eq:F-ellip} -- \eqref{eq:F-0}, $f\in C(\Omega\cap B_2)\cap L^p(\Omega\cap B_2)$, $g\in W^{1,\frac{np}{n-p}}(B_2)$, $\Omega$ is a bounded domain containing the origin, and $u^\e$ is a viscosity solution to \eqref{eq:main-bdry-W1p} with $U$ replaced by $B_2$. 

\begin{lemma}\label{lemma:apprx-W1p-bdry}
Assume that $\partial\Omega\cap B_2$ is $(\delta,2)$-Reifenberg flat, $\|f \|_{L^p(\Omega\cap B_2)} \leq \eta$, $\| Dg \|_{L^{\frac{np}{n-p}}(B_2)} \leq \eta$ for some $\eta>0$, and that $\| u^\e - a\|_{L^\infty(\Omega\cap B_2)} \leq 1$ and $| u^\e(x)- a | \leq |x|$ for all $x\in\Omega\setminus B_1$ for some constant $a\in\R$. Then for any $s>N$, 
$$
|L_s^\e (u^\e,\Omega)\cap \Omega_\e\cap B_1 | \leq c(\delta^{\frac{\alpha}{2}} + \eta)^\mu s^{-\mu},
$$ 
where $\alpha\in(0,1)$, $\mu>0$ and $N>1$ depend only on $n$, $\lambda$ and $\Lambda$, while $\delta>0$ and $c>1$ $N>1$ may depend further on $p$.  
\end{lemma} 

\begin{proof}
Let $\hat a$ denote the integral average of $g$ over $B_2$. By the Poincar\'e inequality, we have $\|g - \hat a\|_{W^{1,\frac{np}{n-p}}(B_2)} \leq c_0\| Dg\|_{L^{\frac{np}{n-p}}(B_2)}\leq c_0\eta$, where $c_0$ depends only on $n$ and $p$. By the Sobolev embedding theorem, $g\in C^{0,2-\frac{n}{p}}(\overline{B_2})$ and 
\begin{equation}\label{eq:hg-Ca}
\| g - \hat a \|_{L^\infty(B_2)} + [g]_{C^{0,2-\frac{n}{p}}(\overline{B_2})} \leq c_0\eta, 
\end{equation}
by taking $c_0$ larger. 

Denote by $S$ the slab $S_{2\delta}(\nu)$ for some unit vector $\nu$ such that $\partial\Omega\cap B_2 \subset S$; such a unit vector exists owing to the Reifenberg flatness of $\partial\Omega\cap B_2$. Also let $E$ be the half-space $E$, so that $\Omega\cap B_2\subset E$.

Since $u^\e = g$ on $\partial\Omega\cap B_2$ and $|u^\e |\leq 1$ on $\overline{\Omega\cap B_2}$, it follows from \eqref{eq:hg-Ca} that $\|g\|_{L^\infty(B_2)} \leq 1 + c_0\eta < 2$. Extend $u^\e$ by $g$ outside $\Omega$. It readily follows that $u^\e \in C(\overline{B_2})$ and $|u^\e - g|\leq 2$ in $B_2$. Consider a Dirichlet boundary value problem, 
$$
\begin{dcases}
F \left( D^2 h^\e,\frac{\cdot}{\e}\right) = 0 &\text{in }E\cap B_{3/2},\\
h^\e = u^\e - g& \text{on }\partial (E\cap B_{3/2}).
\end{dcases}
$$
This problem admits a unique viscosity solution $h^\e \in C(\overline{E\cap B_{3/2}})$ because $F$ and $u^\e - g$ are continuous, and $\partial (E\cap B_{3/2})$ satisfies a uniform exterior sphere condition. Since $|h^\e | \leq 2$ on $\partial (E\cap B_{3/2})$ and $h^\e = 0$ on $\partial E\cap B_2$, we can deduce from Lemma \ref{lemma:bdry-C1a}, the Krylov theory \cite[Corollary 5.7]{CC}, as well as the fact that $\partial E\cap B_{3/2}$ being flat, that for each $x_0\in \partial E\cap B_1$, there exists a linear polynomial $\ell_{x_0}^\e$ for which $|D\ell_{x_0}^\e|\leq c_0$ and $|(h^\e - \ell_{x_0}^\e)(x)| \leq c_0 ( |x-x_0|^{1+\alpha} + \e^{1+\alpha})$ for all $x\in \overline{E\cap B_{3/2}}$, with both $c_0>1$ and $\alpha \in (0,1)$ depending only on $n$, $\lambda$ and $\Lambda$. In particular,
\begin{equation}\label{eq:Lehe-Lq-bdry}
L_{c_0}^\e (h^\e, E_\e\cap B_{3/2}) \cap B_1 = \emptyset,
\end{equation}
by taking $c_0$ larger if necessary, with $E_\e = \{x\in E: \dist(x,\partial E)>0\}$. 

In what follows, we shall denote by $c_0$ a constant which may depend at most on $n$, $\lambda$, $\Lambda$ and $p$, and we shall allow it to vary at each occurrence. 

As $\partial\Omega\cap B_2$ being $(\delta,2)$-Reifenberg flat, one may invoke \cite[Theorem 1.1]{LZ} to deduce, along with $\|u^\e \|_{L^\infty(\Omega\cap B_2)}\leq 1$, $\|f \|_{L^p(\Omega\cap B_2)} \leq \eta$ and $[g]_{C^{0,2-\frac{n}{p}}(B_2)}\leq c_0\eta$ that if $\delta\leq\delta_p$, with $\delta_p$ depending only on $n$, $\lambda$, $\Lambda$ and $p$, then 
$$
|u^\e (x) - g(x_0) | \leq c_0|x-x_0|^{2-\frac{n}{p}},
$$
for all $x\in\Omega\cap B_2$, for each $x_0\in\partial\Omega\cap B_{3/2}$. Combining this inequality with \cite[Lemma 2]{Esc}, we deduce that 
\begin{equation}\label{eq:ue-Ca-bdry}
[u^\e ]_{C^{0,\alpha}(\overline{\Omega\cap B_{3/2}})} \leq c_0. 
\end{equation}

With \eqref{eq:ue-Ca-bdry} at hand, we may now estimate the global H\"older norm of $h^\e$. As $h^\e = 0$ on $\partial E \cap B_{3/2}$, $h^\e = u^\e - g$ on $E\cap \partial B_{3/2}$ and $[u^\e - g]_{C^{0,\alpha}(E\cap \partial B_{3/2})} \leq c_0$, we may apply a variant of \cite[Proposition 4.13]{CC} to derive that 
\begin{equation}\label{eq:he-Ca-bdry}
[h^\e ]_{C^{0,\frac{\alpha}{2}}(\overline{E \cap B_{3/2}})} \leq c_0. 
\end{equation}

Consider another auxiliary function $w^\e =\eta_0^{-1}(u^\e - \hat a - h^\e)$, with $\eta_0$ to be determined later. 
This function is well-defined on $\overline{\Omega\cap B_{3/2}}$, since $\Omega\cap B_{3/2}\subset E\cap B_{3/2}$. Due to Lemma \ref{lemma:visc}, we may compute that  
$$
\begin{dcases}
\cP_- (D^2 w^\e ) \leq \frac{f}{\eta_0} \leq \cP_+ (D^2 w^\e) & \text{in }\Omega\cap B_{3/2},\\
w^\e = \frac{g - \hat a}{\eta_0} &\text{on }\Omega\cap \partial B_{3/2},\\
w^\e = \frac{g -\hat a - h^\e}{\eta_0} &\text{on }\partial\Omega\cap B_{3/2},
\end{dcases}
$$
Thanks to \eqref{eq:he-Ca-bdry}, $h^\e = 0$ on $E\cap B_{3/2}$ and $\partial\Omega\cap B_{3/2}\subset H_\delta(-\nu)\cap E$, we may deduce that $| h^\e |\leq c_0\delta^{\frac{\alpha}{2}}$ on $\partial\Omega\cap B_{3/2}$. This together with \eqref{eq:hg-Ca} yields that with $\eta_0 =  c_0 \max\{\delta^{\frac{\alpha}{2}}, \eta\}$ and $c_0>1$, 
$$
\|w^\e \|_{L^\infty(\partial(\Omega\cap B_{3/2}))} \leq \frac{c_0\delta^{\frac{\alpha}{2}} + c_0\eta}{\eta_0} \leq 1.
$$
With such a choice of $\eta_0$, we also have $\|f \|_{L^p(\Omega\cap B_{3/2})} \leq \eta_0$, whence it follows from the general maximum principle that $\|w^\e\|_{L^\infty(B_{3/2})} \leq c_0$.

Now we can apply Proposition \ref{proposition:W1d} to obtain that for any $t>0$, 
\begin{equation}\label{eq:Ltwe-Ld-bdry}
|L_t (w^\e , \Omega\cap B_{3/2})\cap B_1| \leq ct^{-\mu},
\end{equation}
where $\mu>0$ depends only on $n$, $\lambda$, $\Lambda$ and $p$; this is another step that determines how small the Reifenberg flatness constant $\delta$ should be. Since the set $L_t(w^\e,\cdot)$ is invariant under vertical translation of given function $w^\e$, the conclusion follows immediately from \eqref{eq:Lehe-Lq-bdry},  \eqref{eq:Ltwe-Ld-bdry}, the choice of $\eta_0$ and the assumption that $|u^\e(x) | \leq |x|$ for all $x\in\Omega\setminus B_1$. 
\end{proof} 

The following lemma is the boundary analogue of Lemma \ref{lemma:iter-W1p}.

\begin{lemma}\label{lemma:iter-W1p-bdry}
Suppose that $\partial\Omega\cap B_2$ is $(\delta,2)$-Reifenberg flat, $\|f \|_{L^p(\Omega\cap B_2)} \leq \eta$, $\|D g\|_{L^{\frac{np}{n-p}}(B_{4\sqrt n})} \leq \eta$,  for some $p\in (p_0,n)$ and for some $\eta\in(0,1)$, and $\| u^\e \|_{L^\infty(\Omega\cap B_2)} \leq 1$. Let $L_k^\e$, $B_k$ and $C_k$ denote the sets $L_{m^k}^\e (u^\e ,\Omega\cap B_2)\cap \Omega_\e\cap B_1$, $\{ I_p( |f|^p \chi_{\Omega\cap B_2}) > \eta^p m^{kp}\}$ and respectively $\{ M(|D g|\chi_{B_2}) > \eta^{\frac{np}{n-p}} M^{k\frac{np}{n-p}}\}$. Then 
$$
|L_{k+1}^\e| \leq c(\delta^{\frac{\alpha}{2}} + \eta)^\mu | L_k^\e \cup B_k \cup C_k |,
$$
for any integer $k\geq 1$, where $\alpha\in(0,1)$, $\delta>0$, $\mu>0$, $c>1$ and $m>1$ depend at most on $n$, $\lambda$, $\Lambda$ and $p$. 
\end{lemma}

\begin{proof}
The proof resembles with that of Lemma \ref{lemma:iter-W1d-bdry}. A key difference here is that now we need to be careful of the changes made in boundary data, when rescaling the problem; note that we did not encounter this issue in the proof of Lemma \ref{lemma:iter-W1d-bdry}, since we did not need to see the boundary value at all. Henceforth, we shall proceed with the proof focusing on this issue, and try to skip any argument that only requires a minor modification of what is already shown so far. 

Fix an integer $k\geq 1$. Arguing as in the proof of Lemma \ref{lemma:iter-W1p}, one may use Lemma \ref{lemma:apprx-W1p-bdry}, in place of Lemma \ref{lemma:apprx-W1p}, to deduce that $|L_{k+1}^\e | \leq \eta_0 |B_1|$, with $\eta_0 = c_0(\delta^{\frac{\alpha}{2}} + \eta)^\mu$, where $c_0 > 1$ and $m>1$ are to be chosen later. 

Let $B\subset B_1$ be any ball with center in $\overline\Omega\cap B_2$ and $\rad (B)\leq 1$. Suppose that $|L_{k+1}^\e \cap B| > \eta_0 |B|$. As in the proof of Lemma \ref{lemma:iter-W1d-bdry}, we assert that $\Omega \cap B\subset (L_k^\e \cup B_k\cup C_k)$. If the claim is true, then $|L_{k+1}^\e | \leq c_n\eta_0 | L_k^\e \cup B_k \cup C_k|$ as desired, by Lemma \ref{lemma:vitali}. 

Assume by way of contradiction that $\Omega\cap B\setminus (L_k^\e\cup B_k\cup C_k) \neq\emptyset$. For the same reason as in the proof of Lemma \ref{lemma:apprx-W1p}, we have $2r_B> \e$, by taking $m>1$ large, depending only on $n$. Choose a point $x_B\in\partial\Omega\cap B_1$ such that $2B \subset B_{4r_B}(x_B)$. Select any $\tilde x_B\in \Omega\cap B \setminus (L_k^\e \cup B_k\cup C_k)$. Since $\tilde x_B\not\in \Omega\cap B_{2r_B}(x_B)\setminus L_k^\e$ and $2r_B>\e$, one can find some constant $a$ for which 
\begin{equation}\label{eq:iter-W1p-bdry-1}
|u_B(x) - a | \leq m^k ( 4r_B + |x-x_B|),
\end{equation}
for all $x\in \Omega\cap B_2$. Moreover, it follows from $\tilde x_B\in \Omega\cap B_{2r_B}(x_B)\setminus (B_k\cup C_k)$, one may deduce as in the proof of Lemma \ref{lemma:iter-W1p} that
\begin{equation}\label{eq:iter-W1p-bdry-2}
\begin{aligned}
\int_{\Omega \cap B_{2r_B}(x_B)} |f|^p\,dx \leq c r_B^{n - p} \eta^p, \quad\int_{B_{2r_B}(x_B)} |Dg|^{\frac{np}{n-p}}\,dx \leq cr_B^n \eta^{\frac{np}{n-p}}.
\end{aligned}
\end{equation}

Consider the following rescaled versions of $u^\e$, $f$ and $g$,
$$
\begin{aligned}
&u_B^{\e_B} (x) = \frac{u^\e (x_B + 2r_B x)- a}{2cm^k r_B} ,\quad \e_B = \frac{\e}{2r_B}, \\
&g_B(x)  = \frac{g(x_B + 2r_B x) - a}{2cm^k r_B},\quad f_B(x) = \frac{2r_B}{cm^k} f(x_B + 2r_Bx), 
\end{aligned}
$$
where $c>1$ is a constant to be determined later. In view of \eqref{eq:main-bdry-W1p}, we may compute that 
$$
\begin{dcases}
F_B \left( D^2 u_B^{\e_B}, \frac{\cdot}{\e_B} \right) = f_B& \text{in }\Omega_B \cap B_2,\\
u_B^{\e_B} = g_B  &\text{on }\partial\Omega_B\cap B_2,
\end{dcases}
$$
in the viscosity sense, where $\Omega_B = \frac{1}{2r_B} (-x_B + \Omega\cap B_2)$ and $F_B(P,y) = \frac{2r_B}{cm^k}F(\frac{cm^k}{2r_B} P , y + \frac{x_B}{\e})$. Note that $F_B\in C(\cS^n\times\R^n)$ and it satisfies \eqref{eq:F-ellip} -- \eqref{eq:F-0} for an obvious reason. Selecting a large $c>1$, we observe from \eqref{eq:iter-W1p-bdry-1} and \eqref{eq:iter-W1p-bdry-2} that $u_B$, $f_B$ and $g_B$ verify the hypothesis of Lemma \ref{lemma:apprx-W1p-bdry}. Due to the scaling invariance of the Reifenberg flatness, $\partial\Omega_B\cap B_2$ is also $(\delta,2)$-Reifenberg flat. Hence, all the hypotheses of Lemma \ref{lemma:apprx-W1p-bdry} are verified, from which it follows that
$$
|L_s^{\e_B} (u_B^{\e_B}  , \Omega_B)\cap \Omega_{B,\e_B}\cap B_1| \leq c(\delta^{\frac{\alpha}{2}} + \eta)^\mu s^{-\mu}, 
$$
for any $s>N$, for some $N>1$ depending only on $n$, $\lambda$ and $\Lambda$. To this end, we may follow the argument at the end of the proof of Lemma \ref{lemma:iter-W1p} to arrive at $| L_{k+1} \cap B| \leq \eta_0 |B|$, with suitable choice of $\eta_0$, a contradiction. This finishes the proof. 
\end{proof}

The uniform boundary $W^{1,\frac{np}{n-p}}$-estimates can now be proved as follows.

\begin{proof}[Proof of Theorem \ref{theorem:bdry-W1p}]
One may argue exactly as in the proof of Theorem \ref{theorem:int-W1p}, by substituting Lemma \ref{lemma:iter-W1p} with Lemma \ref{lemma:iter-W1p-bdry}. The additional term, namely the measure of $C_k$ in the notation of the latter lemma, is controlled by the $W^{1,\frac{np}{n-p}}$-regularity of the boundary data $g$, as well as the strong $(q,q)$-type inequality, with $q>1$, for the maximal function. We omit the detail to avoid repeating arguments.
\end{proof}


\section{Uniform $W^{2,p}$-Estimates}\label{section:W2p}


\subsection{Estimates in the Interior}\label{section:int-W2p}

This section is devoted to interior $W^{2,p}$-estimates for viscosity solutions to a certain class of fully nonlinear homogenization problems.

\begin{theorem}\label{theorem:int-W2p}
Let $F\in C(\cS^n\times\R^n)$ be a functional satisfying \eqref{eq:F-ellip} -- \eqref{eq:Fb-W2VMO}, $\Omega\subset\R^n$ be a bounded domain, $f\in C(\Omega)\cap L^p(\Omega)$ for some $p\in(p_0,\infty)$, and $u^\e\in C(\Omega)$ be a viscosity solution to 
\begin{equation}\label{eq:main-int-W2p}
F\left( D^2 u^\e,\frac{\cdot}{\e}\right) = f\quad\text{in }\Omega. 
\end{equation}
Then $H_\Omega^\e (u^\e)\in L_{loc}^p(\Omega)$ and for any subdomain $\Omega'\Subset\Omega$, 
$$
\|  H_\Omega^\e (u^\e)  \|_{L^p(\Omega')} \leq C\left( \frac{\| u^\e \|_{L^\infty(\Omega)}}{\dist(\Omega',\partial\Omega)^{2-\frac{n}{p}}} + \| f \|_{L^p(\Omega)} \right),
$$
where $C>0$ depends only on $n$, $\lambda$, $\Lambda$, $\psi$, $\kappa$ and $p$. Suppose further that there exists some exponent $\bar p\in(p,\infty)$ for which $F(D^2 v,y_0) = a$ has an interior $W^{2,\bar p}$-estimate with constant $\kappa$,\footnote{That is, given any $R>0$ and any $v_0\in C(\partial B_R)$, one can find a viscosity solution $v\in C(\overline{B_R})\cap W_{loc}^{2,\bar p}(B_R)$ to $F(D^2 v, y_0) = 0$ in $B_R$ and $v = v_0$ on $\partial B_R$ such that $$\| D^2 v \|_{L^{\bar p}(B_r)} \leq \kappa(R-r)^{\frac{n}{\bar p} - 2} \| v_0 \|_{L^\infty(\partial B_R)}.$$} uniformly for all $y_0\in\R^n$, and that $\frac{1}{|B|\psi(\rad B)}\int_B |F(P,y) - (F)_B(P)|\leq \kappa|P|$ for all balls $B\subset Q_2$ and any $P\in\cS^n$. Then the same assertion holds with $H_\Omega(u^\e)$, hence with $|D^2 u^\e|$, in place of $H_\Omega^\e( u^\e)$, in which case the constant $C$ may depend also on $\bar p$. 
\end{theorem}

Although Lemma \ref{lemma:apprx} yields an error estimate between $u^\e$ and $\bar u$ in $L^\infty$ norm, we cannot expect $\bar u$ is close to $u^\e$ in the viscosity sense (i.e., $\cP_\pm (D^2 (u^\e - \bar u)) \neq o(1)$ with $\cP_\pm$ being the Pucci extremal operators), since $D^2 u^\e$ is supposed to be rapidly oscillating around $D^2\bar u$ in the small scales. 

In the next lemma, we obtain the closeness between $D^2 u^\e$ and $D^2 \bar u$ in the viscosity sense by incorporating interior correctors. To do so, we shall assume $VMO$-condition, or more exactly $BMO_\psi$-condition for some modulus of continuity $\psi$, for $D^2 \bar u$. This condition replaces the small oscillation condition for standard problems, $F(D^2 u, x) = f$, c.f. \cite[Theorem 7.1]{CC}. Let us remark that $BMO_\psi$-regularity does neither imply nor follow from the boundedness of $D^2 \bar u$, and that it also allows $D^2 \bar u$ to be discontinuous.

With suitable rescaling argument, it suffices to take care of the case where $\Omega = B_{4\sqrt n}$ and $\Omega' = Q_1$. In what follows, we shall always assume that $F$ is a continuous functional satisfying \eqref{eq:F-ellip} -- \eqref{eq:F-0}, $f\in C(B_{4\sqrt n})\cap L^p(B_{4\sqrt n})$ with $p>p_0$, $u^\e\in C(B_{4\sqrt n})$ is a viscosity solution to \eqref{eq:main-int-W2p}, unless stated otherwise. Moreover, we shall let $c$ denote a positive generic constant depending at most on a set of fixed quantities, shown in the statement of each lemmas below, and we allow it to vary at each occurrence. 

\begin{lemma}\label{lemma:apprx-msr}
Suppose that $\Omega\subset\R^n$ is a bounded domain with $B_{4\sqrt n}\subset \Omega$, $\|u^\e \|_{L^\infty(B_{4\sqrt n})} \leq 1$ and $|u^\e(x) | \leq |x|^2$ in $\Omega\setminus B_{2\sqrt n}$. Let $\eta>0$, $s>1$, $p\in(p_0,\infty)$ and $q\in(p,\infty)$ be given constants. Then there exists a constant $\e_{s,\eta} > 0$, depending only on $n$, $\lambda$, $\Lambda$, $\kappa$, $\psi$, $q$, $p$, $\eta$ and $s$, such that if $0<\e < \e_{s,\eta}$ and $\| f\|_{L^{p}(B_{4\sqrt n})} \leq \e_{s,\eta}$, then 
\begin{equation}\label{eq:apprx-msr}
| A_s^\e (u^\e, \Omega)\cap Q_1 | \leq c(\eta s^{-\mu} + s^{-q}),
\end{equation}
where $\mu\in(0,1)$ depends only on $n$, $\lambda$, $\Lambda$, $p$, and $c>1$ depends further on $\psi$, $\kappa$ and $q$.
\end{lemma}





\begin{proof}
Fix $s>1$ and $\eta>0$.  Let $\delta>0$ be a constant to be determined later. Since the hypothesis of Lemma \ref{lemma:apprx} is met, we can find a constant $\e_\delta\in(0,1)$, corresponding to $\delta$, and a function $\bar u\in C(\overline{B_{3\sqrt n}})$ such that $\bar F(D^2 \bar u) = 0$ in $B_{3\sqrt n}$ in the viscosity sense, $\| \bar u\|_{L^\infty(B_{3\sqrt n})} \leq 1$ and $\| u^\e - \bar u \|_{L^\infty(B_{3\sqrt n})} \leq \delta$. By the assumption \eqref{eq:Fb-W2VMO} on $\bar F$, $D^2 \bar u \in BMO_\psi(B_{2\sqrt n})$ and
\begin{equation}\label{eq:D2ub-VMO}
\int_{B_{2\sqrt n}} |D^2 \bar u|\,dx + [D^2 \bar u]_{BMO_\psi(B_{2\sqrt n})} \leq c. 
\end{equation}  
Due to the John-Nirenberg inequality, we have $\| D^2 \bar u\|_{L^q(B_{2\sqrt n})}\leq c$, so by the strong $(q,q)$-type inequality, $\| M( |D^2 \bar u| \chi_{B_{2\sqrt n}}) \|_{L^q(\R^n)} \leq c$, which ensures that 
\begin{equation}\label{eq:D2ub-Lq-max}
|\{M(|D^2 \bar u|\chi_{2\sqrt n}) > s\}| \leq  cs^{-q}. 
\end{equation}
In addition, by \cite[Lemma 2.5]{LZ}, we also have $\|\Theta(\bar u,B_{2\sqrt n})\|_{L^q(B_{2\sqrt n})} \leq c\| D^2 \bar u\|_{L^q(B_{2\sqrt n})} \leq c$ (see Definition \ref{definition:large-grad-Hess} for the definition of $\Theta$), whence it follows from the relation $\{ \Theta(\bar u,B_{2\sqrt n}) > s\} = A_s (\bar u, B_{2\sqrt n})$ that 
\begin{equation}\label{eq:ub-As}
|A_s (\bar u,B_{2\sqrt n})| \leq cs^{-q}. 
\end{equation}

Let $\rho_s\in(0,\frac{1}{4})$ be a constant to be determined later. Our idea is to subdivide $Q_1$ into two groups, say $\cF$ and $\cG$, of dyadic cubes with side-length in between $\rho_s$ and $2\rho_s$ such that $Q\in\cF$ if $|Q\cap \{M(|D^2 \bar u|\chi_{B_{2\sqrt n}}) \leq \frac{s}{c}\}|>0$, and $Q\in\cG$ if $Q\not\in \cF$; here $c>1$ is a constant to be determined by $n$, $\lambda$, $\Lambda$, $\kappa$, $\psi$ and $q$. Thanks to \eqref{eq:D2ub-Lq-max} and \eqref{eq:ub-As} (as well as an obvious fact that $A_s^\e (g,E)\subset A_s(g,E)$ for any $g\in C(E)$) yields that 
\begin{equation}\label{eq:ue-Ase}
\begin{aligned}
&|A_s^\e (u^\e, B_{2\sqrt n})\cap Q_1| &\\
&\leq \sum_{Q\in\cF} |A_s^\e (u^\e, B_{2\sqrt n})\cap Q| + \sum_{Q\in\cG } |Q| \\
&\leq \sum_{Q\in\cF} |A_s^\e (u^\e - \bar u,B_{2\sqrt n})\cap Q| + |A_s(\bar u,B_{2\sqrt n})\cap Q_1| + cs^{-q}  \\
&\leq \sum_{Q\in\cF} |A_s^\e (u^\e - \bar u,B_{2\sqrt n})\cap Q| + cs^{-q}.
\end{aligned}
\end{equation}
Note that we can replace the domain $B_{2\sqrt n}$ in the leftmost side to $\Omega$ by utilizing the assumption that $|u^\e(x)|\leq |x|^2$ for all $x\in\Omega\setminus B_{2\sqrt n}$. Thus, it will be enough to prove that with $\e_{s,\eta}>0$ sufficiently small and $\e<\e_{s,\eta}$, 
\begin{equation}\label{eq:ue-ub-Ase}
|A_s^\e (u^\e - \bar u,B_{2\sqrt n})\cap Q| \leq c\eta s^{-\mu} |Q|,
\end{equation}
for all $Q\in\cF$, where $\mu>0$ is a constant depending only on $n$, $\lambda$ and $\Lambda$. 

To prove \eqref{eq:ue-ub-Ase}, let us fix a cube $Q\in\cF$. By definition, there is some $x_0 \in Q\cap \{M(|D^2 \bar u|\chi_{B_{2\sqrt n}} )\leq \frac{s}{c}\} \neq \emptyset$ and $Q\subset B_{\rho_s\sqrt n}(x_0) \subset B_{2\sqrt n}$. To simplify the exposition, let us write by $B$ the ball $B_{2\rho_s\sqrt n}(x_0)$, by $r_B$ the radius of $B$, i.e., $r_B = 2\rho_s\sqrt n$, and by $P$ the matrix $(D^2 \bar u)_B$, i.e., $P = \frac{1}{|B|}\int_B D^2 \bar u\,dx$. Then by the choice of $x_0$, $|P|\leq \frac{s}{c}$, and by \eqref{eq:D2ub-VMO} along with the John-Nirenberg inequality,
\begin{equation}\label{eq:D2ub-P-Ln}
\frac{1}{|B|} \int_B |D^2 \bar u - P|^n\,dx  \leq c_n\kappa^n\psi(r_B)^n,
\end{equation}
where $c_n>0$ is a constant depending only on $n$. 

Consider an auxiliary function $w:\R^n\to\R$ satisfying
\begin{equation}\label{eq:w-pde}
\begin{dcases}
F (D_y^2 w + P , y ) = \bar F ( P ) \quad\text{in }\R^n,\\
w( y+ k) = w(y)\quad\text{for all }y\in\R^n, k\in\Z^n,\\
w(0) = 0,
\end{dcases}
\end{equation}
in the viscosity sense. According to \cite[Lemma 3.1]{E}, such a periodic viscosity solution exists (in $C^{0,\alpha}(\R^n)$, with $\alpha\in(0,1)$ universal) and unique, and due to \eqref{eq:F-ellip} as well as\eqref{eq:F-0}, it satisfies 
\begin{equation}\label{eq:w-Ca}
\| w \|_{C^{0,\alpha}(\R^n)}  \leq c_0 |P| \leq \frac{c_0s}{c},
\end{equation}
for some constant $c_0>0$ depending only on $n$, $\lambda$ and $\Lambda$.

Consider auxiliary functions $\phi^\e, g : B \to \R$ defined by 
\begin{equation}\label{eq:phie-g}
\phi^\e  = \frac{1}{\kappa\psi(r_B)}\left( u^\e - \bar u - \e^2 w\left(\frac{\cdot}{\e}\right)  \right),\quad g = \frac{r_B^2}{\kappa \psi(r_B)}(f + c|D^2 \bar u - P|)  . 
\end{equation}
Clearly, $\phi\in C(B)$ and $g\in L^n(B)$. By $\| u^\e -\bar u\|_{L^\infty(B_{3\sqrt n})} \leq \delta$  and \eqref{eq:w-Ca} (as well as $\e\in(0,\e_{s,\eta})$ and $B\subset B_{3\sqrt n}$), 
\begin{equation}\label{eq:phie-Linf}
\| \phi^\e \|_{L^\infty(B)} \leq \frac{\delta + c_0c^{-1}s\e_{s,\eta}^2}{\kappa\psi(r_B)}.
\end{equation}
In addition, by $\|f \|_{L^{p}(B_{4\sqrt n})} \leq \e_{s,\eta}$ and \eqref{eq:D2ub-P-Ln}, 
\begin{equation}\label{eq:g-Ln}
\left( \frac{1}{|B|} \int_{B} |g|^n\,dx \right)^{1/n} \leq c_nr_B^2\left( \frac{\e_{s,\eta}}{ r_B \kappa \psi(r_B)}  + 1 \right).
\end{equation}
We claim that $\phi$ is an $L^p$-viscosity solution to
\begin{equation}\label{eq:phi-pde}
\max\{\cP_- (D^2 \phi) , - \cP_+(D^2 \phi) \} \leq |g|\quad\text{in }B. 
\end{equation} 

Suppose for the moment that the claim is true. With sufficiently small $\delta$ and $\e_{s,\eta}$, whose smallness condition depending only on $r_B$, $\kappa$, $\psi(r_B)$, $c_0$ and $c$, we can deduce from \eqref{eq:phie-Linf} and \eqref{eq:g-Ln} that $\|\phi^\e \|_{L^\infty(B)} \leq 1$ and respectively $\| g\|_{L^n(B)} \leq c_n r_B$. Then we can apply (a rescaled form of) \cite[Lemma 7.5]{CC}\footnote{Although this lemma is written for $C$-viscosity solutions and continuous right hand side whose $L^n$-norm is under control, it works equally well for $L^n$-viscosity solutions and measurable right hand side.} to $\phi^\e$ and deduce that for any $t>0$, 
\begin{equation}\label{eq:phie-msr}
|A_t(\phi^\e, B) \cap Q| \leq c_0t^{-\mu} |Q|,
\end{equation}
where $\mu>0$ depends only on $n$, $\lambda$ and $\Lambda$; here we also used that $\dist(Q,\partial B) \geq  \frac{1}{2} r_B$, which is apparent from the choice of $Q$ and $B$. Thus, setting $t = \frac{s}{2\kappa \psi(r_B)}$, we obtain $|A_{\frac{s}{2}}^\e (u^\e - \bar u - \e^2 w(\frac{\cdot}{\e})), B))\cap Q| \leq c_0 (\kappa\psi(r_B))^\mu s^{-\mu} |Q|$. Utilizing \eqref{eq:w-Ca} (as well as a simple observation that $A_{(2\ell + 1)s} (f_1 - f_2,E) \subset A_s (f_1,E)$ for any $f_1,f_2\in C(E)$ with $\|f_2 \|_{L^\infty(E)} \leq \ell \e^2$), we deduce that 
\begin{equation}\label{eq:ue-ub-Ase-re}
\left|A_{(\frac{c_0}{c} + \frac{1}{2})s}^\e (u^\e - \bar u, B)\cap Q \right| \leq c_0\kappa^\mu \psi(r_B)^\mu s^{-\mu} |Q|,
\end{equation}
for any $s>0$. Since $\| u^\e - \bar u \|_{L^\infty(B_{3\sqrt n})} \leq \delta$ with $\delta$ being small depending on $r_B$, one can also replace $B$ above with $B_{2\sqrt n}$. At this point, we select $c \geq 2c_0$, and $\rho_s = \frac{r_B}{2\sqrt n}\in (0,\frac{1}{4})$ as a small constant such that $\kappa\psi(r_B) \leq \eta^{1/\mu}$, so that we arrive at \eqref{eq:ue-ub-Ase}, as desired.  

Thus, we are only left with proving that \eqref{eq:phi-pde}  holds in the $L^{p}$-viscosity sense. Let $\vp$ be a quadratic polynomial such that $D^2 \vp = P$. For the moment, let us denote by $W^\e$ the function $\e^2 (\vp + w)(\frac{\cdot}{\e})$. Clearly, $F(D^2 W^\e,\frac{\cdot}{\e}) = \bar F(P)$, as well as $F(D^2 u^\e ,\frac{\cdot}{\e}) = f$, in $B$ in the viscosity sense, so it follows from Lemma \ref{lemma:visc} that 
$$
\cP_- (D^2 (u^\e - W^\e)) \leq f - \bar F(P)\quad\text{in } B,
$$
in the viscosity sense. Now since $u^\e - \bar u - \e^2 w(\frac{\cdot}{\e}) = u^\e - W^\e - (\bar u - \e^2 \vp(\frac{\cdot}{\e}))$, and $\bar u - \e^2 \vp(\frac{\cdot}{\e}) \in W^{2,p}(B)$, one can deduce, along with $\bar F(D^2 \bar u) = 0$ in $B$ in the $L^{p}$-strong sense, that
$$
\begin{aligned}
\cP_- (D^2 \phi^\e ) &\leq \frac{1}{\kappa\psi(r_B)} ( \cP_-(D^2 (u^\e -W^\e)) - \cP_- (D^2 \bar u - P))  \\
&\leq \frac{1}{\kappa\psi(r_B)}  ( f - \bar F(P) - \cP_- (D^2 \bar u - P) ) \\
&\leq \frac{1}{\kappa\psi(r_B)} ( f - \cP_- (P - D^2 \bar u) - \cP- (D^2 \bar u - P)) \\
& \leq |g|
\end{aligned}
$$
in $B$ in the $L^{p}$-viscosity sense. Similarly, we obtain $\cP_+(D^2 \phi^\e) \geq - |g|$ in $B$ in the $L^{p}$-viscosity sense. This finishes the proof. 
\end{proof}

We are ready to use the by-now standard cube decomposition argument to obtain a geometric decay of $|A_s^\e (u^\e,\Omega)\cap Q_1|$. The idea is the same with \cite{CP} in the sense that we split the set $A_k = A_{m^k}^\e (u^\e,\Omega)\cap Q_1$ into two parts, say $D_k$ and $E_k$, where $D_k$ is the part of $A_k$ intersected by its Calder\'on-Zygmund cube covering whose side-length is at least $\frac{\e}{\e_\eta}$ (these cubes are said to be of high frequency), whereas $E_k = A_k \setminus D_k$, i.e., the part of $A_k$ intersected by the cubes of low frequency. As for $D_k$, we can deduce a geometric decay via Lemma \ref{lemma:apprx-msr}, and this part is almost the same with the argument for standard problem, e.g., \cite[Lemma 7.13]{CC}. As for $E_k$, the above lemma is no longer applicable, but at the same time we cannot argue as in \cite{CP} because we do not assume any structure condition on $F$ so as to ensure sufficient regularity of $u^\e$ in small scales. Here we control $E_k$ directly from the fact that the set $A_k$ (or more exactly $Q_1\setminus A_k$) allows error of order $\e^2$ for quadratic polynomials to touch $u^\e$.


\begin{lemma}\label{lemma:iter-msr}
Let $\eta>0$ and $p\in(p_0,\infty)$ be given. Assume that $0<\e<\e_\eta$, $\| f\|_{L^{p}(B_{4\sqrt n})} \leq \e_\eta$, and $\| u^\e \|_{L^\infty(B_{4\sqrt n})} \leq 1$. Let $A_k$ and $B_k$ denote the sets $A_{m^k}^\e(u^\e,B_{4\sqrt n})\cap Q_1$ and respectively $\{ M(|f|^{p}\chi_{B_{4\sqrt n}} )> \e_\eta^{p} M^{ pk}\} $. Then one has, for each $k\geq k_\eta$,
\begin{equation}\label{eq:iter-msr}
|A_{k+1}| \leq (\eta + cm^{-q}) | A_k \cup B_k | + |B_{k-k_\eta}|,
\end{equation}
where $m>1$ depends at most on $n$, $\lambda$, $\Lambda$, $\psi$, $\kappa$, $p$ and $q$, while $\e_\eta>0$ and $k_\eta>1$ may depend further on $\eta$. 
\end{lemma}

\begin{proof}
As briefly mentioned in the discussion before the statement of this lemma, the set $D_k$ is the part of $A_k$ intersected by its Calder\'on-Zygmund covering whose side-length is no less than $\frac{\e}{\e_\eta}$. More exactly, we choose $M>1$ sufficiently large such that $|A_1| \leq \eta + cm^{-q} < 1$ due to Lemma \ref{lemma:apprx-msr}. We shall fix $\e_\eta$ as the constant $\e_{c_0^{-1}m,\eta}$ from Lemma \ref{lemma:apprx-msr}, with $c_0>1$ to be determined later. As $A_{k+1}\subset A_k\subset\cdots A_1$ for each integer $k\geq 1$, $|A_{k+1}|\leq \eta + cm^{-q}$, whence there exists a Calder\'on-Zygmund covering, denoted by $\cF_{k+1}^\e$, of $A_{k+1}$ corresponding to the level $\eta + cm^{-q}$. Define
$$
D_{k+1} =\bigcup \left\{A_{k+1}\cap Q : Q\in\cF_{k+1}^\e,\, \frac{\diam(Q)}{\sqrt n}  > \frac{\e}{\e_\eta}\right\}. 
$$ 
 
Let us claim that $\tilde Q\subset A_k\cup B_k$ for any cube $Q\in \cF_k^\e$ whose side-length is no less than $\frac{\e}{\e_\eta}$, where $\tilde Q$ is the predecessor of $Q$. Once this is proved, from the fact that $Q$ belongs to the Calder\'on-Zygmund covering of $A_k$ it follows immediately that 
\begin{equation}\label{eq:iter-msr-claim}
|D_{k+1} |\leq (\eta + cm^{-q}) |A_k\cup B_k|. 
\end{equation} 

The proof for the above claim mostly follows the argument for standard problems, e.g., \cite[Lemma 7.12]{CC}, except for the following two points: (i) we need to verify the hypothesis \eqref{eq:Fb-W2VMO} for the effective functional at each iteration step, (ii) the set $A_k$ involves error of order $\e^2$ and may vary along with different scalings in the domain.

Denote by $x_Q$ and $s_Q$ the center and respectively the side-length of $Q$. Since $Q \in \cF_k^\e$, $|A_{k+1}\cap Q| > (\eta + cm^{-q})|Q|$. Assume, by way of contradiction, that $\tilde Q\setminus (A_k\cup B_k) \neq \emptyset$. Set $\Omega_Q = s_Q^{-1} (-x_Q + B_{4\sqrt n})$, the image of $B_{4\sqrt n}$ via the above rescaling. Then due to the assumption $\tilde Q \setminus (A_k\cup B_k)\neq\emptyset$, the functions $u_Q^{\e_Q}\in C(\Omega_Q)$, defined by 
\begin{equation}\label{eq:uiei-fi}
u_Q^{\e_Q} (x) = \frac{ (u^\e - \ell) (x_Q + s_Q x)}{c_0m^k s_Q^2},\quad \e_Q = \frac{\e}{s_Q},
\end{equation}
for certain linear polynomial $\ell$, satisfy that $\| u_Q^{\e_Q}\|_{L^\infty(\Omega_Q)} \leq 1$, $| u_Q^{\e_Q}(x)|\leq |x|^2$ for all $x\in\Omega\setminus B_{2\sqrt n}$, provided that we choose $c_0>1$ to be large (depending only on $n$).  Now as $F(D^2 u^\e,\frac{\cdot}{\e}) = f$ in $B_{4\sqrt n}$ in the viscosity sense, we have
\begin{equation}\label{eq:uiei-pde}
F_Q \left( D^2 u_Q^{\e_Q},\frac{\cdot}{\e_Q}\right) = f_Q\quad\text{in }B_{4\sqrt n}, 
\end{equation}
in the viscosity sense, where $F_Q (P,y) = \frac{1}{c_0m^k} F ( c_0m^k P, \frac{x_0}{\e} + y)$ and $f_Q(x) = \frac{1}{c_0m^k} f( x_0 + s_Q x)$. Clearly, $F_Q \in C(\cS^n\times\R^n)$ satisfies \eqref{eq:F-ellip}, \eqref{eq:F-peri} and \eqref{eq:F-0}. Noting that its effective functional $\bar F_Q : \cS^n\to\R$ is given by $\bar F_Q (P) = \frac{1}{c_0m^k} \bar F(c_0m^k P)$, for any $P\in\cS^n$, $\bar H$ satisfies \eqref{eq:Fb-W2VMO} (with the same modulus of continuity $\psi$ and constant $\kappa$ as $\bar F$ does). Moreover, with $c_0>1$, it follows from the choice of $\tilde Q\setminus B_k\neq\emptyset$ that $\|f_Q \|_{L^\infty(B_{4\sqrt n})} \leq c_0^{-\frac{1}{p}} \e_\eta < \e_\eta$. Therefore, due to the choice of $\e_Q = s_Q^{-1}\e \leq \e_\eta = \e_{c_0^{-1}m,\eta}$, we can apply Lemma \ref{lemma:apprx-msr} to $u_Q^{\e_Q}$ and deduce that 
$$
\left|A_{\frac{m}{c_0}}^{\e_Q}(u_Q^{\e_Q},\Omega_Q)\cap Q_1\right| \leq \eta + cm^{-q}.
$$
Rescaling back, we arrive at $|A_{k+1} \cap Q|\leq (\eta + cm^{-q})|Q|$, a contradiction to $Q\in\cF_k^\e$. This finishes the proof for \eqref{eq:iter-msr-claim}. 

Next, we prove that 
\begin{equation}\label{eq:iter-msr-claim2}
A_{k+1} \setminus D_{k+1} \subset B_{k-k_\eta},
\end{equation}
with $k_\eta = -\log_2 \e_\eta$. Let $Q\in\cF_k^\e$ be a dyadic cube with side-length less than $\frac{\e}{\e_\eta}$. Then $|Q\cap A_{k+1}| > (\eta + cm^{-q})|Q|$, so $|Q\cap A_{k-k_\eta}^\e| > (\eta+cm^{-q})|Q|$. We claim that $\tilde Q\setminus A_{k-k_\eta}^\e \subset B_{k-k_\eta}$, from which \eqref{eq:iter-msr-claim2} follows immediately. Suppose, by way of contradiction, that $\tilde Q \setminus (A_{k-k_\eta}^\e\cup B_{k-k_\eta}^\e)\neq\emptyset$. As we have chosen $k_\eta = -\log_2 \e_\eta$, that $\tilde Q \setminus A_{k-k_\eta}^\e \neq\emptyset$ implies $\tilde Q\cap A_{k+1} = \emptyset$. This is again a contradiction to the fact that $Q\subset \tilde Q$ and $|Q\cap A_{k+1}| > (\eta + cm^{-q})|Q|$. Therefore, $\tilde Q\setminus A_{k-k_\eta}\subset B_{k-k_\eta}$ as desired, and the inclusion in \eqref{eq:iter-msr-claim2} follows. 
\end{proof}

We are now ready to prove the first assertion of Theorem \ref{theorem:int-W2p}.

\begin{proof}[Proof of Theorem \ref{theorem:int-W2p}; the first part]
Fix $p \in (p_0,\infty)$. With suitable rescaling of the problem, it is sufficient to prove the assertions for the case where $\Omega = B_{4\sqrt n}$, $\Omega' = Q_1$, $\|u^\e \|_{L^\infty(B_{4\sqrt n})}\leq 1$ and $\| f \|_{L^p(B_{4\sqrt n})} \leq \e_\eta$, with $\e_\eta$ to be determined. 

Choose $q = 2p$ and $p' = \frac{p + p_0}{2}$. Let $M>1$ be as in Lemma \ref{lemma:iter-msr} with $p$ replaced by $p'$. We can assume that $cM^{p-q} \leq \frac{1}{4}$ by taking $M$ larger, depending on the choice of $p$. Then we select $\eta>0$ as a small constant satisfying $M^p \eta \leq\frac{1}{4}$. Then 
$$
M^p (\eta + cm^{-q}) \leq \frac{1}{2}.
$$ 
Now let $\e_\eta > 0$ and $k_\eta >1$ be as in Lemma \ref{lemma:iter-msr} corresponding to the specific choice of $\eta$. Note that all the constants $\eta$, $M$, $\hat c $, $\e_\eta$ and $k_\eta$ involved in the statement of Lemma \ref{lemma:iter-msr} now depend only on $n$, $\lambda$, $\Lambda$, $\psi$, $\kappa$ and $p$ (as $q$ depending solely on $p$, and $\gamma$ solely on $n$, $\lambda$, $\Lambda$ and $p$).

With \eqref{eq:iter-msr} at hand, we deduce that 
$$
\sum_{k=1}^\infty m^{kp}\alpha_k \leq c + c \sum_{k=k_\eta}^\infty m^{kp} \beta_k,
$$
with $\alpha_k = |A_k|$ and $\beta_k = |B_k|$, where $A_k$ and $B_k$ are as in the statement of Lemma \ref{lemma:iter-msr}. Since $\| |f|^{p'}  \|_{L^{p/ p'}(Q_1)} \leq \| f\|_{L^p(B_{4\sqrt n})}^{p'} \leq \e_\eta^{p'}$, we have $\sum_{k=1}^\infty m^{kp}\beta_k \leq c$, as well. The proof can now be finished, as in \cite[Proposition 7.2]{CC}, and we omit the detail. 
\end{proof}

To prove the second assertion of Theorem \ref{theorem:int-W2p}, we shall modify some of the argument for the proof of Lemma \ref{lemma:apprx-msr} and \ref{lemma:iter-msr} in such a way that we estimate the measure of $A_s(u^\e,\cdot\,)$ instead of $A_s^\e(u^\e,\cdot\,)$. In what follows, we intend to highlight those changes, which are not so trivial, and then omit the detail for the argument that might repeat what is already written so far. 

\begin{proof}[Proof of Theorem \ref{theorem:int-W2p}; the second part]
{\it Step 1:} As the first step, let us replace $A_s^\e(u^\e,\cdot\,)$ with $A_s(u^\e,\cdot\,)$ in the approximation lemma,  i.e., Lemma \ref{lemma:apprx-msr}, and claim that under the additional assumption in the theorem, with $\hat p \in (p_0, \bar p)$ (where $\bar p$ is the exponent for which \eqref{eq:F-W2q} holds), 
\begin{equation}\label{eq:apprx-msr-re}
|A_t(u^\e,\Omega_{4\sqrt n}) \cap Q_1|\leq c \left( \eta t^{-\mu} + \left(\frac{s}{t}\right)^{\hat p} + s^{-q}\right),
\end{equation}
for any $t>s$, for any $\e\in(0,\e_{s,\eta})$. Here $\e_{s,\eta}$ and $c$ will depend on $\hat p$ in addition to the parameters that are specified in the statement of the lemma; in particular, $\e_{s,\eta}$ is independent with the choice of $t>s$.

The additional assumption in Theorem \ref{theorem:int-W2p} (ii) allows us to estimate the corrector function $w$, given as the periodic viscosity solution to \eqref{eq:w-pde}, in space $W^{2,\hat p}(Q_1)$; here $Q_1$ is the period of $w$. In view of \eqref{eq:w-Ca}, we obtain (along with the periodicity of $w$) that $w\in W^{2,\hat p}(Q_1)$ with 
$$
\| D^2 w\|_{L^{\hat p}(Q_1)} \leq cs. 
$$
Here $c$ may depend on $\hat p$, but it is independent of $s$. This can be deduced from the fact that $s^{-1} w$ is a viscosity solution to $\frac{1}{s} F( s (D^2 v + s^{-1}P), y ) = \frac{1}{s} \bar F(P)$ in $Q_2$ with $\| s^{-1} w \|_{L^\infty(Q_2)}\leq c_0$ and $|s^{-1} \bar F(P)| \leq c_0$ (with $c_0>0$ depending only on $n$, $\lambda$ and $\Lambda$), whence $s^{-1} w \in W^{2,\hat p}(Q_1)$ with $\| s^{-1} D^2 w\|_{L^{\hat p}(Q_1)} \leq c$, for some constant $c>0$ depending only on $n$, $\lambda$, $\Lambda$, $\kappa$, $\psi$ and $\hat p$. Thus, it follows from \cite[Lemma 2.5]{LZ} that $\| \Theta(w,Q_1) \|_{L^{\hat p}(Q_1)} \leq cs$, which along with $\{\Theta(w,Q_1)> t\} = A_t(w,Q_1)$ implies that 
\begin{equation}\label{eq:w-At}
|A_t( w,Q_1)|\leq c\left(\frac{s}{t}\right)^{\hat p}. 
\end{equation}

On the other hand, arguing as in the proof of Lemma \ref{lemma:apprx-msr}, one may observe that the auxiliary function $\phi^\e$ as in \eqref{eq:phie-g} continues to satisfy \eqref{eq:phie-msr}. Combining \eqref{eq:phie-msr} with \eqref{eq:w-At}, and employing the periodicity of $w$, we can deduce that 
$$
\begin{aligned}
|A_t(u^\e - \bar u, B)\cap Q| &\leq \left| A_t \left(u^\e -\bar u - \e^2 w\left(\frac{\cdot}{\e}\right), B\right)\cap Q\right| + |A_t (w,Q_1)| |Q| \\
& \leq \left( c_0 \eta t^{-\mu} + c \left(\frac{s}{t}\right)^{\hat p} \right) |Q|,
\end{aligned}
$$ 
with $B$, $r_B$ and $Q$ as in the proof of the lemma. Finally, noting that we can substitute $A_s^\e$ in \eqref{eq:ue-Ase} with $A_t$ for any $t>s$, without any further change in the argument, we arrive at \eqref{eq:apprx-msr-re} as desired. 

\noindent
{\it Step 2:} As the second step, we modify the iteration scheme, Lemma \ref{lemma:iter-msr}. Putting $A_{m^k} (u^\e,B_{4\sqrt n})$ in place of $A_{m^k}^\e(u^\e,B_{4\sqrt n})$ and keeping other notation the same (so now $A_k = |A_{m^k}(u^\e,B_{4\sqrt n})\cap Q_1|$ and $\alpha_k = |A_k|$), we claim (instead of \eqref{eq:iter-msr}) that
\begin{equation}\label{eq:iter-msr-re}
\alpha_k \leq \left( \eta + \frac{c}{N^{\hat p}} + \frac{c N^q}{M^q} \right) (\alpha_{k-1}^\e + \beta_{k-1}) + \beta_{k-k_\eta},
\end{equation}
where $N,M>1$ are some constants depending on $\hat p$ in addition to the quantities involved in the statement of the lemma; yet they are independent on $\eta$. Here one may consider $M = t$ and $N = \frac{t}{s}$ in the notation of Step 1. 

The choice of $M$ and $N$ are quite flexible here; we only need to allow them two satisfy $cM^{-\mu} \leq 1$ and $cN^{-\hat p} + cN^qm^{-q} \leq \frac{1}{4}$ (so that the parenthesis in \eqref{eq:iter-msr-re} is strictly less than $1$ and we can invoke the Calder\'on-Zygmund cube decomposition). 

The first part of the proof of Lemma \ref{lemma:iter-msr} is the same: we prove \eqref{eq:iter-msr-claim}, by utilizing \eqref{eq:apprx-msr-re} in place of  Lemma \ref{lemma:apprx-msr}, where the set $D_k$ is accordingly defined, and the factor $(\eta + cm^{-q})$ is replaced by $(\eta + cN^{-\hat p} + cN^qm^{-q})$. 

As for the second part of the proof, we use the additional assumption in Theorem \ref{theorem:int-W2p} (ii) to prove \eqref{eq:iter-msr-claim2}; let us mention that this step is reminiscent of \cite[Lemma 4.3]{CP}. Let $Q\in\cF_k^\e$ be a cube with side-length less than $\frac{\e}{\e_\eta}$, and assume by way of contradiction that $\tilde Q\setminus (A_{k-k_\eta}\cup  B_{k-k_\eta})\neq\emptyset$. Then we can define $u_Q^{\e_Q}$, $f_Q$, $\Omega_Q$ and $H$ as in the first part of the proof of Lemma \ref{lemma:iter-msr} (see \eqref{eq:uiei-fi} and \eqref{eq:uiei-pde}), with $k$ replaced by $k-k_\eta$, such that we again obtain $F_Q (D^2 u_Q^{\e_Q}, \frac{\cdot}{\e_Q}) = f_Q$ in $B_{4\sqrt n}$, $\| u_Q^{\e_Q} \|_{L^\infty(B_{4\sqrt n})} \leq 1$, $|u_Q^{\e_Q}(x)| \leq |x|^2$ in $\Omega_Q \setminus B_{2\sqrt n}$ and $\|f_Q\|_{L^{p}(B_{4\sqrt n})} \leq \e_\eta$. Note also that $F_Q$ satisfies the additional assumption in the theorem. Here $\e_Q = s_Q^{-1}\e \geq \e_\eta$, where $2^{-i}$ is the side-length of $Q$. This along with $\e_Q\geq \e_\eta$ implies that $\| D^2 u_Q^{\e_Q}\|_{L^{p}(Q_1)} \leq c_\eta$. By \cite[Lemma 2.5]{LZ} (as well as $|u_Q^{\e_Q}(x)| \leq |x|^2$ in $\Omega_Q\setminus B_{2\sqrt n}$), $|A_{M^{k_\eta}} (u_Q^{\e_Q}, \Omega_Q)\cap Q_1| \leq c_\eta M^{-k_\eta p}$. Hence, choosing $k_\eta>1$ sufficiently large such that $c_\eta M^{-k_\eta p} \leq \frac{1}{2} (\eta + cN^{-\hat p} + cm^{-q})$, we arrive at $|A_k\cap Q|\leq \frac{1}{2} (\eta + cN^{-\hat p} + cN^qm^{-q}) |Q|$, a contradiction to $Q\in \cF_k^\e$. This finishes the proof for \eqref{eq:iter-msr-claim2}, and hence \eqref{eq:iter-msr-re} follows. 

To this end, let us choose $\hat p$, $q$, $N$ and $M$ carefully, yet depending only on $p$, $\bar p$ and the universal quantities, such that $\hat p > p$, $q \gg \hat p$ and $1 \ll N \ll M$, and with $\eta \leq 1/(4M^p)$,
$$
M^p \left(\eta + \frac{c}{N^{\hat p}} + \frac{cN^q}{M^q} \right) \leq \frac{1}{2}.
$$
Choose $\hat p = \frac{p + 3\bar p}{4} > \hat{\hat p} = \frac{3p + \bar p}{4}$, and $q = \frac{2p\hat p}{\hat p - \hat{\hat p}} = \frac{4p\hat p}{\bar p- p}$. Next we ask $M>1$, among other conditions determined from the argument above, to satisfy that $cM^{p - \hat{\hat p}} \leq \frac{1}{8}$ and $cM^{-p} \leq \frac{1}{8}$, and then set $N = M^{1-\frac{2p}{q}} = M^{\hat{\hat p}/\hat p} \ll M$. Then $M^p (cN^{-p} + cN^qm^{-q})  \leq c(M^{p-\hat{\hat p}} + M^{-p} )\leq \frac{1}{4}$. Clearly, such a choice of $M$ and $N$ depend only on $p$ and $\bar p$, in addition to the quantities $n$, $\lambda$, $\Lambda$, $\kappa$ and $\psi$. Consequently, we can iterate \eqref{eq:iter-msr-re} as in the proof of Theorem \ref{theorem:int-W2p} (i) and derive that $\sum_{k=1}^\infty m^{kp} \alpha_k \leq c$. This finishes the proof. 
\end{proof} 


\subsection{Estimates near Boundary Layers}\label{section:bdry-W2p}

This section is devoted to the study of uniform boundary $W^{2,p}$-estimates for viscosity solutions to fully nonlinear homogenization problems, where the boundary is an $(n-1)$-dimensional $C^2$-manifold. In what follows, we shall denote by $\Omega_\e$ the set $\{x\in\Omega:\dist(x,\partial\Omega)>\e\}$. 

\begin{theorem}\label{theorem:bdry-W2p}
Let $F \in C(\cS^n\times\R^n)$ be a functional satisfying \eqref{eq:F-ellip} -- \eqref{eq:Fb-W2VMO}, $\Omega\subset\R^n$ be a domain, $U\subset\R^n$ be an open neighborhood of a point of $\partial\Omega$ such that $\partial\Omega\cap U$ is a $C^1$-graph, $f\in L^p(\Omega\cap U)$ and $g\in W^{2,p}(\Omega\cap U)$, for some $p \in (p_0,\infty)$. Suppose that there is a diffeomorphism $\Phi\in C^1 (U;V)$ for which $\Phi(\Omega\cap U) = H(e_n)\cap V$ and $\Phi(\partial\Omega\cap U) = \partial H(e_n)\cap V$. Suppose that $\Phi\in W^{2,n}(U;V)$ if $p_0 < p < n$, $\Phi\in W^{2, n + \sigma}(U;V)$ for some $\sigma>0$ if $p = n$, and $\Phi \in W^{2,p}(U,V)$, if $p > n$. Let $u^\e\in C(\overline\Omega\cap U)$ be a viscosity solution to 
\begin{equation}\label{eq:main-bdry-W2p}
\begin{dcases}
F\left(D^2 u^\e ,\frac{\cdot}{\e}\right) = f &\text{in }\Omega\cap U, \\
u^\e = g &\text{on }\partial\Omega\cap U.
\end{dcases}
\end{equation}
Then $H_{\Omega\cap U}^\e (u^\e) \in L_{loc}^p(\overline\Omega_\e\cap U)$, provided that and for any $U'\Subset U$, 
$$
\begin{aligned}
\| H_{\Omega\cap U}^\e (u^\e)\|_{L^p(\Omega_\e\cap U')} \leq C \left( \frac{\| u^\e  \|_{L^\infty(\Omega\cap U)}}{\dist(U', \partial U)^{2-\frac{n}{p}}} + \| |f| + |D^2 g|\|_{L^p(\Omega\cap U)} \right),
\end{aligned}
$$
where $C>0$ depends at most on $n$, $\lambda$, $\Lambda$, $\psi$, $\kappa$, $p$, $\sigma$ and $\diam(U)$. Under the same additional assumption in Theorem \ref{theorem:int-W2p} with some $\bar p \in (p_0,\infty)$, the same assertion holds with $H_{\Omega\cap U} (u^\e)$, hence with $|D^2 u^\e|$, and $\Omega\cap U$ in place of $D_\e^2 u^\e$ and respectively $\Omega_\e \cap U$, in which case the constant $C$ may also depend on $\bar p$. 
\end{theorem}

Let us begin with a sub-optimal estimate, namely a uniform boundary $W^{2,q}$-estimates, with $q<p$, provided that $f\in L^p$ and $g\in C^{1,1-\frac{n}{p}}$. There are two reasons why we present this proposition separately. Firstly, $g\in C^{1,1-\frac{n}{p}}$ on the boundary layer is a much relaxed assumption than $g\in W^{2,p}$ in a domain containing the boundary layer, although the latter always implies the former. Hence, the proposition itself might be of independent interests. Secondly, our approximation lemma to appear below for the measure of the set of large ``Hessian" will be based on this proposition. It will be a key difference to the interior case (Lemma \ref{lemma:apprx-msr}), as in the case of boundary layer, we do not need to go down to the limit scale and carry over improved regularity of the effective profile. 

\begin{proposition}\label{proposition:bdry-W2p}
Let $F \in C(\cS^n\times\R^n)$ be a functional satisfying \eqref{eq:F-ellip} -- \eqref{eq:Fb-W2VMO}, $\Omega\subset\R^n$ be a domain, $U\subset\R^n$ be an open neighborhood of a point of $\partial\Omega$ such that $\partial\Omega\cap U$ is $(\delta,R)$-Reifenberg flat, $g\in C^{0,\alpha}(\partial\Omega\cap U)$ for some $\alpha\in(0,1]$, $f\in C(\Omega\cap U)\cap L^p(\Omega\cap U)$ for some $p > p_0$, and $u^\e\in C(\overline\Omega\cap U)$ be a viscosity solution to
\begin{equation}\label{eq:main-bdry-W2p}
\begin{dcases}
F\left(D^2 u^\e ,\frac{\cdot}{\e}\right) = f &\text{in }\Omega\cap U, \\
u^\e = g &\text{on }\partial\Omega\cap U.
\end{dcases}
\end{equation}
If $p_0 < p \leq n$, $H_{\Omega\cap U}^\e (u^\e) \in L_{loc}^q(\overline\Omega_\e\cap U)$ for all $q\in(p_0,\min\{p,\frac{n}{2-\alpha}\})$, and for any subdomain $U'\Subset U$, 
$$
\begin{aligned}
\| H_{\Omega\cap U}^\e (u^\e)\|_{L^q(\Omega_\e\cap U')}  \leq C \left( \frac{\| u^\e  \|_{L^\infty(\Omega\cap U)}}{\dist(U', \partial U)^{2-\frac{n}{q}}} + \| f \|_{L^p(\Omega\cap U)} + \| g \|_{C^{0,\alpha}(\partial\Omega\cap U)} \right),
\end{aligned}
$$
where $C>0$ depends only on $n$, $\lambda$, $\Lambda$, $\psi$, $\kappa$, $q$, $\alpha$ and $\diam(U)$. Moreover, if $p > n$, assume that $\partial\Omega\cap U$ is a $C^{1,\alpha}$-hypersurface, and $g\in C^{1,\alpha}(\partial\Omega\cap U)$. Then $H_{\Omega\cap U}^\e(u^\e) \in L_{loc}^q(\overline\Omega_\e\cap U)$ for all $q \in (p_0, \min\{p, \frac{n}{1-\alpha}\})$, and for any subdomain $U'\Subset U$, 
$$
\begin{aligned}
\| H_{\Omega\cap U}^\e (u^\e)\|_{L^q(\Omega_\e\cap U')}  \leq C \left( \frac{\| u^\e  \|_{L^\infty(\Omega\cap U)}}{\dist(U', \partial U)^{2-\frac{n}{q}}} + \| f \|_{L^p(\Omega\cap U)} + \| g \|_{C^{1,\alpha}(\partial\Omega\cap U)} \right).
\end{aligned}
$$
Under the additional assumption in Theorem \ref{theorem:int-W2p} with some $\bar p\in (p_0,\infty)$, the same assertions hold with $H_{\Omega\cap U} (u^\e)$, hence with $|D^2 u^\e|$, and $\Omega\cap U$ in place of $D_\e^2 u^\e$ and respectively $\Omega_\e \cap U$, in which case the constant $C$ may also depend on $\bar p$. 
\end{proposition}

\begin{proof}
The proof is essentially the same with that of Proposition \ref{proposition:bdry-W1p}. After a suitable rescaling argument, it may suffice to prove the case where $U = B_1$, $U' = B_{1/2}$, $\| u^\e\|_{L^\infty(\Omega\cap B_1)} \leq 1$, $\| f \|_{L^p(\Omega\cap B_1)} \leq 1$ and $\| g \|_{C^{0,\alpha}(\partial\Omega\cap B_1)} \leq 1$ if $p\leq n$ and $\| g\|_{C^{1,\alpha}(\partial\Omega\cap B_1)} \leq 1$ if $p > n$. Let $p_n = \min\{ p, \frac{n}{2-\alpha}\}$ if $p < n$, $p_n = \min\{p,\frac{n}{1-\alpha}\}$ if $p > n$, and $p_n = \gamma$ for some $\gamma \in (\frac{n}{2}, \min\{n,\frac{n}{2-\alpha}\})$. By \cite[Theorem 1]{LZ2} if $p\leq n$ or Lemma \ref{lemma:bdry-C1a} if $p > n$, one can find, for each $x_0\in\partial\Omega\cap B_{1/2}$, a linear polynomial $\ell_{x_0}$ (in case $p \leq n$ the linear polynomial $\ell_{x_0}$ is taken by the constant $u^\e(x_0)$) such that 
\begin{equation}\label{eq:he-C1a-bdry}
|(u^\e - \ell_{x_0})(x)| \leq c ( |x-x_0|^{2-\frac{n}{p_n}} + \e^{2 - \frac{n}{p_n}}),
\end{equation}
for all $x\in\overline\Omega\cap B_1$, where $c>0$ depends only on $n$, $\lambda$, $\Lambda$, $\kappa$ and $p_n$.

Now for each ball $B\subset\Omega_\e\cap B_{1/2}$ with $\partial (2B)\cap \partial\Omega\cap B_{1/2} \neq\emptyset$, we can make the following rescaling,
$$
u_B^{\e_B} (x) = \frac{ (u^\e - \ell_{x_{B,0}})(x_B + \rho_B x)}{\rho_B^{2-\frac{n}{p_n}}},
$$ 
of $u^\e$, where $x_B$ is the center of $B$, $\rho_B$ its radius and $x_{B,0}$ the point of intersection between $\partial(2B)$ and $\partial\Omega\cap B_{1/2}$. Then we may repeat the proof of Proposition \ref{proposition:bdry-W1p}, utilizing Theorem \ref{theorem:int-W2p} in place of Theorem \ref{theorem:int-W1p}, to deduce that $\| H_{\Omega_B}^{\e_B} (u_B^{\e_B}) \|_{L^q(B_1)} \leq c$ for any $q\leq p$, with $\Omega_B = \rho_B^{-1} (-x_B + \Omega\cap B_1)$, whence 
$$
\| H_{\Omega\cap B_1}^\e (u^\e) \|_{L^q(B)} \leq c\rho_B^{\frac{n}{q} - \frac{n}{p_n}}. 
$$
Fix any $q < p_n$. Then we can consider the same Besicovitch covering $\cG$, as in the proof of Proposition \ref{proposition:bdry-W1p}, of $\Omega_\e\cap B_{1/2}$ by balls $B$, such that the summation of the right hand side of over all $B\in\cG$ is bounded by a constant $c$. This finishes the proof. 
\end{proof} 

As for the proof of Theorem \ref{theorem:bdry-W2p}, it suffices to consider the case where $f \in L^p$ and $g = 0$, since one may always substitute $u^\e$ with $u^\e - g$ and $f$ with $f + c|D^2 g|$. 

Since our analysis will be of local nature around a boundary point, and will be invariant under translation, we shall work from now on with domains $\Omega$ with $0\in\partial\Omega$, $U = B_2$ and $U' = B_1$. Unless specified otherwise, from now on, $F\in C(\cS^n\times\R^n)$ satisfies \eqref{eq:F-ellip} -- \eqref{eq:Fb-W2VMO}, $\partial\Omega\cap B_2$ is a $C^1$-hypersurface containing the origin, and that there is a diffeomorphism $\Phi \in C^1(B_2;V)$, with $V\subset\R^n$ a neighborhood of the origin, such that $\Phi(0) = 0$, $\Phi(\Omega\cap B_2) = H_0(e_n)\cap V$ and $\Phi (\partial\Omega\cap B_2) = \partial H_0(e_n)\cap V$. We shall call $\Phi$ boundary flattening map (of $\partial\Omega$) around the origin. Moreover, $f\in C(\Omega\cap B_2)\cap L^p(\Omega\cap B_2)$ for some $p > p_0$, and $u^\e \in C(\overline\Omega\cap B_2)$ is a viscosity solution to \eqref{eq:main-bdry-W2p} with $U = B_2$ and $g = -\ell$, with $\ell$ a linear polynomial; we shall discuss later in detail the reason for the involvement of a linear polynomial in the boundary condition.



The difficulty of our analysis arises from the fact that homogenization problems are unfavorable towards boundary flattening argument, as one loses the oscillating pattern by the transformation. In one way or another, one will resort to the fact that the original problem in small scales is homogenized to a ``nice'' effective problem to improve the regularity, whence the level of difficulty remains the same. 

For this reason, we shall study our problem \eqref{eq:main-bdry-W2p} without flattening the boundary. This readily implies some notable changes in the approximation lemma below for the measure of the set of large ``Hessian", compared to the interior case (Lemma \ref{lemma:apprx-msr}) as well as those for standard problems in the setting of flat boundaries (e.g., \cite[Lemma 2.14]{W}). 




\begin{lemma}\label{lemma:apprx-msr-bdry}
Let $\e$, $\delta$, $\alpha$, $\rho$, $p$ and $q$ be constants with $0<\e<1$, $0<\delta \leq \delta_0$, $0< \alpha < 1$, $\rho > 0$ and $p_0 < p  < q < \infty$ be given. Suppose that $| \xi^t D\Phi (0) \xi- 1| \leq \delta$ for any $\xi\in\partial B_1$, $\osc_{B_2} D\Phi \leq \delta$, $\|f \|_{L^p(\Omega\cap B_2)} \leq \delta$, $\| u^\e  \|_{L^\infty(\Omega\cap B_2)} \leq 1$, $|u^\e(x)|\leq |x|^2$ for all $x\in \Omega\setminus B_1$, and $u^\e = -\ell$ on $\partial\Omega\cap B_2$, for some linear polynomial $\ell$. Assume either of the following:
\begin{enumerate}[(i)]
\item $\displaystyle |D\ell| \leq \frac{1}{\rho}$, $\| D^2 \Phi\|_{L^q(B_2)} \leq \delta \rho$ and $q < n$;
\item $\displaystyle |D\ell| \leq \frac{1}{\rho}$, $\| D^2 \Phi \|_{L^q(B_2)} \leq \delta\rho $, and $n< q < \frac{n}{1-\alpha}$;
\item $\displaystyle |D\ell| \leq \frac{1}{\rho^{1-\alpha}}$, $\| D^2 \Phi\|_{L^n(B_2)} \leq \delta\rho$ and $ q> n$. 
\end{enumerate}
Then for any $s>0$, 
$$
|A_s^\e (u^\e, \Omega)\cap \Omega_\e \cap B_1| \leq c( \delta^{\gamma\mu} s^{-\mu} + s^{-q}),
$$
where $\mu > 0$ depends only on $n$, $\lambda$ and $\Lambda$, $\delta_0 > 0$ and $0 < \gamma \leq 1$ depend in addition to $q$ and respectively $\alpha$, and $c>1$ may depend further on $\kappa$ and $\psi$. Nevertheless, none of $\mu$, $\delta_0$, $\gamma$ and $c$ depends on $\rho$ or $\e$. 
\end{lemma}

\begin{proof}
Set $T = \frac{1}{|B_2|}\int_{B_2} D\Phi\,dx$, and let $L_T : \R^n\to\R^n$ be the linear transformation induced by $T$; i.e., $L_T(0)=0$ and $DL_T = T$. In what follows, we shall denote by $c_0$ a constant depending at most on $n$ and $q$, and we shall allow it to vary from one line to another. 

\begin{case}
$|D\ell| \leq \rho^{-1}$, $\|D^2 \Phi\|_{L^q(B_2)} \leq \delta\rho$ and $q<n$. 
\end{case}

By the Poincar\'e inequality, together with $\|D^2 \Phi\|_{L^q(B_2)}\leq \delta \rho$, 
\begin{equation}\label{eq:DPhi-T-Ln}
\int_{B_2} |D\Phi - T|^q\,dx \leq c_0\delta^n \rho^n. 
\end{equation}
Let $L_T:\R^n\to\R^n$ be the linear transformation such that $DL_T = T$. Then by the Sobolev embedding theorem, one can infer that 
\begin{equation}\label{eq:Phi-L-Ca}
[\Phi - L_T ]_{C^{0,2-\frac{n}{q}}(B_2)} \leq c_0\delta \rho.
\end{equation}

Now by the assumption that $\Phi(\partial\Omega\cap B_2) \subset \partial H_0(e_n)$, i.e., $e_n\cdot \Phi = 0$ on $\partial\Omega\cap B_2$, \eqref{eq:Phi-L-Ca} yields that $|e_n\cdot T (x-y)| \leq c_0\delta\rho|x-y|^{2-\frac{n}{q}}$ for any $x,y\in\partial\Omega\cap B_2$. Recall that $T = \frac{1}{|B_2|}\int_{B_2} D\Phi\,dx$. By the assumption on $\Phi$, $|\xi^t T\xi - 1| \leq 2\delta$ for all $\xi\in\partial B_1$, and hence $|Te_n| \geq e_n^t Te_n \geq 1-2\delta > \frac{1}{2}$, provided that $\delta\leq\frac{1}{4}$. Set $\nu = \frac{T e_n}{|T e_n|}$. Then the latter observation implies that 
\begin{equation}\label{eq:bdry-S}
|\nu\cdot (x-y)| \leq 2c_0\delta \rho |x-y|^{2-\frac{n}{q}},
\end{equation}
for any $x,y\in\partial\Omega\cap B_2$.

Let us now turn to the regularity of the linear polynomial $\ell$, for which $u^\e = - \ell$ on $\partial\Omega\cap B_2$. Since $|u^\e |\leq 1$ in $\Omega\cap B_2$, $|\ell| \leq 1$ on $\partial\Omega\cap B_2$. Now we may deduce from $|D\ell|\leq \rho^{-1}$, \eqref{eq:bdry-S} and $0\in \partial\Omega$ that with $S = S_{2c_0\delta \rho}(\nu)$, 
\begin{equation}\label{eq:l-Linf}
\sup_{S\cap B_2} |\ell|  < 1 + \frac{4\delta  \rho}{\rho} \leq 2.
\end{equation}
Therefore, denoting by $P_\nu(e)$ the orthogonal projection of a vector $e$ in direction $\nu$, one can compute that 
\begin{equation}\label{eq:l-C01}
|P_{\nu}(D\ell)| \leq \osc_{S\cap B_2} \ell \leq 4.
\end{equation}
Putting \eqref{eq:l-C01} together with \eqref{eq:bdry-S}, we also derive that for each $\alpha\in(0,1)$, we can take $\delta$ small, depending at most on $n$ and $\beta$, such that 
\begin{equation}\label{eq:l-Ca}
\begin{aligned}
[ \ell ]_{C^{0,2-\frac{n}{q}}(\partial\Omega\cap B_2)} & \leq 2^{\frac{n}{q}-1} |P_\nu(D\ell)| + c_0\delta \rho |D_\nu \ell|   \leq 9,
\end{aligned}
\end{equation}
where the last inequality holds for any small $\delta$, whose smallness condition depends only on $n$ and $q$. 

Let $h^\e \in C(\overline{\Omega\cap B_2})$ be a viscosity solution to
\begin{equation}\label{eq:he-pde}
\begin{dcases}
F\left( D^2 h^\e ,\frac{\cdot}{\e}\right) = 0 &\text{in } \Omega\cap B_2,\\
h^\e = u^\e & \text{on }\partial(\Omega\cap B_2).
\end{dcases}
\end{equation}
The existence of such a viscosity solution is ensured by \cite[Theorem 1]{Sir}, since $\partial\Omega\cap B_2$ is a $C^1$-hypersurface, and that $F\in C(\cS^n\times\R^n)$ and $u^\e \in C(\partial (\Omega\cap B_2))$. By the maximum principle, we have 
\begin{equation}\label{eq:he-Linf}
\| h^\e \|_{L^\infty(\Omega\cap B_2)}  \leq 1. 
\end{equation}
By the assumption on $\Phi$, $\partial\Omega\cap B_2$ is a $C^1$-hypersurface whose Lipschitz norm is less than $c_0\delta$. Thus, by taking $\delta$ smaller if necessary, depending now on $n$, $\lambda$, $\Lambda$ and $q$, we can deduce from Proposition \ref{proposition:bdry-W2p} (with $\alpha = 2-\frac{n}{q}$), \eqref{eq:l-Ca} and \eqref{eq:he-Linf} that $H_{\Omega \cap B_2}^\e( h^\e) \in L^q(\Omega_\e \cap B_1 )$ and
\begin{equation}\label{eq:Tehe-Lq}
\| H_\Omega^\e( h^\e, \Omega\cap B_2) \|_{L^q(\Omega_\e\cap B_1)} \leq c,
\end{equation} 
for some constant $c>0$ depending only on $n$, $\lambda$, $\Lambda$, $\kappa$, $\psi$ and $q$. In particular, by the definition $A_s^\e (h^\e,E) = \{H_E^\e (h^\e) > s\}$, we obtain
\begin{equation}\label{eq:Aehe-Lq}
| A_s^\e ( h^\e, \Omega\cap B_2) \cap \Omega_\e \cap B_1 | \leq cs^{-q},
\end{equation} 
for any $s>0$. 

In comparison of \eqref{eq:he-pde} with \eqref{eq:main-bdry-W2p} (with $g = -\ell$), Lemma \ref{lemma:visc} ensures that $w^\e =  \delta^{-1} (u^\e -  h^\e)$ is a viscosity solution to
\begin{equation}\label{eq:ue-he-pde}
\begin{dcases}
\cP_- (D^2 w^\e) \leq \frac{f}{\delta} \leq \cP_+ (D^2 w^\e) & \text{in }\Omega\cap B_2,\\
w^\e = 0 &\text{on }\partial(\Omega\cap  B_2), 
\end{dcases}
\end{equation} 
Owing to $\| f\|_{L^{p}(\Omega\cap B_2)}\leq \delta$, the generalized maximum principle ensures that $\| w^\e \|_{L^\infty(\Omega\cap B_{3/2})} \leq c$. Then by Proposition \ref{proposition:W2d}, we obtain 
\begin{equation}\label{eq:we-As}
|A_t (w^\e, \Omega\cap B_{3/2})| \leq c_0 t^{-\mu},
\end{equation} 
for any $t>1$, where $\mu \in (0,1)$ depends only on $n$, $\lambda$ and $\Lambda$. Combining \eqref{eq:we-As} with \eqref{eq:Aehe-Lq}, we arrive at 
$$
|A_s (u^\e, \Omega\cap B_2)\cap \Omega_\e\cap B_1| \leq c(\delta^\mu s^{-\mu} + s^{-q}). 
$$
We can then replace $\Omega\cap B_2$ above with $\Omega$ by invoking the inequality $|u^\e(x)| \leq |x|^2$ for all $x\in\Omega\setminus B_1$. 

\begin{case}
$|D\ell| \leq \rho^{-1}$, $\|D^2\Phi\|_{L^q(B_2)} \leq \delta\rho$ and $n \leq q < \frac{n}{1-\alpha}$.
\end{case}

As for this case, by the Sobolev embedding theorem, 
\begin{equation}\label{eq:DPhi-Ca}
[D\Phi]_{C^{0,1-\frac{n}{q}}(B_2)} \leq c_0 \delta\rho.
\end{equation}
Especially, since $|e_n\cdot D\Phi| \geq e_n^t D\Phi e_n \geq 1-2\delta >\frac{1}{2}$ on $B_2$, denoting by $\nu_x$ the unit vector $\frac{e_n\cdot \Phi(x)}{|e_n\cdot D\Phi(x)|}$, we obtain that 
\begin{equation}\label{eq:nu-C1a}
|\nu_x - \nu_y | \leq 2c_0 \delta\rho |x-y|^{1-\frac{n}{q}},
\end{equation}
for any $x,y\in\partial\Omega\cap B_2$. Moreover, one can also deduce from \eqref{eq:DPhi-Ca} that $|\Phi(y) - \Phi(x) - D\Phi(x)\cdot (y-x)| \leq c_0 \delta \rho |y-x|^{2-\frac{n}{q}}$, for any $x,y\in B_2$. Utilizing $e_n\cdot \Phi = 0$ on $\partial\Omega\cap B_2$, we also obtain that 
\begin{equation}\label{eq:bdry-C1a}
|\nu_x\cdot (y-x)| \leq 2c_0 \delta\rho |y-x|^{2 - \frac{n}{q}},
\end{equation}
and In other words, $\partial\Omega\cap B_2$ is a $C^{1,1-\frac{n}{q}}$-hypersurface whose $C^{1,\alpha}$-norm is bounded by $2c_0 \delta\rho$. 

For the rest of the proof, we shall denote by $c_\alpha$ a positive constant depending at most on $n$ and $\alpha$, and it may vary at each occurrence. With \eqref{eq:bdry-C1a} at hand, we claim that $\ell \in C^{1,1-\frac{n}{q}}(\partial\Omega\cap B_2)$ and 
\begin{equation}\label{eq:l-C1a}
\| \ell \|_{C^{1,1-\frac{n}{q}}(\partial\Omega\cap B_2)} \leq c_0.
\end{equation}
Note that the hypothesis of Case 2 is stronger than that of Case 1, whence \eqref{eq:l-Linf} and \eqref{eq:l-C01} continue to hold, with the bounds possibly replaced by $c_0$. Thus, \eqref{eq:l-C01} together with \eqref{eq:bdry-C1a} (with $x_0 = 0$) and $|D\ell| \leq \rho^{-1}$ implies that 
$$
[\ell]_{C^{0,1}(\partial\Omega\cap B_2)} \leq |P_\nu(D\ell)| + c_0 \delta \rho |D_\nu \ell| \leq c_0.
$$
Moreover, we may also compute, via \eqref{eq:nu-C1a} and $|D\ell|\leq \rho^{-1}$, that 
$$
[D_\tau \ell]_{C^{0,1-\frac{n}{q}}(\partial\Omega\cap B_2)} \leq |D\ell| \sup_{x,y\in\partial\Omega\cap B_2}\frac{|\nu_x - \nu_y |}{|x-y|^{1-\frac{n}{q}}} \leq c_0,
$$
where $D_\tau$ is the tangential gradient to $\partial\Omega\cap B_2$. Combining the last two displays altogether with \eqref{eq:l-Linf}, we verify the claim \eqref{eq:l-C1a}. 

Now let $h^\e$ be the viscosity solution to \eqref{eq:he-pde}, as in Case 1. The inequality in \eqref{eq:he-Linf} continues to hold here. However, now with \eqref{eq:bdry-C1a} and \eqref{eq:l-C1a} at hand, Proposition \ref{proposition:bdry-W2p} ensures that \eqref{eq:Tehe-Lq}, hence \eqref{eq:Aehe-Lq} as well, holds for $q \leq \frac{n}{1-\alpha}$. The rest of the proof repeats that of Case 1 verbatim, so it is omitted. 

\begin{case}
$|D\ell|\leq \rho^{\alpha-1}$, $\|D^2 \Phi\|_{L^n(B_2)} \leq \delta\rho$ and $q > n$. 
\end{case}

In the above cases, the auxiliary Dirichlet problem for the approximating function $h^\e$ was imposed on the same domain $\Omega\cap B_2$. Thus, the integrability of $H_{\Omega\cap B_2}^\e(h^\e)$ cannot exceed the exponent determined by the regularity of the boundary layer. On contrast, this last case asks for $q$ to go over the threshold. To achieve this goal, we shall consider another Dirichlet problem, whose boundary layer is much smoother (in fact, a hyperplane). However, the choice of such an auxiliary problem cannot be arbitrary, since the newly obtained approximating function should also be sufficiently close to the original solution on the boundary layer, so that their difference satisfies \eqref{eq:we-As}. To meet the latter requirement, we need $\ell$ to be $C^{0,\alpha}$-regular for some $\alpha>0$, not only on $\partial\Omega\cap B_2$ (as in \eqref{eq:l-Ca}), but also over a slab $S\cap B_2$, with $S = S_{c_0\delta\rho}(\nu)$, that contains $\partial\Omega\cap B_2$. This is where a stronger assumption, $|D\ell| \leq \rho^{\alpha -1}$, is used. 

Let us explain this in more detail. In what follows, let $E$ denote the half-space $H_{-c_0\delta\rho}(\nu)$; note that $\Omega\cap B_2\subset E$. We shall keep writing by $S$ the slab $S_{c_0\delta\rho}(\nu)$; recall from \eqref{eq:bdry-S} that $\partial\Omega\cap B_2\subset S$. 

With $\alpha$ being the exponent for which $|D\ell|\leq \rho^{\alpha-1}$, we may compute, by using \eqref{eq:l-C01}, that
\begin{equation}\label{eq:l-Ca-re}
\begin{aligned}
[\ell]_{C^{0,\alpha}(S\cap B_2)} &\leq 2^{1-\alpha} |P_\nu(D\ell)| + |D_\nu\ell| \sup_{x,y\in S\cap B_2} \frac{|\nu\cdot (x-y)|}{|x-y|^\alpha} \\
&\leq 8 + \rho^{\alpha-1} \sup_{x,y\in S\cap B_2} |\nu\cdot (x - y)|^{1-\alpha} \\
&\leq 9,
\end{aligned}
\end{equation}
where the last inequality again follows by choosing $\delta$ small, depending only on $n$ and $\alpha$. 

Let us remark that $u^\e \in C(\overline\Omega\cap B_2)$, $\cP_- (D^2 u^\e ) \leq f\leq \cP_+(D^2 u^\e)$ in $\Omega\cap B_2$ in the viscosity sense, $\|u^\e\|_{L^\infty(\Omega\cap B_2)} \leq 1$, $\|f\|_{L^p(\Omega\cap B_2)} \leq \delta$ and $u^\e = -\ell$ on $\partial\Omega\cap B_2$. Moreover, $\Omega\cap B_2$ satisfies a uniform exterior cone condition, where the size of the cone is bounded by an absolute constant, because of the assumption on the boundary flattening map $\Phi$. Therefore, one may employ \cite[Theorem 2]{Sir}, along with \eqref{eq:l-Ca-re} (in fact,  \eqref{eq:l-Ca} works equally well here) to deduce that $u^\e \in C^{0,2\gamma}(\overline{\Omega\cap B_{3/2}})$, and 
\begin{equation}\label{eq:ue-Ca}
[ u^\e ]_{C^{0,2\gamma}(\overline{\Omega\cap B_{3/2}})} \leq c,
\end{equation} 
for some $\gamma \in (0,\frac{\alpha}{4}]$ and $c>0$ depending at most on $n$, $\lambda$, $\Lambda$, $p$ and $\alpha$. 

Now let $\phi^\e \in C(\overline{E\cap B_{3/2}})$ be a viscosity solution to 
\begin{equation}\label{eq:phie-pde}
\begin{dcases}
F\left(D^2 \phi^\e ,\frac{\cdot}{\e}\right) = 0 &\text{in }E\cap B_{3/2},\\
\phi^\e = u^\e & \text{on }\Omega\cap\partial B_{3/2},\\
\phi^\e = -\ell & \text{on }\partial (E\cap B_{3/2})\setminus \Omega. 
\end{dcases}
\end{equation} 
Note that  $\phi^\e\in C(\partial(E\cap B_2))$, since $u^\e = -\ell$ on $\partial\Omega\cap B_2$, and that $E\cap B_2$ satisfies a uniform exterior sphere condition with radius $1$. Hence, the existence of a viscosity solution to the above problem is obvious. Also by \eqref{eq:l-Ca-re} and \eqref{eq:ue-Ca}, $\phi^\e \in C^{0,2\gamma}(\partial (E\cap B_{3/2}))$, and thus, it follows from e.g., \cite[Proposition 4.13, 4.14]{CC}, that 
\begin{equation}\label{eq:phie-Ca}
[\phi^\e ]_{C^{0,\gamma}(\overline{E\cap B_{3/2}})} \leq c. 
\end{equation} 

On the other hand, by the maximum principle, $\|u^\e \|_{L^\infty(\Omega\cap B_2)}\leq 1$ and \eqref{eq:l-Linf}, we have $\| \phi^\e \|_{L^\infty(E\cap B_{3/2})} \leq 2$. Moreover, as $\partial E$ being an hyperplane orthogonal to $\nu$,  $|D_\tau \ell| = |P_\nu (D\ell)|$ and $|D_\tau^2 \ell| = 0$ on $\partial E$, with $D_\tau$ being the tangential gradient on $\partial E$. Therefore, it follows from \eqref{eq:l-Linf} and \eqref{eq:l-C01} that 
\begin{equation}\label{eq:l-C11}
\| \ell \|_{C^{1,1}(\partial E\cap B_{3/2})} \leq 6.
\end{equation} 
By \eqref{eq:l-C11} (as well as $\|\phi^\e \|_{L^\infty(E\cap B_2)} \leq 2$, Proposition \ref{proposition:bdry-W2p} (now with $f= 0$, $g = -\ell$, $p = q$ and $\alpha = 1$ there) yields that $H_{E\cap B_{3/2}}^\e (h^\e ) \in L^q(E_\e \cap B_1)$, with the chosen $q \in (n,\infty)$ and $E_\e = \{ x\in E: \dist(x,\partial E)> \e\}$, and hence arguing as in the derivation of \eqref{eq:Tehe-Lq}, we similarly obtain that 
\begin{equation}\label{eq:Hephie-Lq}
|A_s^\e (\phi^\e, E\cap B_{3/2})\cap E_\e\cap B_1| \leq cs^{-q},
\end{equation} 
where $c$ may now depend at most on $n$, $\lambda$, $\Lambda$, $\psi$, $\kappa$ and $q$. 

Set $v^\e = \frac{u^\e - \phi^\e}{c\delta^\gamma}$, with a possibly different $c>1$ to that in the last display, yet depending on the same quantities. Then again by Lemma \ref{lemma:visc}, we may compare \eqref{eq:main-bdry-W2p} (with $g=-\ell$) with \eqref{eq:phie-pde} and compute that
\begin{equation}\label{eq:ve-pde}
\begin{dcases}
\cP_- (D^2 v^\e ) \leq \frac{f}{c\delta^\gamma} \leq \cP_+(D^2 v^\e) & \text{in }\Omega\cap B_{3/2},\\
v^\e = 0 & \text{on }\Omega\cap \partial B_{3/2}, \\
v^\e = \frac{-\ell - \phi^\e }{c\delta^\gamma} &\text{on }\partial\Omega\cap B_{3/2}.
\end{dcases} 
\end{equation}
As $c>1$ and $0< \delta, \gamma < 1$, it follows from the assumption that $\|f \|_{L^p(\Omega\cap B_2)} \leq c\delta^\gamma$. This implies, by the general maximum principle, that
$$
\| v^\e \|_{L^\infty(\Omega\cap B_{3/2})} \leq c_0 + \sup_{\partial\Omega\cap B_{3/2}} \frac{|\phi^\e + \ell|}{c\delta^\gamma}. 
$$
Therefore, once the rightmost term is proved to be bounded by an absolute constant (by choosing $c>1$ large), we may repeat the final part of the proof for Case 1, and derive the desired decay estimate. Since the latter implication is already shown above, let us finish the proof by justifying that $| \phi^\e + \ell | \leq c\delta^\gamma$ on $\partial\Omega\cap B_{3/2}$. 

To this end, let  $x\in \partial\Omega\cap B_{3/2}$ be any, and find $x_0 \in \partial E\cap B_{3/2}$ such that $|x-x_0| \leq c_0\delta \rho$. Such a point $x_0$ always exists because $\partial\Omega\cap B_2\subset S$ with $S$ being a slab with width $c_0\delta\rho$. Then by \eqref{eq:l-Ca-re} and \eqref{eq:phie-Ca}, 
$$
\begin{aligned}
|\phi^\e (x) + \ell(x)| & \leq |\phi^\e (x) + \ell(x_0)| + |\ell(x) - \ell(x_0)| \\
&\leq c( (\delta \rho)^\gamma + (\delta \rho)^\alpha) \\
&\leq c \delta^\gamma,
\end{aligned}
$$
where the last inequality is ensured by $\gamma \leq\frac{\alpha}{4}$, $\delta<1$ and $\rho\leq 1$. This completes the proof.  
\end{proof}




Our next step is to design a suitable iteration argument for the boundary estimate. 

\begin{lemma}\label{lemma:iter-msr-bdry}
Let $\e$, $\delta$, $\alpha$, $\rho$, $p$ and $q$ be constants with $0<\e<1$, $0<\delta \leq \delta_0$, $0< \alpha < 1$, $\rho > 0$ and $p_0 < p  < q < \infty$ be given. Suppose that $| \xi^t D\Phi (0) \xi- 1| \leq \delta$ for any $\xi\in\partial B_1$, $\osc_{B_2} D\Phi \leq \delta$, $\|D^2 \Phi\|_{L^n(B_2)} \leq \delta$, $\|f \|_{L^p(\Omega\cap B_2)} \leq \delta$, $\| u^\e  \|_{L^\infty(\Omega\cap B_2)} \leq 1$, and $u^\e = 0$ on $\partial\Omega\cap B_2$. Set $A_k = A_{m^k}^\e (u^\e,\Omega\cap B_2) \cap \Omega_\e\cap B_1$ and $B_k = \{M(|f|^p \chi_{\Omega\cap B_2}) ) > \delta^p m^{k p}\}$. Assume either of the following: 
\begin{enumerate}[(i)]
\item $C_k = L_{m^{k(1-p/n)}}^\e (u^\e, \Omega\cap B_2)\cap \Omega_\e\cap B_{3/2}$, $D_k = \{ M(|D^2 \Phi|^q\chi_{B_2} )> \delta^q m^{\frac{kpq}{n}}\}$, and $p < q < n$; 
\item $C_k =L_{m^{k\alpha}}^\e(u^\e, \Omega\cap B_2)\cap \Omega_\e\cap B_{3/2}$, $D_k = \{ M(|D^2 \Phi|^q\chi_{B_2}) > \delta^q m^{kn(1-\alpha)q}\}$ and $(1-\alpha)n = p < n <  q < \frac{n}{1-\alpha}$;
\item $\| D^2 \Phi \|_{L^p(B_2)} \leq \delta$, $C_k = \emptyset$, $D_k = \{M(|D^2 \Phi|^p\chi_{B_2} )> \delta m^{kp}\}$ and $n < p< q$.
\end{enumerate}
Then for each integer $1\leq k\leq - c(n,p)\frac{\log\e}{\log m}$, with $c(n,p) = \frac{n}{p}$ for case (i), $c(n,p) = \frac{1}{1-\alpha}$ for case (ii) and $c(n,p) = \frac{p}{n}$ for case (iii),  
$$
|A_{k+1}| \leq c(\delta^{\gamma\mu} m^{-\mu} + m^{-q}) | A_k \cup B_k \cup C_k\cup D_k|,
$$
where $\mu> 0$ depends only on $n$, $\lambda$ and $\Lambda$, $\delta_0\in(0,\frac{1}{4})$ and $\gamma\in(0,1]$ depends at most on $q$ and $\alpha$, while $c> 0$ and $m>1$ may depend further on $\psi$ and $\kappa$. 
\end{lemma}

\begin{proof}
Since $u^\e = 0$ on $\partial\Omega\cap B_2$, we may apply Lemma \ref{lemma:apprx-msr-bdry} (with $\ell = 0$ and $\rho = 1$) to deduce that $ |A_1| \leq c(\delta^{\gamma\mu} m^{-\mu} + m^{-q}) \leq \eta_0 |B_1|$. Let us remark that this initial step applies to all three cases considered in the statement. Fix an integer $k\geq 1$. Since $A_{k+1} \subset A_k \subset \cdots \subset A_1$, we readily obtain $|A_{k+1}| \leq \eta_0 |B_1|$. 

Next, let $B\subset B_1$ be a ball whose center lies in $\Omega_\e\cap B_1$ and $\rad (B) \leq 1$. Assume that $|A_{k+1}\cap B| > \eta_0 |B|$. Our goal is to show that $\Omega_\e\cap B\subset (A_k\cup B_k\cup C_k\cup D_k)$. Thus, by Lemma \ref{lemma:vitali} (which applies to $\Omega_\e\cap B_1$, since $\partial\Omega_\e\cap B_1$ is $(c_0\delta,2)$-Reifenberg flat, which can be easily inferred from the $C^1$-character of $\partial\Omega\cap B_1$), $|A_{k+1}| \leq c\eta_0 |A_k\cup B_k\cup C_k\cup D_k|$. 

To the rest of the proof, we shall assume $\Omega_\e\cap B\setminus (A_k \cup B_k \cup C_k\cup D_k)\neq\emptyset$, and attempt to derive a contradiction against $|A_{k+1} \cap B| > \eta_0 |B|$. Arguing as in Lemma \ref{lemma:iter-W1p}, we observe that 
\begin{equation}\label{eq:e-r}
\e < 2r_B,
\end{equation}
by selecting $m>1$ large, yet depending only on $n$.

Choose any $\tilde x_B\in \Omega_\e\cap B\setminus (A_k\cup B_k \cup C_k\cup D_k)$. We shall only consider the case $B_{2r_B}(\tilde x_B)\setminus \Omega\neq\emptyset$, as the other case can be treated as the interior case. Under this setting, we can find $x_B\in \partial\Omega\cap B_1$ such that $\e < |x_B - \tilde x_B| = \dist(\tilde x_B,\partial\Omega) < 2r_B$.  

For the rest of the proof, we shall denote by $c>1$ a constant depending at most on $n$, $\lambda$, $\Lambda$, $\kappa$, $\psi$, $p$, $q$ and $\sigma$, and we allow it to vary at each occurrence. 

\vspace{0.2cm}
\noindent
{\it Case }1. $C_k = L_{m^{k(1-\frac{p}{n})}}^\e (u^\e, \Omega\cap B_2)\cap \Omega_\e\cap B_{3/2}$, $D_k = \{ M(|D^2 \Phi|^n\chi_{B_2}) > \delta m^{kp}\}$, and $p < q < n$.
\vspace{0.2cm}

By $\tilde x_B\in B\setminus D_k$ and $B_{2r_B}(x_B) \subset B_{4r_B}(\tilde x_B)$, 
\begin{equation}\label{eq:D2Phi-Ln}
\int_{B_{2r_B}(x_B)} |D^2 \Phi|^q\,dx \leq (4r_B)^n (\delta m^{\frac{kp}{n}})^q . 
\end{equation}
Moreover,  thanks to $\tilde x_B\not\in A_k$, there exists a linear polynomial $\ell$ such that for all $x\in\Omega\cap B_2$, 
\begin{equation}\label{eq:ue-l}
|(u^\e - \ell)(x)| \leq \frac{m^k}{2} (|x-\tilde x_B|^2 + \e^2),
\end{equation}
for all $x\in\Omega\cap B_2$. Also observe from $\tilde x_B\not\in C_k$ that $|u^\e (x) - a| \leq m^{k(1-\frac{p}{n})} (|x-\tilde x_B| + \e)$ for all $x\in\Omega\cap B_2$ for some constant $a\in\R$. Putting this observation together with \eqref{eq:ue-l}, $m^{\frac{kp}{n}} \e \leq 1$ and $B_\e(\tilde x_B)\subset \Omega$, we obtain $|\ell(x) - a |\leq m^k\e^2 + 2m^{k(1-\frac{p}{n})}\e \leq 3m^{k(1-\frac{p}{n})}\e$ for all $x\in B_\e(\tilde x_B)$. In particular, we arrive at 
\begin{equation}\label{eq:Dl}
|D\ell| \leq 6m^{k(1-\frac{p}{n})}. 
\end{equation} 
In addition, from $\tilde x_B \not\in B_k$ and $B_{2r_B}(x_B) \subset B_{4r_B}(\tilde x_B)$, we also have 
\begin{equation}\label{eq:f-Lp}
\int_{\Omega\cap B_{2r_B}(x_B)} |f|^p\,dx \leq (4r_B)^n (\delta m^k)^p. 
\end{equation}

Denote by $\Omega_B$ the rescaled domain $\frac{1}{2r_B} (-x_B + \Omega\cap B_2)$, and define $\Phi_B$ by 
\begin{equation}\label{eq:PhiB}
\Phi_B (x) = \frac{1}{2r_B} ( \Phi ( x_B + 2r_B x) - \Phi( x_B) ). 
\end{equation}
Since $\osc_{B_r(x)} D\Phi \leq \delta r$ for any $x\in\partial\Omega\cap B_2$ and any $r\in(0,1)$, we have $\osc_{B_2} D\Phi_B \leq \delta$. Moreover, since $|\xi^t D\Phi(0) \xi - 1 | \leq \delta$ for all $\xi\in\partial B_1$ and $\osc_{B_1} D\Phi \leq \delta$, it follows that $|\xi^t D\Phi_B(0) \xi -1| \leq 2\delta$ for all $\xi\in\partial B_1$. Furthermore, according to \eqref{eq:D2Phi-Ln}, $\| D^2 \Phi_B\|_{L^q(B_2)} \leq 4\delta m^{\frac{kp}{n}}r_B$. 

Consider the following rescaled versions of $u^\e$, $\ell$ and $f$, 
\begin{equation}\label{eq:uBeB-eB-lB-fB}
\begin{aligned}
& u_B^{\e_B} (x) = \frac{ (u^\e - \ell) (x_B + 2r_Bx)}{cm^kr_B^2 },\quad \e_B = \frac{\e}{2r_B} \\
& \ell_B (x) = \frac{\ell (x_B + 2r_Bx)}{cm^kr_B^2},\quad f_B(x) = \frac{f(x_B + 2r_B x)}{cm^k},
\end{aligned}
\end{equation}
with $c>1$ to be determined later. By \eqref{eq:ue-l}, we have $|u_B^{\e_B} (x)| \leq 1$ for all $x\in \Omega_B\cap B_2$ and $|u_B^{\e_B}(x)| \leq |x|^2$ for all $x\in\Omega\setminus B_1$, and by \eqref{eq:f-Lp}, $\|f_B\|_{L^p(\Omega_B\cap B_2)} \leq \delta$, while \eqref{eq:Dl} ensures $|D\ell_B| \leq (m^{\frac{pk}{n}} r_B)^{-1}$, provided that we take $c>1$ larger if necessary. In view of \eqref{eq:main-bdry-W2p} and $u^\e = 0$ on $\partial\Omega\cap B_2$, one may also compute that 
\begin{equation}\label{eq:uBeB-pde}
\begin{dcases}
F_B\left( D^2 u_B^{\e_B}, \frac{\cdot}{\e_B} \right) = f_B &\text{in }\Omega_B\cap B_2,\\
u_B^{\e_B} = -\ell_B &\text{on }\partial\Omega_B\cap B_2,
\end{dcases}
\end{equation} 
in the viscosity sense, where $F_B(P,y) = \frac{1}{cm^k} F( cm^k P, y + \frac{x_B}{\e})$. Obviously, $F_B \in C(\cS^n\times\R^n)$ and it verifies \eqref{eq:F-ellip} -- \eqref{eq:F-0}. One may also check  \eqref{eq:Fb-W2VMO} for $\bar F_B$, for the same reason shown in the proof of Lemma \ref{lemma:iter-msr}. Besides, $\e_B < 1$ because of \eqref{eq:e-r}. 

Summing up all the observations above, $\e_B$, $F_B$, $u_B^{\e_B}$, $\ell_B$, $\partial\Omega_B\cap B_2$, $\Phi_B$ and $f_B$ verify all the hypotheses for the first case of Lemma \ref{lemma:apprx-msr-bdry} (with $\rho = m^{\frac{kp}{n}} r_B$ and $c\delta$ in place of $\delta$ there). Therefore, we obtain that for any $s>0$, 
$$
|A_s^\e (u_B^{\e_B},\Omega_B)\cap \Omega_{B,\e_B}\cap B_1 | \leq c(\delta^\mu s^{-\mu} + s^{-q}),
$$
where $\Omega_{B,\e_B}$ denotes the set  $\{x\in \Omega_B : \dist(x,\partial\Omega_B) > \e_B\}$. Rephrasing the inequality in terms of $u^\e$, we obtain $|A_{c m^k s}(u^\e,\Omega\cap B_2) \cap \Omega_\e\cap B_{2r_B}(x_B)| \leq c(\delta^\mu s^{-\mu} + s^{-q}) r_B^n \leq \eta_0 |B|$. Evaluating this inequality at $s= c^{-1}m$ and using $B\subset B_{2r_B}(x_B)$, we arrive at $|A_{k+1}\cap B| \leq \eta_0|B|$, a contradiction. This finishes the proof for Case 1.

\vspace{0.2cm}
\noindent
{\it Case }2. 
$C_k = L_{cm^{k\alpha}}^\e (u^\e, \Omega\cap B_2)\cap \Omega_\e\cap B_{3/2}$, $D_k = \{ M(|D^2 \Phi|^q \chi_{B_2}) > \delta^q m^{kn(1-\alpha)q}\}$ and $(1-\alpha)n = p < n < q < \frac{n}{1-\alpha}$.
\vspace{0.2cm}

Let us remark that $C_k$ in Case 2 is the same with Case 1 by taking $p =(1-\alpha)n$, while $D_k$ in Case 2 replaces $n$ and $p$ in that of Case 1 with $\frac{n}{1-\alpha}$ and respectively $n$. Keeping this in mind, we follow the lines of the proof for Case 1. Then one may observe that under the new hypothesis in Case 2, \eqref{eq:D2Phi-Ln} is replaced by $\int_{B_{2r_B}(x_B)} |D^2 \Phi|^q\,dx \leq (4r_B)^n \delta^q m^{kn(1-\alpha)q} $, \eqref{eq:ue-l} remains the same, \eqref{eq:Dl} is replaced by $|D\ell|\leq 6m^{k\alpha}$ (here we also need $k\leq \frac{1}{1-\alpha}\log_m\frac{1}{\e}$, which is ensured from the statement of this lemma) and \eqref{eq:f-Lp} is replaced by $\int_{\Omega\cap B_{2r_B}(x_B)} |f|^p\,dx \leq (4r_B)^n (\delta m^k)^{n(1-\alpha)}$. Thus, with $\Phi_B$, $u_B^{\e_B}$, $\e_B$, $\ell_B$, $f_B$, $F_B$ as in \eqref{eq:PhiB}, \eqref{eq:uBeB-eB-lB-fB} and \eqref{eq:uBeB-pde}, one can verify that hypotheses for the second case of Lemma \ref{lemma:apprx-msr-bdry} is satisfied (with $\rho = m^{k(1-\alpha)}r_B$ and $c\delta$ in place of $\delta$ there). The rest of the proof is the same. Let us skip the detail in order to avoid redundant argument. 

\vspace{0.2cm}
\noindent
{\it Case }3. 
$\| D^2 \Phi\|_{L^p(B_2)} \leq \delta$, $C_k = \emptyset$, $D_k\{M(|D^2 \Phi|^n\chi_{B_2} )> \delta^n m^{kn}\}$ and $n < p< q$.
\vspace{0.2cm}

Unlike the first two case, we need to treat the last case differently. By the additional assumption $\|D^2 \Phi\|_{L^p(B_2)} \leq \delta$ and $p > n$, the Sobolev embedding theorem implies that $[D\Phi]_{C^{0,1-n/p}(B_2)} \leq c\delta$. As $\Phi$ being the boundary flattening map of $\partial\Omega\cap B_2$, it follows that $\partial\Omega\cap B_2$ is a $C^{1,1-n/p}$-hypersurface, whose norm is bounded by $c\delta$; since this implication is already rigorously justified in the proof of Lemma \ref{lemma:apprx-msr-bdry}, we shall not repeat it here. Since $\cP_- (D^2 u^\e) \leq f\leq \cP_+(D^2 u^\e)$ in $\Omega\cap B_2$ in the viscosity sense, and $u^\e = 0$ on $\partial\Omega\cap B_2$, one can find a linear polynomial $\ell_{x_B}$, according to \cite[Theorem 1.6]{LZ3} together with $\|u^\e \|_{L^\infty(\Omega\cap B_2)}\leq 1$ and $\|f\|_{L^p(\Omega\cap B_2)} \leq \delta$, that 
\begin{equation}\label{eq:ue-l-C1a-bdry}
|(u^\e - \ell_{x_B})(x)|\leq c|x-x_B|^{1+\alpha},
\end{equation}
for all $x\in\Omega\cap B_2$, where $c>1$ and $\alpha\in(0,1-\frac{n}{p})$ depend at most on $n$, $\lambda$, $\Lambda$ and $p$. 

Next, since $F(D^2 (u^\e - \ell_{x_B}) , \frac{\cdot}{\e}) = f$ in $B_{d_B}(\tilde x_B)$, with $d_B = \dist(\tilde x_B,\partial \Omega) = |x_B -\tilde x_B|$, in the viscosity sense, we may apply Lemma \ref{lemma:int-C1a} to $u^\e$ (with $\Omega = B_{d_B}(\tilde x_B)$ and $x_0 = \tilde x_B$ there), and deduce from \eqref{eq:ue-l-C1a-bdry}, as well as an obvious inequality $(I_{(1-\bar\alpha)p} (|f|^p \chi_{\Omega\cap B_2}) (\tilde x_B))^{\frac{1}{p}} \leq c\|f \|_{L^p(B_{d_B}(\tilde x_B))} \leq c\delta$, that 
\begin{equation}\label{eq:ue-l-C1a-int}
| (u^\e - \ell_{\tilde x_B} )(x) | \leq c( |x-\tilde x_B|^{1+\alpha} + \e^{1+\alpha} ),
\end{equation} 
for any $x\in B_{d_B}(\tilde x_B)$, for some other linear polynomial $\ell_{\tilde x_B}$. Since $d_B = |x_B - \tilde x_B| > \e$, we may compare this with \eqref{eq:ue-l} in $B_\e(\tilde x_B)\subset B_{d_B}(\tilde x_B)$ and utilize $m^k \e^2 \leq \e^{1+\alpha}$ (which is ensured by the choice $k \leq \frac{p}{n}\log_m\e$ and $\alpha\leq 1-\frac{n}{p}$) to deduce that $| (\ell - \ell_{\tilde x_B})(x)| \leq c\e^{1+\alpha}$ for all $x\in B_\e(\tilde x_B)$. Especially,
\begin{equation}\label{eq:Dl-re}
| D\ell - D\ell_{\tilde x_B} | \leq cm^{k\alpha}\e^\alpha.
\end{equation}

On the other hand, we may also derive from $\tilde x_B\not\in D_k$ and $|x_B - \tilde x_B| < 2r_B$ that
\begin{equation}\label{eq:D2Phi-Lp}
\int_{B_{2r_B}(x_B)} |D^2 \Phi|\,dx \leq 4^n \delta^p r_B^n m^{kp},
\end{equation}
Now we define $\Phi_B$, $u_B^{\e_B}$, $\e_B$, $f_B$ and $F_B$ as in \eqref{eq:PhiB}, \eqref{eq:uBeB-eB-lB-fB} and \eqref{eq:uBeB-pde}, but redefine $\ell_B$ by
$$
\ell_B (x) = \frac{(\ell - \ell_{\tilde x_B})(x_B + 2r_B x)}{cm^k r_B^2}. 
$$
By the obvious identity, $u^\e - \ell = u^\e - \ell_{\tilde x_B} - (\ell - \ell_{\tilde x_B})$, we have $u_B^{\e_B} = -\ell_B$ on $\partial\Omega_B\cap B_2$. Moreover, due to \eqref{eq:Dl-re} and \eqref{eq:e-r}, $|D\ell_B| \leq c m^{k(\alpha-1)}r_B^{\alpha-1}$. In addition, by \eqref{eq:D2Phi-Lp}, we have $\|D^2 \Phi_B\|_{L^p(B_2)} \leq 4\delta m^k r_B$. The other properties concerning $u_B^{\e_B}$, $f_B$, $F_B$ and $\Phi_B$ remain the same as in the proof of Case 1. Thus, $\e_B$, $u_B^{\e_B}$, $f_B$, $\ell_B$, $\Phi_B$ and $F_B$ altogether now verify the hypotheses of the last case of Lemma \ref{lemma:apprx-msr-bdry} (with $\rho = m^kr_B$ and $\alpha$ as above). The rest of the argument is the same with that of Case 1, whence it is left out to the reader.
\end{proof}

We are ready to prove the first assertion in Theorem \ref{theorem:bdry-W2p}. 

\begin{proof}[Proof of Theorem \ref{theorem:bdry-W2p}; the first part]
Fix any (finite) exponent $p > p_0$. Consider the case $0\in\partial\Omega$, $U = B_2$, $U' = B_1$, $\|u^\e \|_{L^\infty(\Omega\cap B_2)} \leq 1$, $\|f \|_{L^p(\Omega\cap B_2)} \leq \delta$, $g = 0$ on $B_2$, $\Phi(0) = 0$, $\osc_{B_2} D\Phi \leq \delta$ and $|\xi^t D\Phi(0) \xi - 1| \leq \delta$ for all $\xi\in\partial B_1$. Also assume that $\| D^2 \Phi\|_{L^n(B_2)} \leq \delta$ if $p < n$, $\|D^2 \Phi\|_{L^{n+\sigma}(B_2)} \leq \delta$ if $p = n$, and $\|D^2 \Phi\|_{L^p(B_2)} \leq \delta$ if $p > n$.

Throughout this proof, $c$ will denote a positive constant depending at most on $n$, $\lambda$, $\Lambda$, $\kappa$, $\psi$, $p$ and $\sigma$.

Choose $p' < p$ by $p' = \frac{p_0 + p}{2}$ if $p<n$, $p' = \frac{n}{n+\sigma}$ if $p = n$, with $\sigma$ as in the statement of the theorem, and $p ' = \frac{n + p}{2}$ if $p > n$. Let $\alpha_k$, $\beta_k$, $\gamma_k$ and $\delta_k$ be the measure of $A_k$, $B_k$, $C_k$ and respectively $D_k$ as in Lemma \ref{lemma:iter-msr-bdry} with $p'$ in place of $p$, and $\alpha = \frac{\sigma}{n+\sigma}$ there. Since $u^\e$, $f$, $F$ and $\Phi$ verify the hypotheses of Lemma \ref{lemma:iter-msr-bdry}, depending on the value of $p'$, we obtain, after iteration, that
$$
\alpha_k \leq \eta^k + \sum_{i=1}^k \eta^i (\beta_{k-i} + \gamma_{k-i} + \delta_{k-i}),
$$
with $\eta = c(\delta^{\gamma \mu} m^\mu + m^{-q})$. 

Fix any $q > p$ such that $q<n$ if $p < n$, $q \leq n+\sigma$ if $p = n$ and $ q = 2p$ if $p>n$. We may choose $m$ larger in Lemma \ref{lemma:iter-msr-bdry}, but still depending on the quantities specified in the statement of the lemma, such that $cm^{p-q} \leq \frac{1}{4}$. Then we take $\delta$ sufficiently small such that $c\delta^{\gamma\mu} m^\mu \leq\frac{1}{4}$, which ensures that $m^p \eta \leq\frac{1}{2}$. 

With such a choice of $m$ and $\delta$, one can derive, by following computations in the proof of Theorem \ref{theorem:int-W1p}, that 
\begin{equation}\label{eq:iter-msr-bdry}
\sum_{k=1}^{-c(p',n)\log_m \e} m^{k p}\alpha_k \leq c + c \sum_{k=1}^\infty m^{kp} (\beta_k + \gamma_k +\delta_k). 
\end{equation}
for some $c>0$, depending only on $n$, $\lambda$, $\Lambda$, $\kappa$, $\psi$, $p$ and $\sigma$ (only for the case $p = n$). Hence, it suffices to prove a uniform bound of the rightmost term in the above display.

By the strong $(\frac{p}{p'},\frac{p}{p'})$-type inequality for the maximal function and the assumption that $\|f \|_{L^p(\Omega\cap B_2)} \leq \delta$, one can immediately prove that 
$$
\sum_{k=1}^\infty m^{kp}\beta_k \leq c. 
$$

As for the summability of $m^{kp} \gamma_k$, we only need to take care of the case $p \leq n$, since for the other case, $p>n$, Lemma \ref{lemma:iter-msr-bdry} (iii) assumes $C_k = \emptyset$, i.e., $\gamma_k = 0$. For the case $p\leq n$, it follows from Theorem \ref{theorem:bdry-W1p} (along with $\|u^\e\|_{L^\infty(\Omega\cap B_2)} \leq 1$, $\|f\|_{L^p(\Omega\cap B_2)} \leq \delta$, $g = 0$ on $\partial\Omega\cap B_2$ and $\partial\Omega\cap B_2$ is a $C^1$-hypersurface whose Lipschitz norm is less than $c\delta$) that $G_{\Omega\cap B_2}^\e (u^\e )\in L^{np/(n-p')}(\Omega_\e\cap B_{3/2})$ (note $p' < n$ when $p\leq n$) and 
$$
\int_{\Omega_\e\cap B_2} (G_{\Omega\cap B_2}^\e (u^\e ))^{\frac{np}{n-p'}}\,dx \leq c. 
$$
Thus, writing by $\vp$ the function $G_{\Omega\cap B_2}^\e (u^\e)^{\frac{n}{n-p'}}$, the above display implies that $\int_{\Omega_\e\cap B_2} \vp^p\,dx \leq c$. By the relation between the function $G_E^\e(v)$ and the set $L^\e( v,E)$ (see Definition \ref{definition:large-grad-Hess}), $\{\vp > m^k\} = \{ G_{\Omega\cap B_2}^\e (u^\e) > m^{k(1-\frac{p'}{n})}\} = L_{ m^{k(1-p'/n)}}^\e (u^\e,\Omega\cap B_2)$, so 
$$
\sum_{k=1}^\infty m^{kp} \gamma_k = \sum_{k=1}^\infty m^{kp} |\{ \vp > m^k\}| \leq c. 
$$

Finally, let us verify the summability of $m^{kp}\delta_k$. As for the case $p>n$ (hence $p' > n$), it follows from the strong $(\frac{p}{p'},\frac{p}{p'})$-type inequality for the maximal function and the assumption that $\| D^2 \Phi\|_{L^p(B_2)} \leq \delta$. As for the case $p \leq n$, hence either $p < q < n$ or $p = n < q< n + \sigma$, we invoke strong $(\frac{N}{q},\frac{N}{q})$-type inequality, with $N = n $ if $p < n$ or $N = n+\sigma$ if $p = n$, for the maximal function. Then from the assumption $\| D^2 \Phi \|_{L^N(B_2)} \leq \delta$, we have $\int_{\R^n} M(|D^2 \Phi|^q\chi_{B_2})^{\frac{N}{q}}\,dx\leq c\delta$ and thus,  
$$
\sum_{k=1}^\infty m^{kp} \delta_k = \sum_{k=1} (m^{\frac{p q}{N}})^{\frac{k N}{q}} |\{ M(|D^2 \Phi|^q\chi_{B_2}) > \delta^q (m^{\frac{pq}{N}})^k\}| \leq c. 
$$

In all, we have proved that the rightmost term of \eqref{eq:iter-msr-bdry} is bounded by $c$, uniformly in $\e$. This finishes the proof, for the special case. By a standard rescaling argument, one may recover the case for general $\partial\Omega$, $U$, $U'$, $u^\e$, $f$ yet $g = 0$. Now for non-trivial $g\in W^{2,p}(U)$, we observe that $w^\e = u^\e - g$ satisfies the special case, with $f$ replaced by $f + c_0|D^2 g|$, with $c_0$ depending only on $n$, $\lambda$ and $\Lambda$; for more detail, see the proof of \cite[Theorem 4.5]{W}. 
\end{proof}

\begin{proof}[Proof of Theorem \ref{theorem:bdry-W2p}; the second part] 
Under the additional assumptions, i.e., \eqref{eq:F-C} and \eqref{eq:F-W2q}, we can obtain, in the proof of Lemma \ref{lemma:apprx-msr-bdry}, that  $|D^2h^\e| \in L^q(\Omega\cap B_1)$, or $|D^2\phi^\e| \in L^q(\Omega\cap B_1)$ depending on the value of $q$, uniformly in $\e$, where $h^\e$ and $\phi^\e$ are viscosity solutions to \eqref{eq:he-pde} and respectively \eqref{eq:phie-pde}. Thus \eqref{eq:Tehe-Lq} and respectively \eqref{eq:Hephie-Lq} hold with $A_s^\e$ replaced by $A_s$. As a result, Lemma \ref{lemma:apprx-msr-bdry} now holds with $A_s$ and $\Omega$ in place of $A_s^\e$ and respectively $\Omega_\e$. Moreover, the condition on $k \leq c(n,p)|\log_m \e|$ can also be removed. This readily implies the same improvement the statement of Lemma \ref{lemma:iter-msr-bdry}. We leave out the detail to the reader. 
\end{proof}



\begin{thebibliography}{99999999}

\bibitem[AKM19]{AKM}
S. Armstrong. T. Kuusi and J.-C. Mourrat,
{\it Quantitative Stochastic Homogenization and Large-Scale Regularity},
Vol. 352. Springer, 2019. 

\bibitem[AL89]{AL2}
M. Avellaneda and F.-H. Lin, 
{\it Compactness methods in the theory of homogenization II: equations in non-divergence form}, 
Comm. Pure Appl. Math. {\bf 42} (1989), 139--172.


\bibitem[AL91]{AL3}
M. Avellaneda and F.-H. Lin, 
{\it $L^p$-bounds on singular integrals in homogenization}, 
Comm. Pure Appl. Math. {\bf 44} (1991), 897--910.

\bibitem[BJ16]{BJ}
S.-S. Byun and Y. Jang
{\it Global $W^{1,p}$ estimates for elliptic systems in homogenization problems in Reifenberg domains},
Annali di Matematica {\bf 195} (2016), 2061--2075.


\bibitem[BW04]{BW}
S.-S. Byun and L. Wang,
{\it Elliptic equations with BMO coefficients in Reifenberg domains},
Comm. Pure Appl. Math. {\bf 57} (2004), 1283--1310.

\bibitem[CC95]{CC} 
X. Cabr\'e and L. A. Caffarelli, 
{\it Fully Nonlinear Elliptic Equations}, 
American Math. Society, Providence, RI, 1995.

\bibitem[Caf89]{Caff}
L. A. Caffarelli,
{\it Interior a priori estimates for solutions of fully non-linear equations},
Ann. Math. {\bf 130} (1989), 189--213.


\bibitem[CCKS96]{CCKS}
L. A. Caffarelli, M. G. Crandall, M. Kocan and A. Swiech,
{\it On viscosity solutions of fully nonlinear equations with measurable ingredients}, 
Comm. Pure Appl. Math. {\bf 49} (1996), 365--397.

\bibitem[CH03]{CH}
L. A. Caffarelli and Q. Huang,
{\it Estimates in the generalized Campanato-John-Nirenberg spaces for fully nonlinear elliptic equations},
Duke Math. J. {\bf 118} (2003), 1--17. 

\bibitem[CC98]{CP}
L. A. Caffarelli and I. Peral, 
{\it On $W^{1,p}$ estimates for elliptic equations in divergence form}, 
Comm. Pure Apple. Math. {\bf 51} (1998), 1--21.

\bibitem[CIL92]{CIL}
M. G. Crandall, H. Ishii and P.-L. Lions,
{\it User's guide to viscosity solutions of second order partial differential equations},
Bull. Amer. Math. Soc. {\bf 27} (1992), 1--67. 

\bibitem[CKS\'S96]{CKSS}
M.G. Crandall, M. Kocan, P. Soravia and A. Swiech, 
{\it On the equivalence of various weak notions of solutions of elliptic PDEs with measurable ingredients}, 
in: Progress in Elliptic and Parabolic Partial Differential Equations (Capri, 1994), in: Pitman Res. Notes Math. Ser., vol. 350,
Longman, Harlow, 1996, pp. 136--162.

\bibitem[DKM14]{DKM}
P. Daskalopoulos , T. Kuusi and G. Mingione,
{\it Borderline estimates for fully nonlinear elliptic equations},
Comm. Partial Differential Equations {\bf 39} (2014), 574--590.

\bibitem[Esc93]{Esc}
L. Escauriaza,
{\it $W^{2,n}$ a priori estimates for solutions to fully non-linear equations},
Indiana Univ. Math. J. {\bf 42} (1993), 413--423. 

\bibitem[Eva92]{E}
L. C. Evans, 
{\it Periodic homogenisation of certain fully nonlinear partial differential equations},
Proc. Roy. Soc. Edinburgh Sect. A {\bf 120} (1992), 245--265.

\bibitem[GT01]{GT}
D. Gilbarg and N. S. Trudinger, 
{\it Elliptic Partial Differential Equations of Second Order (reprint of 1998 ed.)},
Classics in Mathematics. Berlin: Springer 2001.

\bibitem[Hua02]{H2}
Q. Huang, 
{\it On the regularity of solutions to fully nonlinear elliptic equations via the Liouville property},
Proc. Amer. Math. Soc. {\bf 130} (2002), 1955--1959.

\bibitem[Hua19]{H}
Q. Huang,
{\it Regularity theory for $L^n$-viscosity solutions to fully nonlinear elliptic equations with asymptotical approximate convexity},
Ann. I. H. Poincar\'e {\bf 36} (2019), 1869--1902. 


\bibitem[KLS13]{KLS}
C. E. Kenig, F. Lin and Z. Shen,
{\it Homogenization of elliptic systems with Neumann boundary conditions},
J. Amer. Math. Soc. {\bf 26} (2013), 901--937. 

\bibitem[KL20]{KL}
S. Kim and K.-A. Lee,
{\it Uniform estimates in periodic Homogenization of fully nonlinear elliptic equations},
preprint available at arXiv:2011.08590

\bibitem[KM08]{KM}
J. Kristensen and C. Melcher, 
{\it Regularity in oscillatory nonlinear elliptic systems},
Math. Z. {\bf 260} (2008), 813--847.

\bibitem[LL86]{LL}
J. M. Lasry and P. L. Lions,
{\it A remark of regularization in Hilbert spaces},
Israel J. Math. {\bf 55} (1986). 

\bibitem[LZ15]{LZ}
D. Li and K. Zhang,
{\it $W^{2,p}$ interior estimates of fully nonlinear elliptic equations},
Bull. London Math. Soc. {\bf 47} (2015), 301--314. 

\bibitem[LZ18]{LZ2}
Y. Lian and K. Zhang,
{\it Boundary H\"older regularity on Reifenberg domains for fully nonlinear elliptic equations},
preprint available at arXiv:1812.11354

\bibitem[LZ20]{LZ3}
Y. Lian and K. Zhang,
{\it Boundary pointwise $C^{1,\alpha}$ and $C^{2,\alpha}$ regularity for fully nonlinear elliptic equations},
J. Differential Equations {\bf 269} (2020), 1172-1191.

\bibitem[LXZ21]{LXZ}
Y. Lian, W. Xu and K. Zhang,
{\it Boundary Lipschitz regularity and the Hopf lemma on Reifenberg domains for fully nonlinear elliptic equations},
Manuscripta Math. {\bf 166} (2021), 343--357. 

\bibitem[MS09]{MS}
C. Melcher and B. Schweizer,
{\it Direct approach to L p estimates in homogenization theory},
Annali di Matematica {\bf 188} (2009), 399--416.



\bibitem[SS14]{SS}
L. Silvestre and B. Sirakov,
{\it Boundary regularity for viscosity solutions of fully nonlinear elliptic equations},
Comm. Partial Differential Equations {\bf 39} (2014), 1694--1717. 


\bibitem[Sir10]{Sir}
B. Sirakov,
{\it Solvability of uniformly elliptic fully nonlinear PDE},
Arch. Rational Mech. Anal. {\bf 195} (2010), 579--607.

\bibitem[She17]{She}
Z. Shen
{\it Boundary estimates in elliptic homogenization},
Anal. PDE {\bf 10} (2017), 653--694.

\bibitem[{\'Swi}97]{S}
A. \'Swiech
{\it $W^{1,p}$-Interior estimates for solutions of fully nonlinear, uniformly elliptic equations},
Adv. Differential Equations {\bf 2} (1997), 1005--1027.


\bibitem[Tei14]{T}
E. V. Teixeira,
{\it Universal moduli of continuity for solutions to fully nonlinear elliptic equations},
Arch. Rational Mech. Anal. {\bf 211} (2014), 911--927. 

\bibitem[Win09]{W}
N. Winter,
{\it $W^{2,p}$ and $W^{1,p}$-estimates at the boundary for solutions of fully nonlinear, uniformly elliptic equations},
Z. Anal. Anwend. {\bf 28} (2009), 129--164. 

\end{thebibliography}
\end{document}